\pdfoutput=1
\RequirePackage{ifpdf}
\ifpdf 
\documentclass[pdftex]{sigma}
\else
\documentclass{sigma}
\fi

\numberwithin{equation}{section}

\newtheorem{Theorem}{Theorem}[section]
\newtheorem*{Theorem*}{Theorem}
\newtheorem{Corollary}[Theorem]{Corollary}
\newtheorem{Lemma}[Theorem]{Lemma}
\newtheorem{Proposition}[Theorem]{Proposition}
 { \theoremstyle{definition}
\newtheorem{Definition}[Theorem]{Definition}

\newtheorem{Remark}[Theorem]{Remark} }

\usepackage{mathtools}

\newcommand{\qbinom}{\genfrac{[}{]}{0pt}{}}

\DeclareMathOperator{\End}{End}

\DeclareMathOperator{\op}{op}

\newcommand{\fg}{\mathfrak{g}}
\newcommand{\wfg}{\widehat{\mathfrak{g}}}
\newcommand{\fn}{\mathfrak{n}}
\newcommand{\fh}{\mathfrak{h}}
\newcommand{\BN}{\mathbb{N}}
\newcommand{\BZ}{\mathbb{Z}}
\newcommand{\BC}{\mathbb{C}}
\newcommand{\gl}{\mathfrak{gl}}
\newcommand{\ssl}{\mathfrak{sl}}
\newcommand{\sso}{\mathfrak{so}}
\newcommand{\ssp}{\mathfrak{sp}}
\newcommand{\uu}{U_{r,s}(\fg)}
\newcommand{\UU}{U_{r,s}(\widehat{\fg})}
\newcommand{\Uu}{U'_{r,s}(\widehat{\fg})}
\newcommand{\Id}{\mathrm{Id}}
\newcommand{\slaws}{\text{standard Lyndon words}}

\newcommand{\II}{\mathcal{I}}

\newcommand\iso{ \vphantom{j^{X^2}}\smash{\overset{\sim}{\vphantom{\rule{0pt}{0.20em}}\smash{\longrightarrow}}} }

\begin{document}
\allowdisplaybreaks

\newcommand{\arXivNumber}{2407.01450}

\renewcommand{\PaperNumber}{064}

\FirstPageHeading

\ShortArticleName{Two-Parameter Quantum Groups and $R$-Matrices: Classical Types}

\ArticleName{Two-Parameter Quantum Groups and $\boldsymbol{R}$-Matrices:\\ Classical Types}

\Author{Ian MARTIN and Alexander TSYMBALIUK}
\AuthorNameForHeading{I.~Martin and A.~Tsymbaliuk}
\Address{Department of Mathematics, Purdue University, West Lafayette, IN, USA}
\Email{\href{mailto:mart2151@purdue.edu}{mart2151@purdue.edu}, \href{mailto:sashikts@gmail.com}{sashikts@gmail.com}}

\ArticleDates{Received January 05, 2025, in final form July 13, 2025; Published online July 31, 2025}

\Abstract{We construct finite $R$-matrices for the first fundamental representation $V$ of two-parameter quantum groups $U_{r,s}(\mathfrak{g})$ for classical $\mathfrak{g}$, both through the decomposition of $V\otimes V$ into irreducibles $U_{r,s}(\mathfrak{g})$-submodules as well as by evaluating the universal $R$-matrix. The latter is crucially based on the construction of dual PBW-type bases of $U^{\pm}_{r,s}(\mathfrak{g})$ consisting of the ordered products of quantum root vectors defined via $(r,s)$-bracketings and combinatorics of standard Lyndon words. We further derive explicit formulas for affine $R$-matrices, both through the Yang--Baxterization technique of [\textit{Internat. J.~Modern Phys.~A}~\textbf{6} (1991), 3735--3779] and as the unique intertwiner between the tensor product of~$V(u)$ and~$V(v)$, viewed as modules over two-parameter quantum affine algebras $U_{r,s}(\widehat{\mathfrak{g}})$ for classical~$\mathfrak{g}$. The latter generalizes the formulas of~[\textit{Comm. Math. Phys.} \textbf{102} (1986), 537--547] for one-parametric quantum affine algebras.}


\Classification{17B37; 16T25}

\section{Introduction}

\subsection{Summary}

Let $\fg$ be a Kac--Moody Lie algebra of finite type. Then, it admits a root space decomposition
\[
\fg=\fn^- \oplus \fh \oplus \fn^+ \qquad \text{with}\quad \fn^{\pm}=\oplus_{\alpha\in \Phi^+} \BC\cdot e_{\pm \alpha}
\]
 corresponding to a polarization of the root system $\Phi=\Phi^+ \sqcup (-\Phi^+)$. The elements $e_{\pm \alpha}$ are called \textit{root vectors}. Thus,
 \[
 U(\fg)=U(\fn^-) \otimes U(\fh) \otimes U(\fn^+)
 \]
 and the ordered products in $\{e_{\pm \alpha}\}_{\alpha\in \Phi^+}$ form a basis of $U\bigl(\fn^\pm\bigr)$ for any total order on $\Phi^\pm$. The root vectors can actually be normalized so that\footnote{Throughout the paper, we use $R^\times$ to denote the set of nonzero elements of any ring $R$.}
\begin{equation}\label{eqn:Lie-root-vectors}
 [e_{\alpha},e_{\beta}] = e_{\alpha}e_{\beta}-e_{\beta}e_{\alpha} \in \BZ^\times \cdot e_{\alpha+\beta}
 \qquad \mathrm{for\ all} \ \alpha, \beta \in \Phi^+ \ \mathrm{such\ that} \ \alpha+\beta \in \Phi^+.
\end{equation}
This inductively recovers all root vectors from the generators $\{e_i\}_{i\in I}$, corresponding to simple roots $\{\alpha_i\}_{i\in I}$. When $\fg$ is a Kac--Moody algebra of affine type, the root subspaces corresponding to imaginary roots are no longer one-dimensional. However, the theory of such algebras and their representations are well-understood due to their alternative realization as central extensions of the loop algebras $L\bar{\fg}$ for $\bar{\fg}$ of finite type:
\begin{equation}\label{eq:affLie-loop}
 0\to \BC \to \fg\to L\bar{\fg}\to 0.
\end{equation}

For any Kac--Moody algebra $\fg$, Drinfeld and Jimbo simultaneously introduced the quantum groups $U_q(\fg)$ which quantize the universal enveloping algebras of $\fg$. Similarly to $U(\fg)$, the quantum groups admit a triangular decomposition $U_q(\fg) = U_q(\fn^-) \otimes U_q(\fh) \otimes U_q(\fn^+)$. Furthermore, $U_q\bigl(\fn^\pm\bigr)$ admit PBW-type bases
\begin{equation}\label{eqn:PBW-quantum-intro}
 U_q\bigl(\fn^\pm\bigr) =
 \bigoplus_{\gamma_1 \geq \dots \geq \gamma_k \in \Phi^+} \BC(q) \cdot e_{\pm \gamma_1} \cdots e_{\pm \gamma_k}
\end{equation}
formed by the ordered products of \textit{$q$-deformed root vectors} $e_{\pm \alpha} \in U_q\bigl(\fn^\pm\bigr)$, defined via Lusztig's braid group action, which requires one to choose a reduced decomposition of the longest element~$w_0$ in the Weyl group $W$ of $\fg$. It is well-known \cite{P} that this choice precisely ensures that the order $\geq$ on $\Phi^+$ is convex, in the sense of Definition~\ref{def:convex}. Moreover, the $q$-deformed root vectors satisfy a $q$-analogue of relation~\eqref{eqn:Lie-root-vectors}, where $\alpha$, $\beta$ and $\alpha+\beta$ are any positive roots satisfying $\alpha<\alpha+\beta<\beta$ and the minimality property~\eqref{eq:min_pair}:
\begin{equation}\label{eqn:quantum-root-vectors}
 [e_{\alpha} , e_{\beta}]_q = e_{\alpha} e_{\beta} - q^{(\alpha, \beta)} e_{\beta} e_{\alpha}
 \in \BZ\bigl[q,q^{-1}\bigr]^\times \cdot e_{\alpha+\beta},
\end{equation}
where $(\cdot, \cdot)$ denotes the scalar product corresponding to the root system of type $\fg$. Therefore, the $q$-deformed root vectors can be defined (up to scalar multiple) as iterated $q$-commutators of the Drinfeld--Jimbo generators $e_i$, using the combinatorics of the root system and the chosen convex order on $\Phi^+$.

There is however a purely combinatorial approach to the construction of PBW-type bases of~$U_q\bigl(\fn^\pm\bigr)$, cf.\ \eqref{eqn:PBW-quantum-intro}, that goes back to the works of~\cite{K1,K2,L,Ro}. To this end, recall Lalonde--Ram's bijection~\cite{LR}:
\begin{equation}\label{eqn:LR-bij-intro}
 \ell \colon\ \Phi^+ \iso \{\text{standard Lyndon words in}\ I\}.
\end{equation}
We note that in the context of~\eqref{eqn:LR-bij-intro}, the notion of standard Lyndon words intrinsically depends on a fixed total order of the indexing set $I$ of simple roots. Furthermore,~\eqref{eqn:LR-bij-intro} gives rise to a~total order on $\Phi^+$ via:
\begin{displaymath}
 \alpha < \beta \quad \Longleftrightarrow \quad \ell(\alpha) < \ell(\beta) \ \text{lexicographically}.
\end{displaymath}
It was shown in~\cite{Ro} (see~\cite[Proposition~26]{L}) that this total order is convex, and hence can be applied to obtain root vectors $e_{\pm \alpha} \in U_q\bigl(\fn^\pm\bigr)$ for all $\alpha\in \Phi^+$, cf.\ \eqref{eqn:quantum-root-vectors}, thus eliminating Lusztig's braid group action.

When $\fg$ is of finite type and $q$ is not a root of unity, the representation theory of $U_q(\fg)$ is completely parallel to that of $\fg$. On the other hand, to develop the representation theory of $U_q(\fg)$ for affine $\fg$ one needs an alternative ``new Drinfeld'' realization \smash{$U^{\mathrm{Dr}}_q(\fg)$} from~\cite{D}, a $q$-analogue of~\eqref{eq:affLie-loop}. The isomorphism
\begin{equation}\label{eq:Dr-to-DJ-intro}
 \Psi\colon\ U^{\mathrm{Dr}}_q(\fg)\iso U_q(\fg) \qquad \text{for affine} \ \fg
\end{equation}
was constructed in~\cite{B} using an affine braid group action, while \smash{$\Psi^{-1}$} was stated in~\cite{D} using \textit{$q$-bracketings}.

One of the key features of quantum groups is that they are actually quasitriangular Hopf algebras. The corresponding universal $R$-matrices $\mathcal{R}\in U_q(\fg)\otimes U_q(\fg)$ (one needs to consider a~completion here) satisfy
\begin{equation}\label{eq:qYB-intro}
 \text{\textit{quantum Yang--Baxter equation:}}
 \quad \mathcal{R}_{12}\mathcal{R}_{13}\mathcal{R}_{23}=\mathcal{R}_{23}\mathcal{R}_{13}\mathcal{R}_{12}.
\end{equation}
In particular, for any two finite-dimensional $U_q(\fg)$-modules $V$, $W$ one obtains a $U_q(\fg)$-module intertwiner
\begin{equation}\label{eq:intertwiner-intro}
 \hat{R}_{VW}=(\rho_W\otimes \rho_V)(\mathcal{R})\circ \tau\colon\ V\otimes W \iso W\otimes V,
\end{equation}
where $\tau\colon V\otimes W\to W\otimes V$ is the flip map $v\otimes w\mapsto w\otimes v$, and $\rho_V\colon U_q(\fg)\to \End(V)$, $\rho_W\colon U_q(\fg)\to \End(W)$.

In fact, quantum groups first appeared in the \emph{quantum inverse scattering method}, the study of exactly solvable statistical models and quantum integrable systems arising through the quantum Yang--Baxter equation~\eqref{eq:qYB-intro}. In this context, one starts with a solution of~\eqref{eq:qYB-intro} or its version with a spectral parameter
\begin{equation}\label{eq:qYB-affine-intro}
 \mathcal{R}_{12}(x)\mathcal{R}_{13}(xy)\mathcal{R}_{23}(y)=
 \mathcal{R}_{23}(y)\mathcal{R}_{13}(xy)\mathcal{R}_{12}(x),
\end{equation}
and defines the algebra \smash{$U^{\mathrm{RTT}}_q(\fg)$} via the so-called \emph{RTT-relations}, see~\cite{FRT}. The isomorphisms
\begin{equation}\label{eq:RTT-to-DJ-intro}
 \Upsilon\colon\ U_q(\fg) \iso U^{\mathrm{RTT}}_q(\fg)
 \qquad \text{for finite type}\ \fg
\end{equation}
and
\begin{equation}\label{eq:RTT-to-Dr-intro}
 \Upsilon\colon\ U^{\mathrm{Dr}}_q(\fg) \iso U^{\mathrm{RTT}}_q(\fg)
 \qquad \mathrm{for\ affine\ type}\ \fg
\end{equation}
were first constructed in~\cite{DF} for types $A_n$, \smash{$A^{(1)}_n$} through the Gauss decomposition of the generating matrices. For other classical Lie algebras and their affinizations, such isomorphisms were first discovered in~\cite{HM} and were revised much more recently in~\cite{JLM2, JLM1}.

The theory of multiparameter quantum groups goes back to the early 90s, see, e.g.,~\mbox{\cite{AST,R,T}}. However, the current interest in the subject stems from the papers~\cite{BW2, BW3, BW1},
which study the two-parameter quantum group $U_{r,s}(\gl_n)$ and provide a further application to pointed finite-dimensional Hopf algebras. In~\cite{BW2}, they developed the theory of finite-dimensional representations in a complete analogy with the one-parameter case, computed the two-parameter $R$-matrix for the first fundamental $U_{r,s}(\gl_n)$-representation, and used it to establish the Schur--Weyl duality between $U_{r,s}(\gl_n)$ and a two-parameter Hecke algebra.

The above works of Benkart and Witherspoon stimulated an increased interest in the theory of two-parameter quantum groups. In particular, the definition and the basic structural results on $\uu$ for other classical simple Lie algebras $\fg$ were provided in~\cite{BGH1,BGH2}. Since then, multiple papers have treated such algebras case-by-case; we refer the reader to~\cite{HP} for a more uniform treatment and complete references.

The generalization of this theory from simple finite-dimensional Lie algebras to affine Lie algebras started with the work~\cite{HRZ} (which however had some gaps in the exposition, see~\cite{Ts} for an alternative treatment of the PBW results stated in~\cite{HRZ} without any proof). Subsequently, some attempts were made to provide a uniform Drinfeld--Jimbo presentation of such algebras, establishing the triangular decomposition and the Drinfeld double construction for them. More importantly, a new Drinfeld realization of these algebras $\UU$ was established on a case-by-case basis for $\fg$ being of type $A_n$ (see~\cite{HRZ}), types $D_n$ and $E_6$ (see~\cite{HZ2}), type $G_2$ (see~\cite{GHZ}), and type~$C_n$ (see~\cite{HZ1}). However, we note a caveat in this treatment: while a surjective homomorphism from the Drinfeld--Jimbo to the new Drinfeld realization is constructed similarly to $\Psi^{-1}$ of~\eqref{eq:Dr-to-DJ-intro}, there is no proper proof of its injectivity. The aforementioned new Drinfeld realization of $\UU$ was used to construct the vertex representations of $\UU$ in an analogy with the one-parameter case (cf.\ \cite{FJ}). Finally, the FRT-formalism for two-parameter quantum groups was carried out for~$U_{r,s}(\gl_n)$ and~$U_{r,s}\bigl(\widehat{\gl}_n\bigr)$ in~\cite{JL1,JL2}, establishing the two-parameter analogues of~\eqref{eq:RTT-to-DJ-intro} and~\eqref{eq:RTT-to-Dr-intro} for types $A_n$ and \smash{$A^{(1)}_n$}.

In this work (followed up by~\cite{MT2}), we develop the FRT-formalism for both~$\uu$,~$\UU$ when $\fg$ is an orthogonal or symplectic Lie algebra. The present note is mostly concerned with the derivation of finite and affine $R$-matrices (denoted by $R$ and $R(u)$), while in~\cite{MT2} we use~$R$,~$R(u)$ to construct~\eqref{eq:RTT-to-DJ-intro} and~\eqref{eq:RTT-to-Dr-intro}, naturally generalizing~\cite{JLM2, JLM1} to the two-parameter setup. The latter, in particular, immediately provides the new Drinfeld realization of $\UU$ (which seemed to be missing for $B_n$). Let us outline the key ingredients.

To derive the finite $R$-matrix $\hat{R}_{VV}$, we factorize the universal $R$-matrix into the ``local'' ones (one factor for each positive root of $\fg$) and then evaluate the result on the first fundamental representation $V$ of $\uu$. Due to the absence of Lusztig's braid group action on $\uu$ (noted first in~\cite{BGH1}), we use the aforementioned construction of orthogonal dual bases of positive and negative subalgebras of $\uu$ through standard Lyndon words and the technique of quantum shuffle algebras, which goes back to~\cite{K1,K2,L,Ro} in the one-parameter setup, to~\cite{CHW} in the super case, and finally to~\cite{BH,BKL,HW1,HW2} and our accompanying note~\cite{MT1} for the two-parameter case. We note, however, that once the explicit formula is obtained, one can directly check that it coincides with the $R$-matrix $\hat{R}_{VV}$ by verifying that it intertwines the $\uu$-action on $V\otimes V$ and acts by the same scalars on the highest weight vectors of $V\otimes V$ as $\hat{R}_{VV}$.

To derive the affine $R$-matrix $\hat{R}_{VV}(u/v)$, we use the \emph{Yang--Baxterization} technique of~\cite{GWX}. Since $\hat{R}_{VV}$ has three distinct eigenvalues (in contrast to the case $\fg=\ssl_n$ when it has only two distinct eigenvalues), there are 12 possible resulting operators, and in each case one needs to check some extra conditions to guarantee that they satisfy the Yang--Baxter equation. Instead, once the explicit formula for the correctly chosen one is derived, it is straightforward to check that it intertwines the $\UU$-actions on the tensor products of evaluation modules $V(u)\otimes V(v)\to V(v)\otimes V(u)$.
The results of~\cite{Jim} then guarantee that the space of all such intertwiners is at most one-dimensional, and therefore the operator $\hat{R}(u/v)$ constructed through the Yang--Baxterization coincides with $\hat{R}_{VV}(u/v)$, thus producing a solution~of~\eqref{eq:qYB-affine-intro} by further composing with the flip map $\tau$.

While we were finishing the present note, closely related preprints~\cite{HXZ1,HJZ,HXZ2} appeared on arXiv. Though we communicated our results to one of those authors back in February 2024, it~is~a~pity they decided not to consolidate our papers to be posted simultaneously. Partially due to this flaw, the present note is separated from a more straightforward part~\cite{MT2} that will be posted later.

\subsection{Outline}

The structure of the present paper is the following:
\begin{itemize}\itemsep=0pt
\item[$\bullet$]
In Section~\ref{ssec:general-def}, we recall  two-parameter quantum groups $U_{r,s}(\fg)$ for simple finite-dimen\-sional Lie algebras~$\fg$, see Definition~\ref{def:general_2param}, and summarize their basic properties (including the pairing of Proposition~\ref{prop:pairing_2param}).

\item[$\bullet$]
In Section~\ref{sec:column_repr}, we explicitly construct the first fundamental representations of $U_{r,s}(\fg)$ for classical $\fg$, see Propositions~\ref{prp:A_rep}--\ref{prp:D_rep}. We further decompose the tensor product $V\otimes V$ into irreducible $U_{r,s}(\fg)$-submodules, see Proposition~\ref{prop:struct}. The proof of the latter is derived, through a reduction to the Lie algebra limit, by providing explicitly the corresponding highest weight vectors.

\item[$\bullet$]
In Section~\ref{sec:R-matrices}, we evaluate the universal intertwiner $\hat{R}_{VV}$ from Theorem~\ref{thm:universal-R} on the first fundamental $\uu$-representations from Section~\ref{sec:column_repr} for $\fg=\sso_{2n+1},\ssp_{2n},\sso_{2n}$, see Theorems~\ref{thm:B_RMatrix}--\ref{thm:D_RMatrix}. This generalizes the corresponding formula of Theorem~\ref{thm:A_RMatrix} discovered first in~\cite{BW2}. Our proof is slightly different though, as we only match the eigenvalues of the three highest weight vectors in $V\otimes V$ featured in Proposition~\ref{prop:struct} (see Lemmas~\ref{lem:B_eigen}--\ref{lem:D_eigen}), and then verify the intertwining property with the action of $f_i$'s (see~Lemma~\ref{lem:B_Intertwining}).

\item[$\bullet$]
In Section~\ref{sec:R-factorized}, we provide an alternative proof of the formulas for $\hat{R}_{VV}\hspace{-0.67pt}$ from Theorems~\mbox{\ref{thm:A_RMatrix}--\ref{thm:D_RMatrix}} by showing that they arise as the product of ``local'' operators parametrized by the positive roots of $\Phi$. Since Lusztig's braid group action has no analogue for $\uu$, we instead use a combinatorial approach to the construction of orthogonal PBW bases, see Theorem~\ref{thm:PBW-general} (which constitutes the main result of~\cite{MT1}). To this end, we construct the (quantum) root vectors iteratively by using the combinatorics of standard Lyndon words, recalled in Sections~\ref{ssec:Lyndon} and~\ref{ssec:PBW theorem and bases}. The factorization formula that results from these considerations is stated in Theorem~\ref{cor:Theta-factorization}. We conclude the section with a case-by-case treatment of each classical series, providing proofs of the formulas \eqref{eq:A_RMatrix}, \eqref{eq:B_RMatrix}, \eqref{eq:C_RMatrix}, \eqref{eq:D_RMatrix} that are more conceptual than those presented in Section~\ref{sec:R-matrices}.

\item[$\bullet$]
In Section~\ref{sec:affine-R}, we introduce the two-parameter quantum affine algebras $\UU$, see Definition~\ref{def:2_param_aff}, which is in agreement with~\cite{HRZ, HZ2,HZ1} for $\fg$ of types $A_n$, $C_n$, $D_n$. We also introduce their counterparts $\Uu$ without \emph{degree generators} and extend the first fundamental $\uu$-representations $\rho$ from Propositions~\ref{prp:A_rep}--\ref{prp:D_rep} to evaluation $\Uu$-representations in Propositions~\ref{prop:A_aff_rep}--\ref{prop:D_aff_rep}. The latter ones are upgraded to $\UU$-modules in Proposition~\ref{prop:affine-actions-withD}. The main results of this section are Theorems~\ref{thm:B_AffRMatrix}--\ref{thm:D_AffRMatrix}, which evaluate the universal intertwiner of $\UU$ on the tensor product of two such representations for $\fg=\sso_{2n+1},\ssp_{2n},\sso_{2n}$. This generalizes the corresponding formula~\eqref{eq:A_AffRMatrix} of Theorem~\ref{thm:A_AffRMatrix}, first discovered in~\cite{JL2}. According to~\cite{Jim}, composing $\hat{R}(z)$ of~(\ref{eq:A_AffRMatrix})--(\ref{eq:D_AffRMatrix}) with a flip map $\tau$ produces solutions of the Yang--Baxter relation with a spectral parameter~\eqref{eq:qYB-affine-intro}. While the proofs are straightforward, the origin of these formulas (whose one-parameter counterparts were discovered in~\cite{Jim}) is postponed till Section~\ref{sec:Baxterization}.

\item[$\bullet$]
In Section~\ref{sec:Baxterization}, we derive the formulas~(\ref{eq:A_AffRMatrix})--(\ref{eq:D_AffRMatrix}) through the \emph{Yang--Baxterization} technique of~\cite{GWX}.
\end{itemize}

\section{Notations and definitions}

Throughout the paper, we will work over an algebraically closed field $\mathbb{K} \supset \mathbb{C}(r,s)$, where $r$ and~$s$ are indeterminates.\footnote{While most of the formulas
are valid for arbitrary parameters $r,s$, we need to assume that $r/s$ is not a root of unity for Section~\ref{sec:R-factorized}, and we further prefer to keep
them algebraically independent for the proof of Proposition~\ref{lem:dim} and construction~\eqref{eq:sigma-antiautom}.}

\subsection{Two-parameter finite quantum groups}\label{ssec:general-def}

Let $E$ be a Euclidean space with a symmetric bilinear form $(\cdot ,\cdot)$, and $\Phi \subset E$ be an indecomposable reduced root system with an ordered set of simple roots $\Pi = \{\alpha_{1},\dots ,\alpha_{n}\}$. Let $\mathfrak{g}$ be the complex simple Lie algebra corresponding to this root system.
Let $C = (c_{ij})_{i,j = 1}^{n}$ be the Cartan matrix of $\fg$, explicitly given by \smash{$c_{ij} = \frac{2(\alpha_{i},\alpha_{j})}{(\alpha_{i},\alpha_{i})}$}, and let \smash{$d_{i} = \frac{1}{2}(\alpha_{i},\alpha_{i})$} where $(\cdot,\cdot)$ is normalized so that the short roots have square length $2$.
We denote the root and weight lattices of $\fg$ by~$Q$ and~$P$,~respectively,
\begin{displaymath}
 \bigoplus_{i=1}^n \BZ \alpha_i=Q\subset P=\bigoplus_{i=1}^n \BZ \varpi_i
 \qquad \mathrm{with}\ (\alpha_i,\varpi_j)=d_i\delta_{ij}.
\end{displaymath}

Having fixed above the order on the set of simple roots $\Pi$, we consider the (modified) Ringel bilinear form $\langle \cdot ,\cdot \rangle$ on $Q$, such that (unless $\{i,j\}=\{n-1,n\}$ in type $D_n$) we have
\begin{equation}\label{eq:ringel-1}
 \langle \alpha_{i},\alpha_{j}\rangle =
 \begin{cases}
 d_{i}c_{ij} & \text{if}\ i < j, \\
 d_{i} & \text{if}\ i = j ,\\
 0 & \text{if}\ i > j,
 \end{cases}
\end{equation}
while in the remaining case of $D_n$-type system, we set
\begin{gather}
 \langle \alpha_{n-1},\alpha_{n} \rangle =
 \langle \varepsilon_{n-1}-\varepsilon_n,\varepsilon_{n-1}+\varepsilon_n\rangle = -1, \nonumber\\
 \langle \alpha_{n},\alpha_{n-1}\rangle =
 \langle \varepsilon_{n-1}+\varepsilon_n,\varepsilon_{n-1}-\varepsilon_n \rangle = 1.\label{eq:ringel-2}
\end{gather}
We note that $(\mu,\nu) = \langle \mu,\nu \rangle + \langle \nu,\mu \rangle$ for any $\mu,\nu \in Q$.

\begin{Remark}
The modification~\eqref{eq:ringel-2} is made to ensure that the algebra $U_{r,s}(\sso_{2n})$ defined below matches its original definition in~\cite{BGH1}.
\end{Remark}

We shall also need two-parameter analogues of $q$-integers and $q$-factorials (cf.\ \cite[equation~(2.2)]{BW2}):
\begin{gather*}
 [m]_{r,s} = \frac{r^{m} - s^{m}}{r - s} = r^{m-1} + r^{m-2}s + \dots + rs^{m-2} + s^{m-1}
 \qquad \mathrm{for\ all} \ m\in \BN,
\\ 
 [m]_{r,s}! = [m]_{r,s}[m-1]_{r,s} \cdots [1]_{r,s} \qquad \mathrm{for} \ m > 0,
 \quad [0]_{r,s}!=1,
\end{gather*}
and most importantly two-parameter analogues of Gaussian binomial coefficients:
\begin{displaymath}
 \qbinom{m}{k}_{r,s} = \frac{[m]_{r,s}!}{[m - k]_{r,s}![k]_{r,s}!}
 \qquad \mathrm{for\ all} \ 0\leq k\leq m.
\end{displaymath}
Finally, in analogy with the one-parameter case (see~\cite[Section~4.2]{J}), we define
\begin{align}
 & r_{\gamma} = r^{(\gamma,\gamma)/2},\ \qquad s_{\gamma} = s^{(\gamma,\gamma)/2}
 \qquad \mathrm{for\ all} \ \gamma\in \Phi, \nonumber\\
 & r_{i} = r_{\alpha_i} = r^{d_{i}},\qquad s_{i} = s_{\alpha_i} = s^{d_{i}}
 \qquad \mathrm{for\ all} \ 1\leq i\leq n.\label{eq:rs_i}
\end{align}

We are ready now to introduce the main actor of this paper, the two-parameter quantum group of~$\fg$, following~\cite{HP}. While this definition is uniform for all types, we will make it more explicit for classical Lie algebras $\fg$ in Section~\ref{ssec:classical-2param}.

\begin{Definition}\label{def:general_2param}
The \textit{two-parameter quantum group} $U_{r,s}(\fg)$ of a simple Lie algebra $\fg$ is the associative $\mathbb{K}$-algebra generated by $\bigl\{e_{i},f_{i},\omega_{i}^{\pm 1},\bigl(\omega_{i}'\bigr)^{\pm 1}\bigr\}_{i=1}^{n}$ with the following defining relations (for all $1\leq i,j\leq n$):
\begin{gather}\label{eq:R1}
 [\omega_i,\omega_j]=\bigl[\omega_i,\omega'_j\bigr]=\bigl[\omega'_i,\omega'_j\bigr]=0, \qquad
 \omega_{i}^{\pm 1}\omega_{i}^{\mp 1} = 1 = \bigl(\omega_{i}'\bigr)^{\pm 1}\bigl(\omega_{i}'\bigr)^{\mp 1},
\\ \label{eq:R2}
 \omega_{i}e_{j} = r^{\langle \alpha_{j},\alpha_{i}\rangle}s^{-\langle \alpha_{i},\alpha_{j}\rangle}e_{j}\omega_{i}, \qquad
 \omega_{i}f_{j} = r^{-\langle \alpha_{j},\alpha_{i}\rangle}s^{\langle \alpha_{i},\alpha_{j}\rangle}f_{j}\omega_{i},
\\ \label{eq:R3}
 \omega_{i}'e_{j} = r^{-\langle \alpha_{i},\alpha_{j}\rangle}s^{\langle \alpha_{j},\alpha_{i}\rangle}e_{j}\omega_{i}', \qquad
 \omega_{i}'f_{j} = r^{\langle \alpha_{i},\alpha_{j}\rangle}s^{-\langle \alpha_{j},\alpha_{i}\rangle}f_{j}\omega_{i}',
\\ \label{eq:R4}
 e_{i}f_{j} - f_{j}e_{i} = \delta_{ij}\frac{\omega_i-\omega'_i}{r_{i} - s_{i}},
\end{gather}
and quantum $(r,s)$-Serre relations
\begin{align}
 &\sum_{k = 0}^{1 - c_{ij}}(-1)^k\qbinom{1 - c_{ij}}{k}_{r_{i},s_{i}}(r_{i}s_{i})^{\frac{1}{2}k(k-1)}(rs)^{k\langle \alpha_{j},\alpha_{i}\rangle}e_{i}^{1 - c_{ij} - k}e_{j}e_{i}^{k}=0
 \qquad \forall i\ne j, \nonumber\\
 &\sum_{k = 0}^{1 - c_{ij}}(-1)^k\qbinom{1 - c_{ij}}{k}_{r_{i},s_{i}}(r_{i}s_{i})^{\frac{1}{2}k(k-1)}(rs)^{k\langle \alpha_{j},\alpha_{i}\rangle} f_{i}^{k}f_{j}f_{i}^{1 - c_{ij} - k}=0
 \qquad \forall i\ne j.\label{eq:R5}
\end{align}
\end{Definition}

\begin{Remark}
We note that this definition does depend on the choice of an order of $\Pi$. We shall make a standard choice for classical $\fg$ in the end of this section, thus matching with the rest of the literature~\cite{BW1,BGH1}.
\end{Remark}

We note that the algebra $U_{r,s}(\mathfrak{g})$ is $Q$-graded via
\begin{displaymath}
 \deg(e_{i})=\alpha_i,\quad \deg(f_i)=-\alpha_i, \quad
 \deg(\omega_{i})=0, \quad \deg\bigl(\omega_{i}'\bigr)=0 \qquad \mathrm{for\ all} \ 1\leq i\leq n.
\end{displaymath}
For $\mu\in Q$, let $U_{r,s}(\fg)_{\mu}$ \big(or simply $(U_{r,s})_{\mu}$\big) denote the degree $\mu$ component of $U_{r,s}(\fg)$ under this $Q$-grading.

Analogously to the one-parameter case (cf.\ \cite[Section~4.11]{J}), there is a Hopf algebra structure on $U_{r,s}(\fg)$, where the coproduct $\Delta$, counit $\epsilon$, and antipode $S$ are defined on generators by the following formulas:
\begin{alignat*}{4}
 &\Delta\bigl(\omega_{i}^{\pm 1}\bigr) = \omega_{i}^{\pm 1} \otimes \omega_{i}^{\pm 1}, \qquad&&
 \epsilon\bigl(\omega_{i}^{\pm 1}\bigr) = 1, \qquad&&
 S\bigl(\omega_{i}^{\pm 1}\bigr) = \omega_{i}^{\mp 1},& \\
 &\Delta\bigl(\bigl(\omega_{i}'\bigr)^{\pm 1}\bigr) = \bigl(\omega_{i}'\bigr)^{\pm 1} \otimes \bigl(\omega_{i}'\bigr)^{\pm 1},\qquad&&
 \epsilon\bigl(\bigl(\omega_{i}'\bigr)^{\pm 1}\bigr) = 1,\qquad&&
 S\bigl(\bigl(\omega_{i}'\bigr)^{\pm 1}\bigr) = \bigl(\omega_{i}'\bigr)^{\mp 1},& \\
 &\Delta(e_{i}) = e_{i} \otimes 1 + \omega_{i} \otimes e_{i},\qquad&&
 \epsilon(e_{i}) = 0,\qquad&&
 S(e_{i}) = -\omega_{i}^{-1}e_{i} ,&\\
 &\Delta(f_{i}) = 1 \otimes f_{i} + f_{i} \otimes \omega_{i}',\qquad&&
 \epsilon(f_{i}) = 0,\qquad&&
 S(f_{i}) = -f_{i}\bigl(\omega_{i}'\bigr)^{-1},&
\end{alignat*}

\begin{Remark}
The simplest way to see how the definition above generalizes the usual Drinfeld--Jimbo one-parametric quantum groups $U_q(\fg)$ (cf.\ \cite[Section~4]{J}) is to work in the numeric setup. To this end, let $r,s\in \BC\setminus \{0\}$ with $r^{2}\neq s^{2}$ and define $U_{r,s}(\fg)$ as in Definition~\ref{def:general_2param}, but now viewed as an algebra over $\mathbb{K}=\BC$. Then, for any $q\in \mathbb{C}$ with $q^{4} \neq 1$, there is a natural Hopf algebra epimorphism
\begin{gather*}
 \pi\colon \ U_{q,q^{-1}}(\fg)\twoheadrightarrow U_q(\fg) \\
\hphantom{\pi\colon} \ \mathrm{given\ by} \
 e_{i} \mapsto E_{i},\ f_{i} \mapsto F_{i},\
 \omega_{i} \mapsto K_{i},\ \omega_{i}' \mapsto K_{i}^{-1} \ \mathrm{for\ all} \ 1\leq i\leq n.
\end{gather*}
Moreover, the kernel of $\pi$ is the two-sided ideal $\mathfrak{I}$ generated by \smash{$\bigl\{\omega'_i-\omega_{i}^{-1}\bigr\}_{i=1}^{n}$}. Thus, we get
\[
 \mathrm{Hopf\ algebra\ isomorphism:}\quad U_{q,q^{-1}}(\fg)/\mathfrak{I} \iso U_{q}(\fg).
\]
\end{Remark}

Let us also define several subalgebras of $U_{r,s}(\fg)$:
\begin{itemize}\itemsep=0pt
\item the ``positive'' subalgebra $U_{r,s}^{+}(\fg)$, generated by $\{e_{i}\}_{i=1}^{n}$,

\item the ``negative'' subalgebra $U_{r,s}^{-}(\fg)$, generated by $\{f_{i}\}_{i=1}^{n}$,

\item the ``Cartan'' subalgebra $U_{r,s}^{0}(\fg)$, generated by
\smash{$\bigl\{\omega_{i}^{\pm 1},(\omega'_{i})^{\pm 1}\bigr\}_{i=1}^{n}$},

\item the ``non-negative subalgebra'' $U_{r,s}^{\ge}(\fg)$, generated by
\smash{$\bigl\{e_{i},\omega_{i}^{\pm 1}\bigr\}_{i=1}^{n}$},

\item the ``non-positive subalgebra'' $U_{r,s}^{\le}(\fg)$, generated by
\smash{$\bigl\{f_{i},(\omega'_{i})^{\pm 1}\bigr\}_{i=1}^{n}$}.

\end{itemize}
When $\fg$ is clear from the context, we will use $U_{r,s}$ instead of $U_{r,s}(\fg)$, and similarly for the above subalgebras.

For any $\mu=\sum_{i=1}^{n} k_{i}\alpha_{i}\in Q$, we define $\omega_\mu,\omega'_\mu\in U_{r,s}^{0}(\fg)$ via
\[
 \omega_{\mu} = \omega_{1}^{k_{1}}\omega_{2}^{k_{2}} \cdots \omega_{n}^{k_{n}}, \qquad
 \omega'_\mu = \bigl(\omega_{1}'\bigr)^{k_{1}}\bigl(\omega_{2}'\bigr)^{k_{2}} \cdots \bigl(\omega_{n}'\bigr)^{k_{n}}.
\]

\subsection{Hopf pairing}

One of the basic structural properties of $U_{r,s}(\fg)$ is that it may be realized as a Drinfel'd double of its subalgebras $U^{\le}_{r,s}(\fg)$ and $U^{\ge}_{r,s}(\fg)$ with respect to the Hopf algebra pairing. This was established case-by-case in the literature, and we shall rather just refer to~\cite{HP}:

\begin{Proposition}\label{prop:pairing_2param}
There exists a unique bilinear pairing
\begin{equation}\label{eq:Hopf-parity}
 (\cdot,\cdot)\colon\ U_{r,s}^{\le}(\fg) \times U_{r,s}^{\ge}(\fg) \to \mathbb{K}
\end{equation}
satisfying the following structural properties:
{\samepage\begin{gather}
 \bigl(yy',x\bigr) = \bigl(y \otimes y',\Delta(x)\bigr), \qquad \bigl(y, xx'\bigr) = \bigl(\Delta(y),x' \otimes x\bigr)\nonumber\\
 \forall x,x' \in U_{r,s}^{\ge}(\fg),\quad y,y' \in U_{r,s}^{\le}(\fg),\label{eq:Hopf-properties}
\end{gather}
as well as being given on the generators by}
\begin{align}
 & (f_{i},\omega_{j}) = 0, \qquad \bigl(\omega_{i}', e_{i}\bigr) = 0, \qquad
 (f_{i},e_{j}) = \delta_{ij}\frac{1}{s_{i} - r_{i}} \qquad \mathrm{for\ all} \ 1\leq i,j\leq n , \nonumber\\
 & \bigl(\omega_{\lambda}',\omega_{\mu}\bigr) = r^{\langle \lambda,\mu\rangle}s^{-\langle \mu,\lambda \rangle}
 \qquad \mathrm{for\ all} \ \lambda,\mu\in Q.\label{eq:generators-parity}
\end{align}
\end{Proposition}

We can also formally extend this pairing to the weight lattice $P$ as follows:
\begin{gather}
 \bigl(\omega_{\lambda}',\omega_{\mu}\bigr) = \prod_{i,j = 1}^{n}\bigl(\omega_{i}',\omega_{j}\bigr)^{\lambda_{i}\mu_{j}} \quad \ \
\mathrm{for\ any\ weights} \
 \lambda = \sum_{i = 1}^{n} \lambda_{i}\alpha_{i}\in P,\ \mu = \sum_{i = 1}^{n} \mu_{i}\alpha_{i} \in P.\label{eq:weight-pairing}
\end{gather}
Although $\lambda_{i}$, $\mu_{i}$ may not be integers, the expression above still makes sense because $\mathbb{K}$ is algebraically closed.

\begin{Remark}
As mentioned in the beginning of this subsection, the above pairing allows for the realization of any two-parameter quantum group $U_{r,s}(\fg)$ as a Drinfel’d double of its Hopf subalgebras $U^{\leq}_{r,s}(\fg)$, $U^{\geq}_{r,s}(\fg)$ with respect to the pairing $(\cdot,\cdot)$ of~\eqref{eq:Hopf-parity}.
\end{Remark}

Let us list several basic properties of this pairing that will be needed later:
\begin{itemize}\itemsep=0pt
\item
First, if $x\in U_{r,s}(\fg)_{\mu}$ and $\nu\in Q$, then we have
\begin{displaymath}
 \omega_{\nu}x\omega_{\nu}^{-1} = \bigl(\omega_{\mu}',\omega_{\nu}\bigr)x, \qquad
 \omega_{\nu}'x(\omega_{\nu}')^{-1} = \bigl(\omega_{\nu}',\omega_{\mu}\bigr)^{-1}x.
\end{displaymath}

\item
Second, the pairing $(\cdot,\cdot)$ is of homogeneous degree zero, i.e.,
\begin{equation}\label{eq:pairing-orthogonal}
 (y,x)=0 \qquad \mathrm{for} \
 x\in U_{r,s}^{\ge}(\fg)_{\mu},\ y\in U_{r,s}^{\le}(\fg)_{-\nu}
 \ \mathrm{with} \ \mu\ne \nu.
\end{equation}

\item
Third, similarly to the one-parameter case (cf.\ \cite[Sections~6.14 and 6.15]{J}), we have
\begin{gather}\label{eq:coprod-1}
 \Delta(x) \in x\otimes 1 + \bigoplus_{0<\nu<\mu} U_{r,s}^{+}(\fg)_{\mu - \nu}\omega_{\nu} \otimes U_{r,s}^{+}(\fg)_{\nu} + \omega_{\mu} \otimes x,
\\ \label{eq:coprod-2}
 \Delta(y) \in y \otimes \omega_{\mu}' + \bigoplus_{0<\nu<\mu} U_{r,s}^{-}(\fg)_{-\nu} \otimes U_{r,s}^{-}(\fg)_{-(\mu - \nu)}\omega_{\nu}' + 1 \otimes y
\end{gather}
for any $x\in U_{r,s}^{+}(\fg)_{\mu}$ and $y \in U_{r,s}^{-}(\fg)_{-\mu}$. Here, we use the
standard order $\leq$ on $Q$:
\[
 \nu\leq \mu \Longleftrightarrow
 \mu-\nu=\sum_{i} k_i\alpha_i \qquad \mathrm{with} \ k_i\in \BZ_{\geq 0}.
\]
Then, combining the properties~\eqref{eq:pairing-orthogonal}--\eqref{eq:coprod-2} with the defining
properties~\eqref{eq:Hopf-properties} and~\eqref{eq:generators-parity}, we obtain
\begin{align*}
 & \bigl(\omega_{\mu}'y,x\bigr) = \bigl(\omega_{\mu}',\omega_{\nu}\bigr)(y,x), \qquad
 (y\omega_{\mu}',x) = (y,x), \\
 & (y,\omega_{\nu}x) = (\omega_{\mu}',\omega_{\nu})(y,x), \qquad
 (y,x\omega_{\nu}) = (y,x)
\end{align*}
for any $x\in U_{r,s}^{+}(\fg)_{\nu}$ and $y \in U_{r,s}^{-}(\fg)_{-\mu}$.
\end{itemize}

\subsection{Classical types}\label{ssec:classical-2param}

Since in this paper we are only interested in the classical Lie algebras $\fg$, it will be helpful to have more explicit formulas for the bilinear form, avoiding the use of the form $\langle \cdot,\cdot \rangle$ on $Q$ defined in~\eqref{eq:ringel-1} and~\eqref{eq:ringel-2}. To this end, let us first recall the explicit realization of the classical root systems as well as specify the choice of simple roots for them:
\begin{itemize}\itemsep=0pt
\item
\emph{$A_{n}$-type} (corresponding to $\fg\simeq \ssl_{n+1}$).
Let $\{\varepsilon_{i}\}_{i=1}^{n+1}$ be an orthonormal basis of $\mathbb{R}^{n+1}$. Then, we have
\begin{align*}
 & \Phi_{A_{n}} = \{\varepsilon_{i} - \varepsilon_{j}\mid 1\leq i\ne j\leq n+1\}
 \subset \mathbb{R}^{n+1}, \qquad
 \Pi_{A_{n}} = \{\alpha_{i} = \varepsilon_{i} - \varepsilon_{i + 1}\}_{i=1}^{n}.
\end{align*}

\item
\emph{$B_{n}$-type} (corresponding to $\fg\simeq \sso_{2n+1}$).
Let $\{\varepsilon_{i}\}_{i=1}^{n}$ be an orthogonal basis of $\mathbb{R}^{n}$ with $(\varepsilon_{i},\varepsilon_{i})=2$ for all $i$. Then, we have
\begin{align*}
 & \Phi_{B_{n}} = \{\pm \varepsilon_{i} \pm \varepsilon_{j}\mid 1 \le i < j \le n\} \cup
 \{\pm \varepsilon_{i}\mid 1\leq i\leq n\} \subset \mathbb{R}^{n}, \\
 & \Pi_{B_{n}} = \{\alpha_{i} = \varepsilon_{i} - \varepsilon_{i + 1}\}_{i=1}^{n-1}
 \cup \{\alpha_{n} = \varepsilon_{n}\}.
\end{align*}

\item
\emph{$C_{n}$-type} (corresponding to $\fg\simeq \ssp_{2n}$).
Let $\{\varepsilon_{i}\}_{i=1}^{n}$ be an orthonormal basis of $\mathbb{R}^{n}$. Then, we have
\begin{align*}
 & \Phi_{C_{n}} = \{\pm \varepsilon_{i} \pm \varepsilon_{j}\mid 1 \le i < j \le n\}
 \cup \{\pm 2\varepsilon_{i}\mid 1\leq i\leq n\} \subset \mathbb{R}^{n}, \\
 & \Pi_{C_{n}} = \{\alpha_{i} = \varepsilon_{i} - \varepsilon_{i + 1}\}_{i=1}^{n-1}
 \cup \{\alpha_{n} = 2\varepsilon_{n}\}.
\end{align*}

\item
\emph{$D_{n}$-type} (corresponding to $\fg\simeq \sso_{2n}$).
Let $\{\varepsilon_{i}\}_{i=1}^{n}$ be an orthonormal basis of $\mathbb{R}^{n}$. Then, we have
\begin{align}
 & \Phi_{D_{n}} = \{\pm \varepsilon_{i} \pm \varepsilon_{j}\mid 1 \le i < j \le n\}
 \subset \mathbb{R}^{n}, \nonumber\\
 & \Pi_{D_{n}} = \{\alpha_{i} = \varepsilon_{i} - \varepsilon_{i + 1}\}_{i=1}^{n-1}
 \cup \{\alpha_{n} = \varepsilon_{n-1} + \varepsilon_{n}\}.\label{eq:D-system}
\end{align}
\end{itemize}

\begin{Remark}
We note that $(\varepsilon_i,\varepsilon_i)=2$ in type $B_n$ in agreement with our scaling of the form~$(\cdot,\cdot)$.
\end{Remark}

Then we have the following respective formulas for the pairing of Cartan elements
(we note that while the second formula follows from the first one in each case, it will be convenient to use both later on):
\begin{itemize}\itemsep=0pt
\item
\emph{$A_{n}$-type}
\begin{align*}
 &\bigl(\omega_{\lambda}',\omega_{i}\bigr)= r^{( \varepsilon_{i},\lambda )}s^{( \varepsilon_{i +1},\lambda)}, \qquad
 \bigl(\omega_{i}',\omega_{\lambda}\bigr)= r^{-(\varepsilon_{i + 1},\lambda )}s^{-( \varepsilon_{i},\lambda)}.
\end{align*}

\item
\emph{$B_{n}$-type}
\begin{align}
 &\bigl(\omega_{\lambda}',\omega_{i}\bigr)=
 \begin{cases}
 r^{( \varepsilon_{i},\lambda )}s^{( \varepsilon_{i + 1},\lambda )}
 & \mathrm{if}\ 1\leq i<n ,\\
 r^{( \varepsilon_{n},\lambda ) }(rs)^{-\lambda_{n}}
 & \mathrm{if}\ i = n,
 \end{cases} \nonumber\\
 &\bigl(\omega_{i}',\omega_{\lambda}\bigr)=
 \begin{cases}
 r^{-( \varepsilon_{i + 1},\lambda)}s^{-( \varepsilon_{i},\lambda)}
 & \mathrm{if}\ 1\leq i<n ,\\
 s^{-( \varepsilon_{n},\lambda ) }(rs)^{\lambda_{n}}
 & \mathrm{if}\ i = n.
 \end{cases} \label{eq:B-pairing}
\end{align}

\item
\emph{$C_{n}$-type}
\begin{align}
 &\bigl(\omega_{\lambda}',\omega_{i}\bigr)=
 \begin{cases}
 r^{( \varepsilon_{i},\lambda) }s^{( \varepsilon_{i + 1},\lambda )}
 & \mathrm{if}\ 1\leq i<n ,\\
 r^{2( \varepsilon_{n},\lambda ) }(rs)^{-2\lambda_{n}}
 & \mathrm{if}\ i=n,
 \end{cases} \nonumber\\
 &\bigl(\omega_{i}',\omega_{\lambda}\bigr)=
 \begin{cases}
 r^{-( \varepsilon_{i + 1},\lambda )}s^{-( \varepsilon_{i},\lambda)}
 & \mathrm{if}\ 1\leq i<n, \\
 s^{-2( \varepsilon_{n},\lambda )}(rs)^{2\lambda_{n}}
 & \mathrm{if}\ i = n.
 \end{cases} \label{eq:C-pairing}
\end{align}

\item
\emph{$D_{n}$-type}
\begin{align*}
 &\bigl(\omega_{\lambda}',\omega_{i}\bigr)=
 \begin{cases}
 r^{( \varepsilon_{i},\lambda )}s^{( \varepsilon_{i + 1},\lambda )}
 & \mathrm{if}\ 1\leq i<n, \\
 r^{( \varepsilon_{n-1},\lambda)}s^{-(\varepsilon_{n},\lambda)}(rs)^{-2\lambda_{n - 1}}
 & \mathrm{if}\ i = n,
 \end{cases} \\
 &\bigl(\omega_{i}',\omega_{\lambda}\bigr)=
 \begin{cases}
 r^{-( \varepsilon_{i + 1},\lambda )}s^{-(\varepsilon_{i},\lambda )}
 & \mathrm{if}\ 1\leq i<n ,\\
 r^{( \varepsilon_{n},\lambda )}s^{-(\varepsilon_{n-1},\lambda )}(rs)^{2\lambda_{n-1}}
 & \mathrm{if}\ i = n.
\end{cases}
\end{align*}
\end{itemize}

\section{Column-vector representations}\label{sec:column_repr}

Let $N=n+1$ for $A_n$-type, $N=2n+1$ for $B_n$-type, and $N=2n$ for $C_n$- and $D_n$-types. Let $V = \mathbb{K}^{N}$ with standard basis vectors $\{v_i\}_{i=1}^N$. In this section, we construct an action of $\uu$ on~$V$ for any classical $\fg$ by specifying explicitly $\rho\colon \uu \to \End(V)$, and further decom\-pose~${V\otimes V}$ into irreducible submodules.

\subsection{First fundamental representations}

While the $A$-type representation goes back to~\cite{BW2}, we decided to include it since it serves as a~prototype for the other classical types, new in the literature.

\begin{Proposition}[{\cite[Section~1]{BW2}}]\label{prp:A_rep}
The following defines a representation $\rho\colon  U_{r,s}(\ssl_{n+1})\to \End(V)$ with
\begin{align*}
 & \rho(e_{i}) = E_{i,i+1} , \qquad \rho(f_{i}) = E_{i+1,i} , \\
 & \rho(\omega_{i}) =
 rE_{ii} + sE_{i+1,i+1} + \sum_{1\leq j\leq n+1}^{j\ne i,i+1} E_{jj} , \qquad
 \rho\bigl(\omega'_{i}\bigr) =
 sE_{ii} + rE_{i+1,i+1} + \sum_{1\leq j\leq n+1}^{j \neq i,i + 1} E_{jj} .
\end{align*}
\end{Proposition}

In what follows, we shall use the involution $'$ on the indexing set $\{1, \dots, N\}$ defined via
\[
 i':= N + 1 - i \qquad \mathrm{for\ all} \ 1\leq i\leq N.
\]

\begin{Proposition}[type $B_n$]\label{prp:B_rep}
The following defines a representation $\rho\colon U_{r,s}(\sso_{2n + 1})\to \End(V)$ with
\begin{gather*}
 \rho(e_{i}) = E_{i,i + 1} - E_{(i + 1)',i'},\qquad 1\leq i\leq n , \\
 \rho(f_{i}) =
 \begin{cases}
 E_{i + 1, i} - (rs)^{-2}E_{i', (i + 1)'} & \mathrm{if}\ 1\leq i<n, \\
 \bigl(r^{-1} + s^{-1}\bigr)E_{n + 1,n} - \bigl(r^{-1} + s^{-1}\bigr)E_{n', n + 1} & \mathrm{if}\ i = n,
 \end{cases} \\
 \rho(\omega_{i}) =
 \begin{cases}
 \displaystyle r^{2}E_{ii} + s^{2}E_{i + 1, i + 1} + r^{-2}E_{i'i'} + s^{-2}E_{(i + 1)', (i + 1)'} \\
 \displaystyle \qquad{} + \sum_{1\leq j\leq n}^{j \neq i,i + 1} \bigl(E_{jj}+E_{j'j'}\bigr)
 + E_{n+1,n+1}
 & \mathrm{if}\ 1\leq i<n, \\
 \displaystyle{rs^{-1}E_{nn} + E_{n + 1,n+1} + r^{-1}sE_{n'n'} + \sum_{j = 1}^{n-1}\bigl(r^{-1}s^{-1}E_{jj} + rsE_{j'j'}\bigr)} & \mathrm{if}\ i = n,
 \end{cases} \\
 \rho\bigl(\omega_{i}'\bigr) =
 \begin{cases}
 \displaystyle s^{2}E_{ii} + r^{2}E_{i + 1, i + 1} + s^{-2}E_{i'i'} + r^{-2}E_{(i + 1)', (i + 1)'} \\
 \displaystyle\qquad{}+ \sum_{1\leq j\leq n}^{j \neq i,i + 1} \bigl(E_{jj} + E_{j'j'}\bigr)
 + E_{n+1,n+1}
 & \mathrm{if}\ 1\leq i<n ,\\
 \displaystyle{r^{-1}sE_{nn} + E_{n + 1,n+1} + rs^{-1}E_{n'n'} + \sum_{j = 1}^{n-1}\bigl(r^{-1}s^{-1}E_{jj} + rsE_{j'j'}\bigr)} & \mathrm{if}\ i = n.
\end{cases}
\end{gather*}
\end{Proposition}

\begin{proof}
The proof is straightforward, as we just need to verify that the above linear operators satisfy the defining relations~\eqref{eq:R1}--\eqref{eq:R5}. The relation~\eqref{eq:R1} is obvious since all the operators~$\rho(\omega_i)$,~$\rho(\omega'_j)$ act diagonally in the basis $\{v_k\}_{k=1}^N$. To check the first equality of~\eqref{eq:R2}, we note that both sides act trivially on $v_k$ unless $k\in \{j+1,j'\}$. In the latter case, one needs to compare the ratios of eigenvalues of $\rho(\omega_i)$ on
$v_{j}$ and $v_{j+1}$, or $v_{(j+1)'}$ and $v_{j'}$, to the pairing $(\omega'_j,\omega_i)$. The other three relations of~\eqref{eq:R2} and \eqref{eq:R3} are verified analogously. The relation $[\rho(e_i),\rho(f_j)]=0$ for $i\ne j$ is obvious as both $\rho(e_i)\rho(f_j)$ and $\rho(f_j)\rho(e_i)$ then act trivially on $V$. On the other hand, for $1\leq i<n$, the commutator $[\rho(e_i),\rho(f_i)]$ acts diagonally in the basis $\{v_k\}_{k=1}^N$, with nonzero eigenvalues $1$, $-1$, $r^{-2}s^{-2}$, $-r^{-2}s^{-2}$ only for $k=i,i+1,(i+1)',i'$, respectively, which exactly coincide with the eigenvalues of \smash{$\frac{\rho(\omega_i)-\rho(\omega'_i)}{r^2-s^2}$}. Likewise, the commutator $[\rho(e_n),\rho(f_n)]$ acts diagonally in the basis $v_k$ with nonzero eigenvalues $r^{-1}+s^{-1}$, $-r^{-1}-s^{-1}$ only for $k=n,n'$, respectively, which exactly coincide with the eigenvalues of \smash{$\frac{\rho(\omega_n)-\rho(\omega'_n)}{r-s}$}. Finally, the Serre relations~\eqref{eq:R5} hold as each summand acts trivially on all~$v_k$.\looseness=1
\end{proof}

\begin{Proposition}[type $C_n$]\label{prp:C_rep}
The following defines a representation $\rho\colon  U_{r,s}(\ssp_{2n})\to \End(V)$ with
\begin{gather*}
 \rho(e_{i}) =
 \begin{cases}
 E_{i,i + 1} - E_{(i + 1)',i'} & \mathrm{if}\ 1\leq i<n ,\\
 E_{nn'} & \mathrm{if}\ i=n,
 \end{cases} \\
 \rho(f_{i}) =
 \begin{cases}
 E_{i+ 1,i} - (rs)^{-1}E_{i',(i + 1)'} & \mathrm{if}\ 1\leq i<n ,\\
 (rs)^{-1}E_{n'n} & \mathrm{if}\ i=n ,
 \end{cases} \\
 \rho(\omega_{i}) =
 \begin{cases}
 \displaystyle rE_{ii} + sE_{i + 1, i+1} + r^{-1}E_{i'i'} + s^{-1}E_{(i +1)',(i+1)'} \\
 \displaystyle\qquad{} + \sum_{1\leq j\leq n}^{j \neq i,i + 1} \bigl(E_{jj} + E_{j'j'}\bigr)
 & \mathrm{if}\ 1\leq i<n ,\\
 \displaystyle{rs^{-1}E_{nn} + r^{-1}sE_{n'n'} + \sum_{j = 1}^{n-1}\bigl(r^{-1}s^{-1}E_{jj} + rsE_{j'j'}\bigr)}
 & \mathrm{if}\ i = n,
 \end{cases} \\
 \rho(\omega'_{i}) =
 \begin{cases}
 \displaystyle sE_{ii} + rE_{i + 1, i+1} + s^{-1}E_{i'i'} + r^{-1}E_{(i +1)',(i+1)'} \\
 \displaystyle \qquad{} + \sum_{1\leq j\leq n}^{j \neq i,i + 1} \bigl(E_{jj} + E_{j'j'}\bigr)
 & \mathrm{if}\ 1\leq i<n, \\
 \displaystyle{r^{-1}sE_{nn} + rs^{-1}E_{n'n'} + \sum_{j = 1}^{n-1}\bigl(r^{-1}s^{-1}E_{jj} + rsE_{j'j'}\bigr)}
 & \mathrm{if}\ i = n.
\end{cases}
\end{gather*}
\end{Proposition}

\begin{proof}
The proof is analogous to that of Proposition~\ref{prp:B_rep}; we leave details to the reader.
\end{proof}

\begin{Proposition}[type $D_n$]\label{prp:D_rep}
For $n \ge 2$, the following defines a representation $\rho\colon  U_{r,s}(\sso_{2n}) \to \End(V)$ with
\begin{gather*}
 \rho(e_{i}) =
 \begin{cases}
 E_{i,i+1} - E_{(i + 1)',i'} & \mathrm{if}\ 1\leq i<n, \\
 (rs)^{-1}E_{n - 1, n'} - E_{n,(n-1)'} & \mathrm{if}\ i=n, \\
 \end{cases} \\
 \rho(f_{i}) =
 \begin{cases}
 E_{i + 1,i} - (rs)^{-1}E_{i', (i + 1)'} & \mathrm{if}\ 1\leq i<n, \\
 E_{n', n-1} - E_{(n-1)',n} & \mathrm{if}\ i=n ,\\
 \end{cases} \\
 \rho(\omega_{i}) =
 \begin{cases}
 \displaystyle rE_{ii} + sE_{i + 1,i + 1} + r^{-1}E_{i'i'} + s^{-1}E_{(i + 1)', (i+1)'} \\
 \displaystyle\qquad{}+\sum_{1\leq j\leq n}^{j \neq i,i + 1} \bigl(E_{jj} + E_{j'j'}\bigr)
 & \mathrm{if}\ 1\leq i<n ,\\
 \displaystyle s^{-1}E_{n-1,n-1} + rE_{nn} + sE_{(n-1)',(n-1)'} + r^{-1}E_{n'n'} \\
 \displaystyle\qquad{}+ \sum_{j = 1}^{n-2}\bigl(r^{-1}s^{-1}E_{jj} + rsE_{j'j'}\bigr)
 & \mathrm{if}\ i = n ,\\
 \end{cases} \\
 \rho(\omega'_{i}) =
 \begin{cases}
 \displaystyle sE_{ii} + rE_{i + 1,i + 1} + s^{-1}E_{i'i'} + r^{-1}E_{(i + 1)',(i + 1)'} \\
 \displaystyle\qquad{}+\sum_{1\leq j\leq n}^{j \neq i,i + 1} \bigl(E_{jj} + E_{j'j'}\bigr)
 & \mathrm{if}\ 1\leq i<n ,\\
 \displaystyle r^{-1}E_{n-1,n-1} + sE_{nn} + rE_{(n-1)',(n-1)'} + s^{-1}E_{n'n'} \\
 \displaystyle\qquad{}+ \sum_{j = 1}^{n-2}\bigl(r^{-1}s^{-1}E_{jj} + rsE_{j'j'}\bigr)
 & \mathrm{if}\ i = n. \\
\end{cases}	 	
\end{gather*}
\end{Proposition}

\begin{proof}
The proof is analogous to that of Proposition~\ref{prp:B_rep}; we leave details to the reader.
\end{proof}

The classification of finite-dimensional $U_{r,s}(\fg)$-modules for classical $\fg$ is completely parallel to that of one-parameter quantum groups (cf.\ \cite[Section~5]{J}). For $A$-type, this has been carried out in~\cite{BW2}, while for $BCD$-types this constitutes the major result of~\cite{BGH2}. Let us recall only the notions relevant to the rest of this section.
A vector $v$ in a $U_{r,s}(\mathfrak{g})$-module $V$ is said to have \textit{weight} $\lambda \in P$ if
\begin{displaymath}
 \omega_{i}v = (\omega_{\lambda}',\omega_{i})v \qquad \mathrm{and} \qquad
 \omega_{i}'v = (\omega_{i}',\omega_{\lambda})^{-1}v \qquad \mathrm{for\ all} \ 1\leq i\leq n.
\end{displaymath}
Let $V[\lambda]$ denote the subspace of all weight $\lambda$ vectors in $V$.
The following result is straightforward.

\begin{Lemma}\quad
\begin{enumerate}\itemsep=0pt
\item[$(a)$] For the vector representation $V$ from Propositions~$\ref{prp:B_rep}$ in type $B_n$, we have
\begin{gather*}
 V = V[0] \oplus \bigoplus_{i=1}^{n}(V[\varepsilon_{i}] \oplus V[-\varepsilon_{i}])\\
\mathrm{with} \
 V[0] = \mathbb{K} v_{n+1}, \ V[\varepsilon_{i}] = \mathbb{K}v_{i},\ V[-\varepsilon_{i}] = \mathbb{K}v_{i'}
 \ \text{for} \  1\leq i\leq n.
\end{gather*}

\item[$(b)$] For the vector representation $V$ from Proposition~$\ref{prp:C_rep}$ in type $C_n$, we have
\[
 V = \bigoplus_{i=1}^{n} (V[\varepsilon_{i}] \oplus V[-\varepsilon_{i}])
 \qquad \mathrm{with} \
 V[\varepsilon_{i}] = \mathbb{K}v_{i},\ V[-\varepsilon_{i}] = \mathbb{K}v_{i'} \ \mathrm{for} \ 1\leq i\leq n .
\]

\item[$(c)$] For the vector representation $V$ from Proposition~$\ref{prp:D_rep}$ in type $D_n$, we have
\[
 V = \bigoplus_{i=1}^{n} (V[\varepsilon_{i}] \oplus V[-\varepsilon_{i}])
 \qquad \mathrm{with} \
 V[\varepsilon_{i}] = \mathbb{K}v_{i},\ V[-\varepsilon_{i}] = \mathbb{K}v_{i'} \ \mathrm{for} \ 1\leq i\leq n .
\]
\end{enumerate}
\end{Lemma}

We note that $\varepsilon_1=\varpi_1$ is the first fundamental weight of $\fg$. Since the vector $v_1$ is clearly annihilated by all $\rho(e_i)$ and $\rho\bigl(U^0_{r,s}(\fg)\bigr)$ acts diagonally on the basis $\{v_k\}_{k=1}^N$ with distinct joint eigenvalues, we conclude the following corollary.

\begin{Corollary}
For $\fg$ being one of the Lie algebras $\sso_{2n+1}$, $\ssp_{2n}$, $\sso_{2n}$, the $U_{r,s}(\fg)$-module $V$ constructed respectively in Propositions~$\ref{prp:B_rep}$--$\ref{prp:D_rep}$ is isomorphic to $L(\varpi_1)$, the first fundamental $\uu$-representation.
\end{Corollary}

\subsection{Decomposition of the tensor square}

In the rest of this section, we shall study the decomposition of the tensor square $V\otimes V$ into irreducible $U_{r,s}(\fg)$-submodules. In type $A_n$, this has been carried out in~\cite[Proposition~5.3]{BW2}.

\begin{Proposition}[{\cite{BW2}}]
The $U_{r,s}(\ssl_{n+1})$-module $V\otimes V$ decomposes into the direct sum of two irreducible modules $V\otimes V\simeq L(2\varepsilon_1)\oplus L(\varepsilon_1+\varepsilon_2)$ with the highest weight vectors $v_1\otimes v_1$ and~${v_1\otimes v_2 - r v_2\otimes v_1}$.
\end{Proposition}

In contrast, we shall show in Proposition~\ref{prop:struct} that $V\otimes V$ actually decomposes into the direct sum of three irreducible $U_{r,s}(\fg)$-submodules in the remaining classical types $B_n$, $C_n$, $D_n$ (with $n\geq 2$ for type~$C_n$ and $n\geq 3$ for type~$D_n$). To this end, we start by establishing the following preliminary result (which is of independent interest).

\begin{Proposition}\label{lem:dim}
If $\lambda \in P^+$ is a dominant integral weight, then the irreducible $U_{r,s}(\fg)$-module $L(\lambda)$ of the highest weight $\lambda$ has the same dimension as the corresponding irreducible $\fg$-module.
\end{Proposition}

\begin{proof}
The idea of the proof is to specialize $r$ to $q$ and $s$ to $q^{-1}$, and then appeal to the analogous result in the one-parameter case. The argument presented below closely follows that of~\cite[Sections~5.12--5.15]{J}.

Let us first set up some notation. Let \smash{$M(\lambda) = U_{r,s}(\mathfrak{g})\otimes_{U^{\geq}_{r,s}(\mathfrak{g})} \mathbb{K}$} be the Verma $\uu$-module of highest weight $\lambda\in P^+$, and let $V = L(\lambda)$ be its finite-dimensional irreducible quotient, with highest weight vector denoted by $v_\lambda$ (cf.\ \cite{BW2,BGH2}). Let $\mathbb{F}=\BC(r)$ and $A = \mathbb{F}\bigl[s,s^{-1}\bigr]$. For any sequence $\II = (i_{1},\dots ,i_{k})\in \{1, 2, \dots, n\}^k$, set $f_{\II} = f_{i_{1}} \cdots f_{i_{k}}$, so that $V$ is spanned by all the vectors $f_{\II}v_{\lambda}$. Since $V$ is finite-dimensional and $Q$-graded, we note that $f_{\II}v_{\lambda} = 0$ for all but finitely many sequences $\II$. Therefore the $A$-submodule
\[
 V_{A} = \sum_{\II} A f_{\II}v_{\lambda}
\]
is finitely generated, and is clearly also torsion-free. Thus, $V_{A}$ is a free $A$-module since $A$ is a PID. Furthermore, one can easily check that the natural map $\mathbb{K} \otimes_{A} V_{A} \to V$ is an isomorphism (since~$\mathbb{K}$ contains the field of fractions of $A$), so that $\mathrm{rk}_A{V_A}=\dim_{\mathbb{K}}{V}$.
On the other hand, we have an isomorphism $\mathbb{F} \iso \mathbb{C}(q)$ given by $r \mapsto q$, which may be further extended to a~homomorphism~${A \to \mathbb{C}(q)}$ via $s \mapsto q^{-1}$, making $\mathbb{C}(q)$ into an $A$-module. Consider the module $\overline{V} = \mathbb{C}(q) \otimes_{A} V_{A}$, which has an obvious $\mathbb{C}(q)$-vector space structure. Note that any basis of~$V_{A}$ over $A$ yields a basis of $\overline{V}$ over $\mathbb{C}(q)$. Thus, it suffices to make $\overline{V}$ into a $U_{q}(\mathfrak{g})$-module in such a~way that $\overline{V} \simeq L(\lambda)$ as $U_{q}(\mathfrak{g})$-modules, where $U_{q}(\mathfrak{g})$ is the Drinfeld--Jimbo quantum group over~$\mathbb{C}(q)$ (as the irreducible highest weight modules over $\fg$ and $U_q(\fg)$ have the same dimension, cf.~\cite[Theorem~5.15]{J}).

Let us first show that $V_{A}$ is actually stable under the action of all generators $f_{i}$, $e_{i}$, $\omega_{i}^{\pm 1}$, \smash{$\bigl(\omega_{i}'\bigr)^{\pm 1}$}. For $f_{i}$, this is obvious. Moreover, since all the pairings \smash{$(\omega_{\mu}', \omega_{i})^{\pm 1}$}, \smash{$\bigl(\omega_{i},\omega_{\mu}'\bigr)^{\pm 1}$} belong to~$A$, the fact that $V_{A}$ is stable under $\omega_{i}^{\pm 1}$, \smash{$\bigl(\omega_{i}'\bigr)^{\pm 1}$} follows from the fact that $V_{A}$ has a weight space decomposition. Finally, to prove stability of $V_A$ under $e_{i}$, we first note that if $v \in V_{\mu}$, then
\[
 \frac{\omega_{i} - \omega_{i}'}{r_i - s_i} v =
 \frac{r^{\langle \mu,\alpha_{i}\rangle}s^{-\langle \alpha_{i},\mu \rangle} - r^{-\langle \alpha_{i},\mu\rangle}s^{\langle \mu,\alpha_{i}\rangle}}{r_i - s_i} v =
 (rs)^{-\langle \alpha_{i},\mu\rangle}\frac{r^{(\alpha_{i},\mu)} - s^{(\alpha_{i},\mu)}}{r_i - s_i} v,
\]
which belongs to $V_{A}$ since $(\alpha_{i},\mu)$ is divisible by $d_i$. Evoking $e_{j}v_{\lambda} = 0$ for all $j$, we thus get
\[
 e_{j}f_{i_1}\cdots f_{i_k} v_{\lambda} =
 \sum_{1\leq s\leq k}^{i_s=j} f_{i_{1}} \cdots f_{i_{s-1}} \frac{\omega_{i_{s}} - \omega_{i_{s}}'}{r_i - s_i}
 f_{i_{s + 1}} \cdots f_{i_{k}} v_{\lambda},
\]
which is in $V_{A}$ by what we have already proved.

It remains to note that the factors $r^{\langle \alpha_{j},\alpha_{i}\rangle}s^{-\langle \alpha_{i},\alpha_{j}\rangle}$ and $r^{-\langle \alpha_{i},\alpha_{j}\rangle}s^{\langle \alpha_{j},\alpha_{i}\rangle}$ specialize to $q^{(\alpha_i,\alpha_j)}$ and $q^{-(\alpha_i,\alpha_j)}$, respectively, under the specialization $r\mapsto q$, $s\mapsto q^{-1}$. Therefore, the operators $e_{i}$, $f_{i}$, $\omega_{i}$, $\omega_{i}'$ on $\overline{V}$ satisfy the same relations as the generators $E_{i}$, $F_{i}$, $K_{i}$, $K_{i}^{-1}$ of $U_{q}(\fg)$, so that $\overline{V}$ is a $U_{q}(\mathfrak{g})$-module. Since $\overline{V}$ is also a finite-dimensional highest weight module with the highest weight $\lambda$, we have $\overline{V} \simeq L(\lambda)$ as $U_{q}(\fg)$-modules.
\end{proof}

\begin{Proposition}\label{prop:struct}
For $\fg$ being one of the Lie algebras $\sso_{2n + 1}, \ssp_{2n}, \sso_{2n}$ $($where $n \geq 2$ if $\fg=\mathfrak{sp}_{2n}$ and $n\geq 3$ if $\fg=\mathfrak{so}_{2n})$ and $V$ being the corresponding $U_{r,s}(\fg)$-representation from Propositions~$\ref{prp:B_rep}$--$\ref{prp:D_rep}$, we~have the following decomposition into irreducibles:
\begin{equation}\label{eq:irr-decomp}
 V \otimes V \simeq L(2\varepsilon_{1})\oplus L(\varepsilon_{1} + \varepsilon_{2}) \oplus L(0).
\end{equation}
\end{Proposition}

\begin{proof}
We shall only present complete details for $\mathfrak{g} = \mathfrak{so}_{2n + 1}$, the other cases being analogous.

\textit{Type $B_n$}.
Let us first show that the following are highest weight vectors for the $U_{r,s}(\sso_{2n+1})$-action on~${V \otimes V}$:
\begin{gather}
 w_{1} = v_{1} \otimes v_{1}, \nonumber\\
 w_{2} = v_{1} \otimes v_{2} - \bigl(\omega_{\varepsilon_{1}}',\omega_{1}\bigr)v_{2} \otimes v_{1} =
 \begin{cases}
 v_{1} \otimes v_{2} - rs^{-1}v_{2} \otimes v_{1} & \text{if}\ n = 1, \\
 v_{1} \otimes v_{2} - r^2 v_{2} \otimes v_{1} & \text{if}\ n > 1 ,
 \end{cases} \nonumber\\
 w_{3} = \sum_{i=1}^{n} r^{2(i - 1)}v_{i} \otimes v_{i'} + r^{2n - 1}s^{-1}v_{n + 1}\otimes v_{n+1} +
 \sum_{i = 1}^{n}r^{2n - 1}s^{2(i - n) - 1}v_{i'} \otimes v_{i}.\label{eq:B_hwv}
\end{gather}
The equality $e_{i}(w_{1})=0$ for all $i$ follows immediately from
\[
 e_{i}\cdot v_{j} =
 \begin{cases}
 v_{i} & \text{if}\ j = i + 1 ,\\
 -v_{(i + 1)'} & \text{if}\ j = i', \\
 0 & \text{otherwise}.
 \end{cases}
\]
Likewise, since $e_{i}(v_{2})=0$ unless $i=1$ and $e_1(v_2)=v_1$, $\omega_1(v_1)=\bigl(\omega_{\varepsilon_{1}}',\omega_{1}\bigr)v_1$, we also obtain
\[
 e_{i}(w_{2}) =
 \omega_{i}(v_{1}) \otimes e_{i}(v_{2}) - (\omega_{\varepsilon_{1}}',\omega_{1}) e_{i}(v_{2}) \otimes v_{1} = 0
 \qquad \mathrm{for\ all} \ 1\leq i\leq n.
\]
It remains to check that $w_3$ is a highest weight vector. Indeed, for $k<n$ we get
\begin{align*}
 e_{k}(w_{3})&{}=
 r^{2k}v_{k}\otimes v_{(k+1)'} - r^{2n-1}s^{2(k-n)-1} v_{(k+1)'}\otimes v_{k}
 \\
 &\quad{}- r^{2(k-1)} \omega_k(v_{k})\otimes v_{(k+1)'} + r^{2n-1}s^{2(k-n)+1} \omega_k\bigl(v_{(k + 1)'}\bigr) \otimes v_{k} \\
 &{}= r^{2k}v_{k} \otimes v_{(k + 1)'} - r^{2(k - 1)}r^{2}v_{k} \otimes v_{(k + 1)'} \\
 &\quad{}-
 r^{2n - 1}s^{2(k - n)-1}v_{(k + 1)'} \otimes v_{k} + r^{2n - 1}s^{2(k - n) + 1}s^{-2}v_{(k + 1)'} \otimes v_{k} = 0,
\end{align*}
and similarly for $k = n$, we have
\begin{displaymath}
\begin{split}
 e_{n}(w_{3}) &=
 r^{2n-1}s^{-1} v_{n} \otimes v_{n+1} - r^{2n-1}s^{-1} v_{n+1} \otimes v_{n}
 \\
 &\quad- r^{2(n-1)} \omega_n(v_{n}) \otimes v_{n+1} + r^{2n-1}s^{-1} \omega_n(v_{n+1}) \otimes v_{n} \\
 &= r^{2n - 1}s^{-1}v_{n} \otimes v_{n+1}- rs^{-1}\cdot r^{2(n-1)}v_{n} \otimes v_{n+1} \\
 &\quad - r^{2n - 1}s^{-1}v_{n+1} \otimes v_{n} + r^{2n - 1}s^{-1}v_{n+1} \otimes v_{n} = 0.
\end{split}
\end{displaymath}
It is clear that $w_{1}$, $w_{2}$, $w_{3}$ are linearly independent, and since $w_{1}$ has weight $2\varepsilon_{1}$, $w_{2}$ has weight ${\varepsilon_{1} + \varepsilon_{2}}$, and $w_{3}$ has weight $0$, the result of the Proposition~will follow if we can match the corresponding dimensions:
\begin{equation}\label{eq:dim-sum}
 \dim{L(\varepsilon_{1})} + \dim{L(\varepsilon_{1} + \varepsilon_{2})} + \dim{L(0)} = (2n + 1)^{2}.
\end{equation}
In view of Proposition~\ref{lem:dim}, we can prove this by using the Weyl dimension formula. In type $B_{n}$, we have
\begin{displaymath}
 \rho = \frac{1}{2}\sum_{\alpha \in \Phi^{+}} \alpha =
 \sum_{i = 1}^{n}\biggl(n + \frac{1}{2} - i\biggr)\varepsilon_{i},
\end{displaymath}
and therefore
\begin{displaymath}
 \bigl(\tfrac{1}{2}\rho,\varepsilon_{i} - \varepsilon_{j}\bigr) = j - i, \qquad
 \bigl(\tfrac{1}{2}\rho,\varepsilon_{i} + \varepsilon_{j}\bigr) = 2n + 1 - i - j, \qquad
 \bigl(\tfrac{1}{2}\rho,\varepsilon_{i}\bigr) = n + \tfrac{1}{2} - i.
\end{displaymath}
Consequently, we obtain:
\begin{displaymath}
 \prod_{\alpha \in \Phi^{+}}\bigl(\tfrac{1}{2}\rho ,\alpha\bigr) =
 \Biggl( \prod_{i = 1}^{n} \bigl(n + \tfrac{1}{2} - i \bigr) \Biggr)
 \biggl( \prod_{1 \le i < j \le n}(j - i)(2n + 1 - j - i) \biggr).
\end{displaymath}
A similar computation yields
\begin{gather*}
 \prod_{\alpha \in \Phi^{+}}\bigl(\varepsilon_{1} + \tfrac{1}{2}\rho,\alpha\bigr) \\
 \qquad{}=
 \Biggl( \prod_{i = 2}^{n}(1 + i)(2n + 2 - i)\bigl(n + \tfrac{1}{2} - i\bigr)\Biggr)
 \bigl( n + \tfrac{3}{2} \bigr) \biggl( \prod_{2 \le i < j \le n}(j - i)(2n + 1 - i - j) \biggr),
\end{gather*}
and thus we get
\begin{displaymath}
 \dim{L(2\varepsilon_{1})} =
  \frac{\prod_{\alpha \in \Phi^{+}}(2\varepsilon_{1} + \rho,\alpha)}{\prod_{\alpha \in \Phi^{+}}(\rho,\alpha)} =
  \frac{\prod_{\alpha \in \Phi^{+}}\bigl(\varepsilon_{1} + \tfrac{1}{2}\rho,\alpha\bigr)}{\prod_{\alpha \in \Phi^{+}}\bigl(\tfrac{1}{2}\rho,\alpha\bigr)} =
  n(2n + 3).
\end{displaymath}
Similarly, we obtain
\begin{displaymath}
 \dim{L(\varepsilon_{1} + \varepsilon_{2})} = n(2n + 1).
\end{displaymath}
This completes the proof of the equality~\eqref{eq:dim-sum} and hence of the Proposition, since $\dim{L(0)} = 1$.

\textit{Type $C_n$}.
Highest weight vectors for the $U_{r,s}(\ssp_{2n})$-action on $V\otimes V$ are
\begin{gather}
 w_{1} = v_{1} \otimes v_{1},\qquad
 w_{2} = v_{1} \otimes v_{2} - \bigl(\omega_{\varepsilon_{1}}',\omega_{1}\bigr)v_{2} \otimes v_{1} =
 v_{1} \otimes v_{2} - r v_{2} \otimes v_{1} , \nonumber\\
 w_{3} = \sum_{i = 1}^{n}\bigl(r^{i - 1}v_{i} \otimes v_{i'} - r^{n}s^{i - n- 1}v_{i'} \otimes v_{i}\bigr),\label{eq:C_hwv}
\end{gather}
and we have
\[
 \dim{L(2\varepsilon_{1})} = n(2n + 1) , \qquad
 \dim{L(\varepsilon_{1} + \varepsilon_{2})} = (2n + 1)(n - 1), \qquad \dim{L(0)}=1 .
\]

\textit{Type $D_n$}.
Highest weight vectors for the $U_{r,s}(\sso_{2n})$-action on $V\otimes V$ are
\begin{gather}
 w_{1} = v_{1} \otimes v_{1},\qquad
 w_{2} = v_{1} \otimes v_{2} - \bigl(\omega_{\varepsilon_{1}}',\omega_{1}\bigr)v_{2} \otimes v_{1} =
 v_{1} \otimes v_{2} - r v_{2} \otimes v_{1}, \nonumber\\
 w_{3} = \sum_{i = 1}^{n}\bigl(r^{i - 1}v_{i} \otimes v_{i'} + r^{n - 1}s^{i - n}v_{i'} \otimes v_{i}\bigr),\label{eq:D_hwv}
\end{gather}
and we have
\begin{displaymath}
 \dim{L(2\varepsilon_{1})} = (2n - 1)(n + 1), \qquad
 \dim{L(\varepsilon_{1} + \varepsilon_{2})} = n(2n - 1), \qquad \dim{L(0)}=1 . \tag*{\qed}
\end{displaymath}
\renewcommand{\qed}{}
\end{proof}

\section[R-matrices]{$\boldsymbol{R}$-matrices}\label{sec:R-matrices}

In this section, we evaluate the $\uu$-module isomorphism $V \otimes V \iso V \otimes V$ for $\fg$ one of~$\sso_{2n + 1}$, $\ssp_{2n}$, $\sso_{2n}$ and their first fundamental representation $V$ from Section~\ref{sec:column_repr}, arising through the universal $R$-matrix. This produces two-parametric solutions of the quantum Yang--Baxter equation for classical Lie algebras, cf.\ \eqref{eq:qYB-intro}:
\begin{equation}\label{eq:qYB}
 R_{12}R_{13}R_{23}=R_{23}R_{13}R_{12}.
\end{equation}

\subsection{Universal construction}\label{ssec:universal-R}

Let us first recall the general construction~\eqref{eq:intertwiner-intro} of a $\uu$-module isomorphism $V \otimes W \to W \otimes V$ arising through the universal $R$-matrix (see~\cite[Section~4]{BW1} for $U_{r,s}(\ssl_n)$ and~\cite[Section~3]{BGH2} for other classical $U_{r,s}(\fg)$, both modeled after the treatment of one-parameter case in~\cite[Section~7]{J}). To~this end, we pick dual bases $\bigl\{x_{i}^{\mu}\bigr\}$ and $\bigl\{y_{i}^{\mu}\bigr\}$ of $U^{+}_{r,s}(\fg)_{\mu}$ and $U^{-}_{r,s}(\fg)_{-\mu}$ with respect to the Hopf pairing~\eqref{eq:Hopf-parity}, and set
\begin{equation}\label{eq:Theta}
 \Theta = 1 + \sum_{\mu > 0}\Theta_{\mu} \qquad \text{with} \
 \Theta_{\mu} = \sum_{i} y_{i}^{\mu} \otimes x_{i}^{\mu}.
\end{equation}
Let $\tau\colon V \otimes W \to W \otimes V$ be the flip map $v \otimes w \mapsto w \otimes v$.
Finally, consider $f\colon P \times P \to \mathbb{K}^{\times}$ satisfying
\begin{alignat}{3}
 &f(\lambda + \nu,\mu) = f(\lambda,\mu)f(\nu,\mu), \qquad&&
 f(\lambda,\mu + \nu) = f(\lambda,\mu)f(\lambda,\nu) , &\nonumber \\
 &f(\lambda,\alpha_{i}) = \bigl(\omega_{i}',\omega_{\lambda}\bigr)^{-1}, \qquad&&
 f(\alpha_{i},\mu) = \bigl(\omega_{\mu}',\omega_{i}\bigr)^{-1}&\label{eq:f}
\end{alignat}
for all $\lambda,\mu \in P$, $\nu \in Q$, and $\alpha_{i} \in \Pi$. Then, for any two $\uu$-modules $V$ and $W$ with weight space decomposition, we define a linear map \smash{$\widetilde{f}\colon V \otimes W \to V \otimes W$} via \smash{$\widetilde{f}(v \otimes w) = f(\lambda,\mu) \cdot v \otimes w$} if $v \in V[\lambda]$, $w \in W[\mu]$. The following is standard.

\begin{Theorem}\label{thm:universal-R}
For any finite-dimensional $\uu$-modules $V$ and $W$, the map
\begin{equation}\label{eq:R-int}
 \hat{R}_{VW}=\Theta \circ \widetilde{f} \circ \tau\colon\ V \otimes W \to W \otimes V
\end{equation}
is an isomorphism of $\uu$-modules.
\end{Theorem}

Let \smash{$R_{VW}=\hat{R}_{WV} \circ \tau=\Theta \circ \widetilde{f}\colon V\otimes W\to V\otimes W$}. Given finite-dimensional $\uu$-modules $V_1$, $V_2$, $V_3$, define three endomorphisms of $V_1\otimes V_2\otimes V_3$: $R_{12}=R_{V_1,V_2}\otimes \mathrm{Id}_{V_3}$, $R_{23}=\mathrm{Id}_{V_1}\otimes R_{V_2,V_3}$, $R_{13}=(\mathrm{Id}\otimes \tau)R_{12}(\mathrm{Id}\otimes \tau)$. We likewise define linear operators \smash{$\hat{R}_{12}$}, \smash{$\hat{R}_{23}$}, \smash{$\hat{R}_{13}$}. According to~\cite{BW1,BGH2}, modelled after~\cite{J}, we have
\begin{align}
 & R_{12} R_{13} R_{23} = R_{23} R_{13} R_{12}\colon\
 V_1\otimes V_2\otimes V_3 \to V_1\otimes V_2\otimes V_3,\nonumber \\
 & \hat{R}_{12} \hat{R}_{23} \hat{R}_{12} = \hat{R}_{23} \hat{R}_{12} \hat{R}_{23}\colon\
 V_1\otimes V_2\otimes V_3 \to V_3\otimes V_2\otimes V_1.\label{eq:qYB-two}
\end{align}
In particular, we obtain a whole family of solutions of the quantum Yang--Baxter equation.

\begin{Corollary}
For any finite-dimensional $\uu$-module $V$, the operator $R_{VV}=\hat{R}_{VV} \circ \tau$ satisfies~\eqref{eq:qYB}.
\end{Corollary}

\subsection[Explicit R-matrices]{Explicit $\boldsymbol{R}$-matrices}

For the representation $V$ of Proposition~\ref{prp:A_rep}, the explicit formula for $\hat{R}_{VV}$ was obtained in~\cite[Sec\-tion~5]{BW2}.

\begin{Theorem}[type $A_{n}$]\label{thm:A_RMatrix}
The $U_{r,s}(\ssl_{n+1})$-module isomorphism $\hat{R}_{VV}\colon V \otimes V \iso V \otimes V$ from Theorem~$\ref{thm:universal-R}$ for the $U_{r,s}(\ssl_{n+1})$-module $V$ from Proposition~$\ref{prp:A_rep}$ coincides with the following operator:
\begin{gather}
 \hat{R} =
 \sum_{i=1}^{n+1} E_{ii} \otimes E_{ii}
 + r \sum_{1\leq i<j\leq n+1} E_{ji} \otimes E_{ij}+ s^{-1} \sum_{1\leq i<j\leq n+1} E_{ij} \otimes E_{ji}\nonumber\\
 \hphantom{\hat{R} =}{}
 + \bigl(1-rs^{-1}\bigr) \sum_{1\leq i<j\leq n+1} E_{jj} \otimes E_{ii} .\label{eq:A_RMatrix}
\end{gather}
\end{Theorem}

The main results of this section generalize the above formula to the other classical series.

\begin{Theorem}[type $B_{n}$]\label{thm:B_RMatrix}
The $U_{r,s}(\sso_{2n+1})$-module isomorphism $\hat{R}_{VV}\colon V \otimes V \iso V \otimes V$ from Theorem~$\ref{thm:universal-R}$ for the $U_{r,s}(\sso_{2n+1})$-module $V$ from Proposition~$\ref{prp:B_rep}$ coincides with the following operator:
\begin{align}
 \hat{R} &{}=
 r^{-1}s \sum_{1\leq i\leq 2n+1}^{i\ne n+1} E_{ii} \otimes E_{ii} +
 E_{n + 1,n + 1} \otimes E_{n + 1,n+1} +
 rs^{-1} \sum_{1\leq i\leq 2n+1}^{i\ne n+1} E_{ii'} \otimes E_{i'i} \nonumber\\
 &\quad{}+ \sum_{1\leq i,j\leq 2n+1}^{j\ne i,i'} a_{ij}E_{ij} \otimes E_{ji} +
 \bigl(r^{2} - s^{2}\bigr)(rs)^{-1} \sum_{i = 1}^{n}\bigl(r^{2(n-i)+1}s^{2(i-n)-1}-1\bigr) E_{i'i'} \otimes E_{ii}\nonumber \\
 &\quad{}+ \bigl(s^{2} - r^{2}\bigr)(rs)^{-1} \sum_{i > j}^{j \neq i'} E_{ii} \otimes E_{jj} +
 \bigl(r^{2} - s^{2}\bigr)(rs)^{-1} \sum_{i < j}^{j \neq i'} t_{i}t_{j}^{-1}E_{i'j} \otimes E_{ij'},\label{eq:B_RMatrix}
\end{align}
with the constants $t_i$ and $a_{ij}$ given explicitly by
\begin{gather}\label{eq:B_coeffs}
 t_{i} =
 \begin{cases}
 s^{2(i - n) - 1} & \text{if}\ i < n + 1, \\
 s^{-1} & \text{if}\ i = n + 1, \\
 r^{2(n +1 - i) + 1} & \text{if}\ i > n + 1,
 \end{cases}\qquad
 a_{ij} =
 \begin{cases}
 (rs)^{-\sigma_{i}\sigma_{j}} & \text{if}\ i < j,j'\ \text{or}\ i > j,j', \\
 (rs)^{\sigma_{i}\sigma_{j}} & \text{if}\ j < i < j'\ \text{or}\ j' < i < j,
 \end{cases}
\end{gather}
where we set
\[
 \sigma_{i} =
 \begin{cases}
 -1 & \text{if}\ i < n + 1, \\
 0 & \text{if}\ i = n+ 1, \\
 1 & \text{if}\ i > n + 1.
 \end{cases}
\]
\end{Theorem}

\begin{Theorem}[type $C_{n}$]\label{thm:C_RMatrix}
The $U_{r,s}(\ssp_{2n})$-module isomorphism $\hat{R}_{VV}\colon V \otimes V \iso V \otimes V$ from Theorem~$\ref{thm:universal-R}$ for the $U_{r,s}(\ssp_{2n})$-module $V$ from Proposition~$\ref{prp:C_rep}$ coincides with the following operator:
\begin{align}
 \hat{R} &{}=
 r^{-1/2}s^{1/2}\sum_{i = 1}^{2n} E_{ii} \otimes E_{ii} +
 r^{1/2}s^{-1/2}\sum_{i = 1}^{2n} E_{ii'} \otimes E_{i'i}+ \sum_{1\leq i,j\leq 2n}^{j \ne i,i'} a_{ij}E_{ij} \otimes E_{ji} \nonumber\\
 &\quad{}
 + (s - r)(rs)^{-1/2}\sum_{i = 1}^{n} \bigl(r^{n+ 1 - i}s^{i - n - 1} + 1\bigr)E_{i'i'} \otimes E_{ii} \nonumber\\
 &\quad{} + (s - r)(rs)^{-1/2}\sum_{i > j}^{j \neq i'} E_{ii} \otimes E_{jj}
 + (r - s)(rs)^{-1/2}\sum_{i < j}^{j \neq i'} t_{i}t_{j}^{-1}E_{i'j} \otimes E_{ij'} ,\label{eq:C_RMatrix}
\end{align}
with the constants $t_i$ and $a_{ij}$ given explicitly by
\begin{equation}\label{eq:C_coeffs}
 t_{i} = \
 \begin{cases}
 s^{i - n - 1} & \text{if}\ i \le n, \\
 -r^{n - i} & \text{if}\ i > n,
 \end{cases} \qquad
 a_{ij} =
 \begin{cases}
 (rs)^{-\frac{1}{2}\sigma_{i}\sigma_{j}} & \text{if}\ i < j,j'\ \text{or}\ i > j,j', \\
 (rs)^{\frac{1}{2}\sigma_{i}\sigma_{j}} & \text{if}\ j < i < j'\ \text{or}\ j' < i < j,
 \end{cases}
\end{equation}
where we set
\begin{displaymath}
 \sigma_{i} =
 \begin{cases}
 1 & \text{if}\ i \le n, \\
 -1 & \text{if}\ i > n.
 \end{cases}
\end{displaymath}
\end{Theorem}

\begin{Theorem}[type $D_{n}$]\label{thm:D_RMatrix}
The $U_{r,s}(\sso_{2n})$-module isomorphism $\hat{R}_{VV}\colon V \otimes V \iso V \otimes V$ from Theorem~$\ref{thm:universal-R}$ for the $U_{r,s}(\sso_{2n})$-module $V$ from Proposition~$\ref{prp:D_rep}$ coincides with the following operator:
\begin{align}
 \hat{R} &{}=
 r^{-1/2}s^{1/2}\sum_{i = 1}^{2n} E_{ii} \otimes E_{ii}
 + r^{1/2}s^{-1/2}\sum_{i = 1}^{2n}E_{ii'} \otimes E_{i'i}+ \sum_{1\leq i,j\leq 2n}^{j \neq i,i'} a_{ij}E_{ij} \otimes E_{ji} \nonumber\\
 &\quad{}
 + (s - r)(rs)^{-1/2}\sum_{i = 1}^{n} \bigl(1 - r^{n- i}s^{i - n}\bigr) E_{i'i'} \otimes E_{ii} \nonumber\\
 &\quad{} + (s -r)(rs)^{-1/2}\sum_{i > j}^{j \neq i'} E_{ii} \otimes E_{jj}
 + (r - s)(rs)^{-1/2}\sum_{i < j}^{j \neq i'} t_{i}t_{j}^{-1}E_{i'j} \otimes E_{ij'} ,\label{eq:D_RMatrix}
 \end{align}
with the constants $t_i$ and $a_{ij}$ given explicitly by
\begin{equation}\label{eq:D_coeffs}
 t_{i} =
 \begin{cases}
 s^{i - n} & \text{if}\ i \le n, \\
 r^{n + 1 - i} & \text{if}\ i > n,
 \end{cases} \qquad
 a_{ij} =
 \begin{cases}
 (rs)^{-\frac{1}{2}\sigma_{i}\sigma_{j}} & \text{if}\ i < j,j'\ \text{or}\ i > j,j', \\
 (rs)^{\frac{1}{2}\sigma_{i}\sigma_{j}} & \text{if}\ j < i < j'\ \text{or}\ j' < i < j,
 \end{cases}
\end{equation}
where we set
\begin{displaymath}
 \sigma_{i} =
 \begin{cases}
 1 & \text{if}\ i \le n, \\
 -1 & \text{if}\ i > n.
 \end{cases}
\end{displaymath}
\end{Theorem}

\begin{Remark}
Although the proofs of Theorems~\ref{thm:B_RMatrix}--\ref{thm:D_RMatrix} share many similarities with that of Theorem~\ref{thm:A_RMatrix}, there are slight differences. Indeed, instead of explicitly determining bases of each irreducible component from the decomposition~\eqref{eq:irr-decomp} and computing the action of $\hat{R}$ on these basis vectors, we rather verify that operators $\hat{R}$ from~\eqref{eq:B_RMatrix},~\eqref{eq:C_RMatrix} and~\eqref{eq:D_RMatrix} do commute with the action of all $f_i$ acting on $V\otimes V$, and act on the three highest weight vectors by the desired scalars. This approach is essential to our construction of $\hat{R}(z)$ in Section~\ref{sec:affine-R}.
\end{Remark}

\subsection{Proofs of explicit formulas}

To prove theorems in the previous subsection, we shall first study the operators $\hat{R}$ featured in~\eqref{eq:B_RMatrix},~\eqref{eq:C_RMatrix} and~\eqref{eq:D_RMatrix}. Our first technical result evaluates $\hat{R}$-action on the three highest weight vectors of $V\otimes V$, see Proposition~\ref{prop:struct}.

\begin{Lemma}[type $B_{n}$]\label{lem:B_eigen}
The highest weight vectors $w_{1}, w_{2}, w_{3}\in V\otimes V$ from~\eqref{eq:B_hwv} are eigenvectors of the operator $\hat{R}$ from~\eqref{eq:B_RMatrix}, with respective eigenvalues $\lambda_{1}=r^{-1}s$, $\lambda_{2} = -rs^{-1}$, $\lambda_{3} = r^{2n}s^{-2n}$.
\end{Lemma}

\begin{Lemma}[type $C_{n}$]\label{lem:C_eigen}
The highest weight vectors $w_{1}, w_{2}, w_{3}\in V\otimes V$ from~\eqref{eq:C_hwv} are eigenvectors for the operator $\hat{R}$ from~\eqref{eq:C_RMatrix}, with respective eigenvalues $\lambda_{1} = r^{-1/2}s^{1/2}$, $\lambda_{2} = -r^{1/2}s^{-1/2}$, $\lambda_{3} = -r^{n + 1/2}s^{-n-1/2}$.
\end{Lemma}

\begin{Lemma}[type $D_{n}$]\label{lem:D_eigen}
The highest weight vectors $w_{1}, w_{2}, w_{3}\in V\otimes V$ from~\eqref{eq:D_hwv} are eigenvectors of the operator $\hat{R}$ from~\eqref{eq:D_RMatrix}, with respective eigenvalues $\lambda_{1} = r^{-1/2}s^{1/2}$, $\lambda_{2} = -r^{1/2}s^{-1/2}$, $\lambda_{3} = r^{n-1/2}s^{-n+1/2}$.
\end{Lemma}

As the proofs of all three results are completely analogous, we shall prove only the first one.

\begin{proof}[Proof of Lemma \ref{lem:B_eigen}]
For $w_{1}$, we clearly get
\begin{displaymath}
 \hat{R}(v_{1} \otimes v_{1}) = r^{-1}s v_{1} \otimes v_{1} .
\end{displaymath}
Likewise, if $n>1$, then for $w_{2}$ we obtain
\begin{align*}
 \hat{R}\bigl(v_{1} \otimes v_{2} - r^{2}v_{2} \otimes v_{1}\bigr) &{}=
 rsv_{2} \otimes v_{1} - rs^{-1}v_{1} \otimes v_{2} - \bigl(s^{2} - r^{2}\bigr)rs^{-1}v_{2} \otimes v_{1} \\
 &{}=
 -rs^{-1}\bigl(v_{1} \otimes v_{2} - r^{2}v_{2} \otimes v_{1}\bigr).
\end{align*}
If $n = 1$, then we have $a_{12} = a_{21} = 1$, and so we get
\begin{displaymath}
\begin{split}
 \hat{R}\bigl(v_{1} \otimes v_{2} - rs^{-1}v_{2} \otimes v_{1}\bigr) &{}= v_{2} \otimes v_{1} - rs^{-1}v_{1} \otimes v_{2} - rs^{-1}\bigl(r^{-1}s - rs^{-1}\bigr)v_{2} \otimes v_{1} \\
 &{}= -rs^{-1}\bigl(v_{1} \otimes v_{2} - rs^{-1}v_{2} \otimes v_{1}\bigr).
\end{split}
\end{displaymath}
Finally, for $w_{3}$ we obtain
\begin{gather*}
 \hat{R} \Biggl( \sum_{i = 1}^{n}r^{2(i - 1)}v_{i} \otimes v_{i'} +
 r^{2n -1}s^{-1}v_{n+1} \otimes v_{n+1} +
 \sum_{i = 1}^{n}r^{2n - 1}s^{-2(n - i) - 1}v_{i'} \otimes v_{i} \Biggr) \\
 \qquad{} = \sum_{i = 1}^{n} r^{2i - 1}s^{-1}v_{i'} \otimes v_{i} +
 r^{2n - 1}s^{-1}v_{n+1} \otimes v_{n+1} +
 \sum_{i = 1}^{n}r^{2n}s^{2(i - n) - 2} v_{i} \otimes v_{i'} \\
 \qquad\hphantom{=}{} + \bigl(rs^{-1} - r^{-1}s\bigr)\sum_{i = 1}^{n}\bigl(r^{2(n - i) + 1}s^{2(i - n) - 1} - 1\bigr)r^{2n - 1}s^{2(i - n) - 1} v_{i'} \otimes v_{i} \\
 \qquad\hphantom{=}{} + \bigl(rs^{-1} - r^{-1}s\bigr)\sum_{i = 1}^{n} r^{2n - 1}s^{2(i - n) - 1} v_{i'} \otimes v_{i} \\
 \qquad\hphantom{=}{} + \bigl(rs^{-1} - r^{-1}s\bigr)\sum_{i = 1}^{n}\sum_{j = 1}^{i - 1} r^{2(i - 1)}s^{2(j - i)} v_{j'} \otimes v_{j} \\
 \qquad\hphantom{=}{} + \bigl(rs^{-1} - r^{-1}s\bigr) \sum_{i = 1}^{n} \sum_{1\leq j\leq (i+1)'}^{j \neq i}
 t_{j}t_{i'}^{-1}r^{2n - 1}s^{2(i - n) - 1} v_{j'} \otimes v_{j} ,
\end{gather*}
where the last three summands arise from the action of the last sum in~\eqref{eq:B_RMatrix}.

The last two summands simplify as follows:
\begin{align*}
 \sum_{i = 1}^{n}\sum_{j = 1}^{i - 1} r^{2(i - 1)}s^{2(j - i)} v_{j'} \otimes v_{j} &{}=
 \sum_{j = 1}^{n-1} \Biggl( \sum_{i = j + 1}^{n}r^{2(i - 1)}s^{2(j - i)} \Biggr) v_{j'}\otimes v_{j}\\
 &{} =
 \sum_{j = 1}^{n}r^{2j}s^{2(j - n)}[n-j]_{r^{2},s^{2}}v_{j'} \otimes v_{j}
\end{align*}
and
\begin{gather*}
 \sum_{i = 1}^{n}\sum_{1\leq j\leq (i + 1)'}^{j \neq i} t_{j}t_{i'}^{-1}r^{2n - 1}s^{2(i - n) - 1}v_{j'} \otimes v_{j} \\
 \qquad{}= \sum_{i = 1}^{n}\sum_{j = 1}^{n} r^{4n - 2i}s^{2(i + j - 2n - 1)}v_{j'} \otimes v_{j} -
 \sum_{i = 1}^{n} r^{4n - 2i}s^{4(i - n) - 2} v_{i'} \otimes v_{i} \\
 \qquad\hphantom{=}{} + \sum_{i = 1}^{n}\sum_{j = i + 1}^{n} r^{2(j - i + n) - 1}s^{2(i - n) - 1}v_{j} \otimes v_{j'}
 + \sum_{i = 1}^{n} r^{4n -2i}s^{2(i - n - 1)} v_{n + 1} \otimes v_{n + 1} \\
 \qquad{} = \sum_{j = 1}^{n}\bigl(r^{2n}s^{2j - 4n}[n]_{r^{2},s^{2}} - r^{4n - 2j}s^{4(j - n) - 2}\bigr) v_{j'} \otimes v_{j} \\
 \qquad\hphantom{=}{} + \sum_{j = 1}^{n} r^{2n + 1}s^{1 - 2n}[j - 1]_{r^{2},s^{2}} v_{j} \otimes v_{j'}+ r^{2n}s^{-2n}[n]_{r^{2},s^{2}}v_{n + 1} \otimes v_{n + 1}.
\end{gather*}

Thus, the coefficient of $v_{i'} \otimes v_{i}$ in $\hat{R}(w_3)$ for $1\leq i\leq n$ equals
\begin{gather*}
 r^{2i - 1}s^{-1} +
 \bigl(rs^{-1} - r^{-1}s\bigr)\bigl(r^{2(n - i) + 1}s^{2(i - n) - 1} - 1\bigr)r^{2n - 1}s^{2(i-n) - 1}\\
 \quad\hphantom{=}{} + \bigl(rs^{-1} - r^{-1}s\bigr)r^{2n - 1}s^{2(i - n) - 1}+r^{2i - 1}s^{2(i - n) - 1}\bigl(r^{2(n - i)} - s^{2(n - i)}\bigr) \\
 \quad\hphantom{=}{} + r^{2n - 1}s^{2i - 4n - 1}\bigl(r^{2n} - s^{2n}\bigr) - r^{4n - 2i - 1}s^{4(i - n) - 3}\bigl(r^{2} - s^{2}\bigr) \\
 \quad{} = r^{4n - 1}s^{2i - 4n - 1} = r^{2n}s^{-2n}\cdot r^{2n - 1}s^{2(i - n) - 1},
\end{gather*}
as desired.
Likewise, the coefficient of $v_{i} \otimes v_{i'}$ in $\hat{R}(w_3)$ for $1\leq i\leq n$ equals $r^{2(n + i - 1)}s^{-2n} = r^{2n}s^{-2n}\cdot r^{2(i - 1)}$, while the coefficient of $v_{n+ 1} \otimes v_{n+1}$ equals $r^{4n - 1}s^{-2n - 1} = r^{2n}s^{-2n}\cdot r^{2n - 1}s^{-1}$. This completes the proof.
\end{proof}

Next, we verify that the operators $\hat{R}$ are indeed $\uu$-module homomorphisms.

\begin{Lemma}\label{lem:B_Intertwining}
The operators $\hat{R}\colon V \otimes V \to V \otimes V$ from~\eqref{eq:B_RMatrix},~\eqref{eq:C_RMatrix} and~\eqref{eq:D_RMatrix} are isomorphisms of $\uu$-modules.
\end{Lemma}

\begin{proof}
In each case, it suffices to verify $\Delta(f_{k})\hat{R} = \hat{R}\Delta(f_{k})\in \End(V \otimes V)$ for all $k$, since Lemmas~\ref{lem:B_eigen}--\ref{lem:D_eigen} then imply that $\hat{R}$ acts as a nonzero scalar on each irreducible component of $V \otimes V$ (see Proposition~\ref{prop:struct}). We will present this verification only in type $B_{n}$, since the arguments in the remaining cases are similar.

To make the computations more manageable, it will be helpful to break the operator $\hat{R}$ from~\eqref{eq:B_RMatrix} into the following six pieces:
\begin{gather*}
 R_{1} = r^{-1}s \sum_{1\leq i\leq 2n+1}^{i \neq n + 1} E_{ii} \otimes E_{ii}
 + E_{n + 1,n + 1} \otimes E_{n + 1,n+1}, \\
 R_{2} = rs^{-1} \sum_{1\leq i\leq 2n+1}^{i \neq n + 1} E_{ii'} \otimes E_{i'i}, \\
 R_{3} = \sum_{1\leq i,j\leq 2n+1}^{j \neq i,i'} a_{ij}E_{ij} \otimes E_{ji}, \\
 R_{4} = \bigl(r^{2} - s^{2}\bigr)(rs)^{-1} \sum_{i = 1}^{n}\bigl(r^{2(n - i) + 1}s^{2(i - n) - 1} - 1\bigr) E_{i'i'} \otimes E_{ii}, \\
 R_{5} = \bigl(s^{2} - r^{2}\bigr)(rs)^{-1} \sum_{i > j}^{j \neq i'} E_{ii} \otimes E_{jj},\\
 R_{6} = \bigl(r^{2} - s^{2}\bigr)(rs)^{-1} \sum_{i < j}^{j \neq i'} t_{i}t_{j}^{-1}E_{i'j} \otimes E_{ij'}.
\end{gather*}

Now, for $k < n$, the matrix of $\Delta(f_{k})$ is
\begin{align*}
 \Delta(f_{k}) &{}= 1 \otimes E_{k + 1,k} - (rs)^{-2}1 \otimes E_{k',(k + 1)'} + s^{2}E_{k + 1,k} \otimes E_{kk} \\
 &\quad{} + r^{2}E_{k + 1,k} \otimes E_{k + 1,k + 1} + s^{-2}E_{k + 1,k} \otimes E_{k'k'} + r^{-2}E_{k + 1,k} \otimes E_{(k + 1)',(k+ 1)'} \\
 &\quad{} + \sum_{1\leq j\leq n}^{j\ne k,k + 1} \bigl(E_{k + 1,k} \otimes E_{jj} + E_{k + 1,k} \otimes E_{j'j'}\bigr) + E_{k + 1,k} \otimes E_{n + 1,n + 1} \\
 &\quad{} - r^{-2}E_{k',(k + 1)'} \otimes E_{kk} - s^{-2}E_{k',(k + 1)'} \otimes E_{k + 1,k + 1} - r^{-2}s^{-4}E_{k',(k + 1)'} \otimes E_{k'k'} \\
 &\quad{} - r^{-4}s^{-2}E_{k',(k + 1)'} \otimes E_{(k + 1)',(k + 1)'} \\
 &\quad{} -(rs)^{-2}\sum_{1\leq j\leq n}^{j\ne k,k + 1} \bigl(E_{k',(k + 1)'} \otimes E_{jj} + E_{k',(k + 1)'} \otimes E_{j'j'}\bigr) \\
 &\quad{} - (rs)^{-2}E_{k',(k + 1)'} \otimes E_{n + 1,n+1}.
\end{align*}
Thus, by direct computation, we get
\begin{gather*}
 R_{1}\Delta(f_{k}) = r^{-1}sE_{k + 1,k + 1} \otimes E_{k + 1,k} - r^{-3}s^{-1}E_{k'k'} \otimes E_{k',(k + 1)'} \\[0.2mm]
 \hphantom{R_{1}\Delta(f_{k}) =}{}
 + rsE_{k+1,k} \otimes E_{k+1,k + 1} - r^{-3}s^{-3}E_{k',(k + 1)'} \otimes E_{k'k'},
 \\[0.2mm]
 \Delta(f_{k})R_{1} = r^{-1}sE_{kk} \otimes E_{k + 1,k} - r^{-3}s^{-1}E_{(k + 1)',(k + 1)'} \otimes E_{k',(k + 1)'} \\[0.2mm]
 \hphantom{\Delta(f_{k})R_{1} =}{}
 + r^{-1}s^{3}E_{k + 1,k} \otimes E_{kk} - r^{-5}s^{-1}E_{k',(k + 1)'} \otimes E_{(k + 1)',(k + 1)'},
\\[0.2mm]
 R_{2}\Delta(f_{k}) = rs^{-1}E_{k + 1,(k + 1)'} \otimes E_{(k + 1)',k} - r^{-1}s^{-3}E_{k',k} \otimes E_{k,(k + 1)'} \\[0.2mm]
 \hphantom{R_{2}\Delta(f_{k}) =}{}
 + (rs)^{-1}E_{(k + 1)',k} \otimes E_{k + 1,(k + 1)'} - (rs)^{-1}E_{k,(k + 1)'} \otimes E_{k'k},
 \\[0.2mm]
 \Delta(f_{k})R_{2} = rs^{-1}E_{k'k} \otimes E_{k + 1,k'} - r^{-1}s^{-3}E_{k + 1,(k + 1)'} \otimes E_{k',k + 1} \\[0.2mm]
 \hphantom{\Delta(f_{k})R_{2}=}{} + rs^{-3}E_{k + 1,k'} \otimes E_{k'k} - rs^{-3}E_{k',k + 1} \otimes E_{k + 1,(k + 1)'},
\\[0.2mm]
 R_{3}\Delta(f_{k}) - \Delta(f_{k})R_{3} \\[0.2mm]
 \qquad{}= rsE_{k + 1,k} \otimes E_{kk} + (rs)^{-1}E_{k + 1,k'} \otimes E_{k'k} - (rs)^{-1}E_{k',k + 1} \otimes E_{k + 1,(k + 1)'} \\[0.2mm]
 \quad\qquad{} - (rs)^{-3}E_{k',(k + 1)'} \otimes E_{(k + 1)',(k + 1)'} + r^{-1}sE_{kk} \otimes E_{k + 1,k} + rs^{-1}E_{k'k} \otimes E_{k + 1,k'} \\[0.2mm]
 \quad\qquad{} - rsE_{k + 1,k} \otimes E_{k + 1,k+1} - (rs)^{-1}E_{(k + 1)',k} \otimes E_{k + 1,(k + 1)'} + (rs)^{-1}E_{k,(k + 1)'} \otimes E_{k'k} \\[0.2mm]
 \quad\qquad{} + (rs)^{-3}E_{k',(k + 1)'} \otimes E_{k'k'} - r^{-1}s^{-3}E_{k + 1,(k + 1)'} \otimes E_{k',k + 1} \\[0.2mm]
 \quad\qquad{} - r^{-3}s^{-1}E_{(k + 1)',(k + 1)'} \otimes E_{k',(k + 1)'}- rs^{-1}E_{k + 1,k+1} \otimes E_{k + 1,k} \\[0.2mm]
 \quad\qquad{} - r^{-1}sE_{k + 1,(k + 1)'} \otimes E_{(k + 1)',k} + r^{-3}s^{-1}E_{k'k} \otimes E_{k,(k +1)'}+ r^{-1}s^{-3}E_{k'k'} \otimes E_{k',(k + 1)'},
\\[0.2mm]
 R_{4}\Delta(f_{k}) = \bigl(r^{2} - s^{2} \bigr)(rs)^{-1} \bigl(r^{2(n - k - 1) + 1}s^{2(k + 1 - n) - 1} - 1 \bigr)E_{(k + 1)',(k + 1)'} \otimes E_{k + 1,k} \\[0.2mm]
 \hphantom{R_{4}\Delta(f_{k}) =}{}
 - \bigl(r^{2} - s^{2} \bigr)r^{-3}s^{-1} \bigl(r^{2(n - k) + 1}s^{2(k - n) - 1} - 1 \bigr)E_{k',(k + 1)'} \otimes E_{kk},
 \\[0.2mm]
 \Delta(f_{k})R_{4} = \bigl(r^{2} - s^{2} \bigr)(rs)^{-1} \bigl(r^{2(n - k) + 1}s^{2(k - n) - 1} - 1 \bigr)E_{k'k'} \otimes E_{k + 1,k} \\[0.2mm]
 \hphantom{\Delta(f_{k})R_{4} = }{}
 - \bigl(r^{2} - s^{2} \bigr)r^{-1}s^{-3} \bigl(r^{2(n - k) - 1}s^{2(k - n) + 1} - 1 \bigr)E_{k',(k + 1)'} \otimes E_{k + 1,k+1},
\\[0.2mm]
 \Delta(f_{k})R_{5} - R_{5}\Delta(f_{k}) \\[0.2mm]
 \qquad{}= \bigl(rs^{-1} - r^{-1}s \bigr)E_{k'k'} \otimes E_{k + 1,k} + \bigl(r^{-1}s - rs^{-1} \bigr)E_{k + 1,k + 1} \otimes E_{k + 1,k} \\[0.2mm]
 \quad\qquad{} + \bigl(r^{-1}s - rs^{-1} \bigr)E_{(k + 1)',(k + 1)'} \otimes E_{k + 1,k} - \bigl(r^{-3}s^{-1} - r^{-1}s^{-3} \bigr)E_{k'k'} \otimes E_{k',(k + 1)'} \\[0.2mm]
 \quad\qquad{} - \bigl(r^{-1}s^{3} - rs \bigr)E_{k + 1,k} \otimes E_{kk} + \bigl(r^{-1}s^{-1} - rs^{-3} \bigr)E_{k',(k + 1)'} \otimes E_{k + 1,k + 1} \\[0.2mm]
 \quad\qquad{} - \bigl(r^{-3}s - r^{-1}s^{-1} \bigr)E_{k',(k + 1)'} \otimes E_{kk} + \bigl(r^{-5}s^{-1} - r^{-3}s^{-3} \bigr)E_{k',(k + 1)'} \otimes E_{(k + 1)',(k + 1)'},
\\[0.2mm]
 R_{6}\Delta(f_{k}) = \bigl(r^{2} - s^{2} \bigr)(rs)^{-1}r^{2(n - k) - 1}s^{2(k - n) -1}E_{k',(k + 1)'} \otimes E_{kk} \\[0.2mm]
 \hphantom{R_{6}\Delta(f_{k}) =}{}
 - \bigl(r^{2} -s^{2} \bigr)(rs)^{-1}s^{2(k - n) + 1}r^{2(n - k) - 1}E_{(k + 1)',(k + 1)'} \otimes E_{k + 1,k} \\[0.2mm]
 \hphantom{R_{6}\Delta(f_{k}) = }{}
 - \bigl(rs^{-1} - r^{-1}s \bigr)E_{k + 1,(k + 1)'} \otimes E_{(k + 1)',k} + \bigl(r^{-1}s^{-3} - r^{-3}s^{-1} \bigr)E_{k'k} \otimes E_{k,(k + 1)'},
 \\[0.2mm]
 \Delta(f_{k})R_{6} = \bigl(r^{2} - s^{2} \bigr)(rs)^{-1}r^{2(n - k) - 1}s^{2(k - n) - 1}E_{k',(k + 1)'} \otimes E_{k + 1,k + 1} \\[0.2mm]
 \hphantom{\Delta(f_{k})R_{6} =}{}
 - \bigl(r^{2} - s^{2} \bigr)(rs)^{-1}r^{2(n - k) + 1}s^{2(k - n) - 1}E_{k'k'} \otimes E_{k + 1,k} \\[0.2mm]
 \hphantom{\Delta(f_{k})R_{6} =}{}
 + \bigl(rs^{-3} - r^{-1}s^{-1} \bigr)E_{k',k + 1} \otimes E_{k + 1,(k + 1)'} - \bigl(rs^{-3} - r^{-1}s^{-1} \bigr)E_{k + 1,k'} \otimes E_{k'k}.
\end{gather*}
In particular, we obtain
\begin{gather*}
 R_{4}\Delta(f_{k}) + R_{6}\Delta(f_{k}) \\
 \qquad{}= -\bigl(rs^{-1} - r^{-1}s\bigr)E_{(k + 1)',(k + 1)'} \otimes E_{k + 1,k} + \bigl(r^{-1}s^{-1} - r^{-3}s\bigr)E_{k',(k + 1)'} \otimes E_{kk} \\
 \qquad\quad\ {}{-} \bigl(rs^{-1} - r^{-1}s\bigr)E_{k + 1,(k + 1)'} \otimes E_{(k + 1)',k} + \bigl(r^{-1}s^{-3} - r^{-3}s^{-1}\bigr)E_{k'k} \otimes E_{k, (k + 1)'}
\end{gather*}
and
\begin{gather*}
 \Delta(f_{k})R_{4} + \Delta(f_{k})R_{6} \\
 \qquad{}=-\bigl(rs^{-1} - r^{-1}s\bigr)E_{k'k'} \otimes E_{k + 1,k} + \bigl(rs^{-3} - r^{-1}s^{-1}\bigr)E_{k',(k + 1)'} \otimes E_{k + 1,k + 1} \\
 \quad\qquad\ {} {+}\bigl(rs^{-3} - r^{-1}s^{-1}\bigr)E_{k',k + 1} \otimes E_{k + 1,(k + 1)'} - \bigl(rs^{-3} - r^{-1}s^{-1}\bigr)E_{k + 1,k'} \otimes E_{k'k}.
\end{gather*}

From the computations above, we finally get
\begin{gather}
 R_{3}\Delta(f_{k}) - \Delta(f_{k})R_{3} + R_{1}\Delta(f_{k}) + R_{2}\Delta(f_{k}) - \Delta(f_{k})R_{1} - \Delta(f_{k})R_{2} \label{eq:R123}\\
 \qquad{}= \bigl(r^{-1}s^{-3} - r^{-3}s^{-1}\bigr)E_{k'k'} \otimes E_{k',(k+1)'} + \bigl(rs^{-1} - r^{-1}s\bigr)E_{k + 1,(k + 1)'} \otimes E_{(k + 1)',k} \nonumber\\
 \quad\qquad{} + \bigl(r^{-3}s^{-1} - r^{-1}s^{-3}\bigr)E_{k'k} \otimes E_{k,(k + 1)'} + \bigl(rs - r^{-1}s^{3}\bigr)E_{k + 1,k} \otimes E_{kk} \nonumber\\
 \quad\qquad{} + \bigl(r^{-1}s^{-1} - rs^{-3}\bigr)E_{k + 1,k'} \otimes E_{k'k} + \bigl(rs^{-3} - r^{-1}s^{-1}\bigr)E_{k',k + 1} \otimes E_{k + 1,(k + 1)'} \nonumber\\
 \quad\qquad{} + \bigl(r^{-5}s^{-1} - r^{-3}s^{-3}\bigr)E_{k',(k + 1)'} \otimes E_{(k + 1)',(k + 1)'} + \bigl(r^{-1}s - rs^{-1}\bigr)E_{k + 1,k + 1} \otimes E_{k + 1,k}\nonumber
\end{gather}
and
\begin{gather}
 \Delta(f_{k})R_{5} - R_{5}\Delta(f_{k}) + \Delta(f_{k})R_{4} + \Delta(f_{k})R_{6} - R_{4}\Delta(f_{k}) - R_{6}\Delta(f_{k})\label{eq:R456}\\
 \qquad{}=\bigl(r^{-1}s - rs^{-1}\bigr)E_{k + 1,k + 1} \otimes E_{k + 1,k} +\bigl(r^{-1}s^{-3} - r^{-3}s^{-1}\bigr)E_{k'k'} \otimes E_{k',(k + 1)'} \nonumber\\
 \quad\qquad{} + \bigl(rs - r^{-1}s^{3}\bigr)E_{k + 1,k} \otimes E_{kk} + \bigl(r^{-5}s^{-1} - r^{-3}s^{-3}\bigr)E_{k',(k + 1)'} \otimes E_{(k + 1)',(k + 1)'}\nonumber \\
 \quad\qquad{} + \bigl(rs^{-3} - r^{-1}s^{-1}\bigr)E_{k',k + 1} \otimes E_{k + 1,(k + 1)'} + \bigl(r^{-1}s^{-1} - rs^{-3}\bigr)E_{k + 1,k'} \otimes E_{k'k}\nonumber \\
 \quad\qquad{} + \bigl(rs^{-1} - r^{-1}s\bigr)E_{k + 1,(k + 1)'} \otimes E_{(k + 1)',k} + \bigl(r^{-3}s^{-1} - r^{-1}s^{-3}\bigr)E_{k'k} \otimes E_{k,(k + 1)'}.\nonumber
\end{gather}
The right-hand sides of \eqref{eq:R123} and \eqref{eq:R456} are obviously equal, which implies $\hat{R}\Delta(f_{k}) = \Delta(f_{k})\hat{R}$ for $k < n$.

The computation for $k = n$ is similar. This completes the proof of the lemma.
\end{proof}

Now we are ready to prove the main results of this section, Theorems~\ref{thm:B_RMatrix}--\ref{thm:D_RMatrix}. We present full details only for the first one, since the other two are completely analogous (details are left to the interested~reader).

\begin{proof}[Proof of Theorem \ref{thm:B_RMatrix}]
According to Theorem~\ref{thm:universal-R} and Proposition~\ref{prop:struct}, the action of $\hat{R}_{VV}$ on the tensor product $V\otimes V$ is uniquely determined by the eigenvalues of the highest weight vectors $w_1$, $w_2$, $w_3$ from~\eqref{eq:B_hwv} with respect to~$\hat{R}_{VV}$. We shall now verify that these eigenvalues are precisely equal to $\lambda_1$, $\lambda_2$, $\lambda_3$ from Lemma~\ref{lem:B_eigen}, which thus completes the proof due to Lemma~\ref{lem:B_Intertwining}.

The eigenvalue $\widetilde{\lambda}_1$ of the $\hat{R}_{VV}$-action on $w_1$ is equal to $f(\varepsilon_{1},\varepsilon_{1})$. Since $\varepsilon_1=\alpha_1+\dots+\alpha_n$, \smash{$f(\varepsilon_1,\alpha_i)=\bigl(\omega'_i,\omega_{\varepsilon_1}\bigr)^{-1}$}, and a computation using \eqref{eq:B-pairing} yields \smash{$\bigl(\omega'_1,\omega_{\varepsilon_1}\bigr)=s^{-2}$}, \smash{$\bigl(\omega'_n,\omega_{\varepsilon_1}\bigr)=rs$}, \smash{$\bigl(\omega'_i,\omega_{\varepsilon_1}\bigr)=1$} if $1<i<n$, we thus obtain \smash{$\widetilde{\lambda}_1=r^{-1}s=\lambda_1$}.

The eigenvalue \smash{$\widetilde{\lambda}_2$} of the \smash{$\hat{R}_{VV}$}-action on $w_2$ equals the coefficient of $v_1\otimes v_2$ in \smash{$\hat{R}_{VV}(w_2)$}, and the latter appears only from applying $\widetilde{f}\circ \tau$ to the multiple of $v_2\otimes v_1$. Thus, we have \smash{$\widetilde{\lambda}_{2} = -rs^{-1}f(\varepsilon_{1},0) = -rs^{-1} = \lambda_{2}$} if $n = 1$. On the other hand, if $n > 1$, then \smash{$\widetilde{\lambda}_2=-r^2f(\varepsilon_1,\varepsilon_2)$}, and since $\varepsilon_2=\alpha_2+\dots+\alpha_n$, a similar calculation to the one above yields $\smash{\widetilde{\lambda}_2}=-r^2\cdot (rs)^{-1}=-rs^{-1}=\lambda_2$.

The eigenvalue $\widetilde{\lambda}_3$ of the $\hat{R}_{VV}$-action on $w_3$ equals the coefficient of $v_1\otimes v_{1'}$ in $\hat{R}_{VV}(w_3)$. The latter appears only from applying $\widetilde{f}\circ \tau$ to the multiple of $v_{1'}\otimes v_1$, thus $\widetilde{\lambda}_3=r^{2n-1}s^{1-2n} f(\varepsilon_1,-\varepsilon_1)$. As $f(\varepsilon_1,-\varepsilon_1)=f(\varepsilon_1,\varepsilon_1)^{-1}=rs^{-1}$, we thus get $\widetilde{\lambda}_3=r^{2n}s^{-2n}=\lambda_3$.
\end{proof}

\begin{Remark}\label{rem:easy-finite}
The above proofs of Theorems~\ref{thm:B_RMatrix}--\ref{thm:D_RMatrix} are quite elementary, but they require knowing the correct formulas for $\hat{R}$ in the first place. In the next section, we provide the conceptual origin of these formulas by factorizing them into an ordered product of ``local'' operators, one for each positive root of $\fg$.
\end{Remark}

\section{PBW bases, orthogonality, and factorization}\label{sec:R-factorized}

In this section, we present the factorization formulas for $\hat{R}$ from~\eqref{eq:B_RMatrix},~\eqref{eq:C_RMatrix} and~\eqref{eq:D_RMatrix}. In the absence of Lusztig's braid group action on $\uu$, one rather needs to use the combinatorial construction of orthogonal dual bases of $U^+_{r,s}(\fg)$ and $U^-_{r,s}(\fg)$, based on the combinatorics of standard Lyndon words, cf.\ \cite{K1,K2,L,Ro} (the details are presented in~\cite{MT1,MT2}).

\subsection{Standard Lyndon words}\label{ssec:Lyndon}

Let $I=\{1, 2, \dots, n\}$ be a finite ordered alphabet parametrizing the simple roots of $\fg$, and let~$I^*$ be the set of all finite length words in the alphabet $I$. For $u=[i_1 \dots i_k]\in I^*$, we define its \emph{length} by $|u|=k$. We introduce the \emph{lexicographical order} on $I^*$ in a standard way:
\[
 [i_1 \dots i_k] < [j_1 \dots j_l]
\quad \text{if}\ \begin{cases}
\text{$i_1=j_1, \dots, i_a=j_a, i_{a+1} < j_{a+1}$ for some $a \geq 0$}\\
\text{or $i_1=j_1, \dots, i_k=j_k$ and $k < l$}.
\end{cases}
\]
For a word $w = [i_1 \dots i_k]\in I^*$, the subwords
\begin{displaymath}
 w_{a|} = [i_1 \dots i_a] \qquad \text{and} \qquad w_{|a} = [i_{k-a+1} \dots i_k]
\end{displaymath}
with $0\leq a\leq k$ will be called a \emph{prefix} and a \emph{suffix} of $w$, respectively.
We call such a prefix or a~suffix \emph{proper} if $0<a<k$. We start with the following important definition:

\begin{Definition}
A word $w$ is \textit{Lyndon} if it is smaller than all of its proper suffixes:
\begin{displaymath}
 w < w_{|a} \qquad \mathrm{for\ all} \ 0<a<|w|.
\end{displaymath}
\end{Definition}

We recall the following two basic facts from the theory of Lyndon words:

\begin{Proposition}[{\cite[Proposition~5.1.3]{Lo}}]
Any Lyndon word $\ell$ has a factorization
\begin{equation}\label{eqn:st-factor}
 \ell = \ell_1 \ell_2
\end{equation}
defined by the property that $\ell_1$ is the longest proper prefix of $\ell$ which is also a Lyndon word.
Then, $\ell_2$ is also a Lyndon word.
\end{Proposition}

The factorization~\eqref{eqn:st-factor} is called a \textit{costandard factorization} of a Lyndon word.

\begin{Proposition}[{\cite[Proposition~5.1.5]{Lo}}]
Any word $w$ has a unique factorization
\begin{equation}\label{eqn:canon-factor}
 w = \ell_1 \dots \ell_k,
\end{equation}
where $\ell_1 \geq \dots \geq \ell_k$ are all Lyndon words.
\end{Proposition}

The factorization~\eqref{eqn:canon-factor} is called a \textit{canonical factorization}.

Let $\fn^+$ be a Lie subalgebra of $\fg$ generated by all $\{e_i\}_{i=1}^n$. The \textit{standard bracketing} of a Lyndon word~$\ell$ (with respect to the Lie algebra $\fn^+$) is defined inductively by the following procedure:
\begin{itemize}\itemsep=0pt\samepage
\item[$\bullet$]
$e_{[i]}=e_i\in \fn^+$ for $i \in I$,

\item[$\bullet$]
$e_{[\ell]} = [e_{[\ell_1]}, e_{[\ell_2]}]\in \fn^+$, where $\ell=\ell_1\ell_2$ is
the costandard factorization~\eqref{eqn:st-factor}.

\end{itemize}
The following definition is due to \cite{LR}.

\begin{Definition}
A Lyndon word $\ell$ is called \textit{standard Lyndon} if $e_{\ell}$ cannot be expressed as
a~linear combination of $e_m$ for various Lyndon words $m>\ell$.
\end{Definition}

The major importance of this definition is due to the following result.

\begin{Theorem}[{\cite{LR}}]
The set $\bigl\{e_{[\ell]}\mid \ell\mathrm{-standard\ Lyndon\ word}\bigr\}$ provides a basis of $\fn^+$.
\end{Theorem}

Due to a root space decomposition $\fn^+=\bigoplus_{\alpha\in \Phi^+} \fg_\alpha$ with all $\fg_\alpha$ being 1-dimensional, we get
\begin{displaymath}
 \ell \colon\ \Phi^+ \iso \{\slaws \},
\end{displaymath}
the so-called \emph{Lalonde--Ram} bijection, evoked in~\eqref{eqn:LR-bij-intro}. This bijection was described explicitly in~\cite{L}.

\begin{Proposition}[{\cite[Proposition~25]{L}}]\label{prop:Leclerc algorithm}
The bijection $\ell$ is inductively given as follows:
\begin{itemize}\itemsep=0pt
\item[$\bullet$]
for simple roots, we have $\ell(\alpha_i)=[i]$,

\item[$\bullet$]
for other positive roots, the value of $\ell(\alpha)$ is determined using Leclerc's algorithm:
\begin{equation*}
 \ell(\alpha) =
 \max\{ \ell(\gamma_1)\ell(\gamma_2) \mid
 \alpha=\gamma_1+\gamma_2 ,\, \gamma_1,\gamma_2\in \Phi^+ ,\, \ell(\gamma_1) < \ell(\gamma_2) \}.
\end{equation*}
\end{itemize}
\end{Proposition}

We shall also need one more important property of $\ell$. To the end, let us recall the following.

\begin{Definition}\label{def:convex}
A total order on the set of positive roots $\Phi^+$ is \textit{convex} if
\begin{displaymath}
 \alpha < \alpha+\beta < \beta
\end{displaymath}
for all $\alpha < \beta \in \Phi^+$ such that $\alpha+\beta$ is also a root.
\end{Definition}

The following result is~\cite[Proposition~28]{L}, where it is attributed to~\cite{Ro}
(see also~\cite[Proposition~2.34]{NT}).

\begin{Proposition}\label{prop:fin.convex}
Consider the order on $\Phi^+$ induced from the lexicographical order on standard Lyndon words:
\begin{equation}\label{eqn:induces}
 \alpha < \beta \quad \Longleftrightarrow \quad \ell(\alpha) < \ell(\beta) \ \mathrm{lexicographically}.
\end{equation}
This order is convex.
\end{Proposition}

Finally, recall that given any convex order on $\Phi^+$, a pair $(\alpha,\beta)$ of positive roots is called a~\textit{minimal pair} for $\gamma=\alpha+\beta\in \Phi^+$ if
\begin{equation}\label{eq:min_pair}
 \alpha < \beta \ \mathrm{and} \ \nexists
 \alpha < \alpha' < \gamma < \beta' < \beta \ \mathrm{such\ that} \ \alpha' + \beta' = \gamma.
\end{equation}
The following result goes back to~\cite{L, Ro} (cf.\ \cite[Proposition~2.38]{NT}).

\begin{Proposition}\label{prop:minimial-via-standard}
For any $\gamma\in \Phi^+$, consider the costandard factorization~\eqref{eqn:st-factor}, so that $\ell_1=\ell(\alpha)$ and $\ell_2=\ell(\beta)$. Then $(\alpha,\beta)$ is a minimal pair for $\gamma=\alpha+\beta$.
\end{Proposition}

\subsection{Convex orders for classical types}\label{ssec:ABCD-orders}

Let $\Phi$ be a root system of classical type, and choose the order of $I$ exactly as in Section~\ref{ssec:classical-2param}. Combining Propositions~\ref{prop:Leclerc algorithm} and~\ref{prop:fin.convex}, we obtain the following explicit convex orders on the sets~$\Phi^+$ of positive roots:
\begin{itemize}\itemsep=0pt
\item
\emph{Type $A_{n}$}:
\begin{equation}\label{eq:A_order}
 \alpha_{1} < \alpha_{1} + \alpha_{2} < \dots < \alpha_{1} + \dots + \alpha_{n} < \alpha_{2} <
 \dots < \alpha_{n-1} < \alpha_{n-1} + \alpha_{n} < \alpha_{n}.
\end{equation}

\item
\emph{Type $B_{n}$}:
\begin{align}
\alpha_{1} &{}< \alpha_{1} + \alpha_{2} < \dots < \alpha_{1} + \dots + \alpha_{n}\nonumber \\
 & {} < \alpha_{1} + \dots + \alpha_{n-1} + 2\alpha_{n} < \dots
 <\alpha_{1} + 2\alpha_{2} + \dots + 2\alpha_{n} \nonumber \\
 &{}<\alpha_{2} < \dots < \alpha_{n-1} < \alpha_{n-1} + \alpha_{n} < \alpha_{n-1} + 2\alpha_{n} < \alpha_{n}.\label{eq:B_order}
\end{align}

\item
\emph{Type $C_{n}$}:
\begin{align}
 \alpha_{1} &{}< \alpha_{1} + \alpha_{2} < \dots < \alpha_{1} + \dots + \alpha_{n-1} <
 2\alpha_{1} + \dots + 2\alpha_{n-1} + \alpha_{n} < \alpha_{1} + \dots + \alpha_{n} \nonumber\\
 &{}< \alpha_{1} + \dots + \alpha_{n-2} + 2\alpha_{n-1} + \alpha_{n} < \dots <
 \alpha_{1} + 2\alpha_{2} + \dots + 2\alpha_{n-1} + \alpha_{n} \nonumber\\
 &{}<\alpha_{2} < \dots < \alpha_{n-1} < 2\alpha_{n-1} + \alpha_{n} < \alpha_{n-1} + \alpha_{n} < \alpha_{n}.\label{eq:C_order}
\end{align}

\item
\emph{Type $D_{n}$}:
\begin{align}
 \alpha_{1} & < \alpha_{1} + \alpha_{2} < \dots < \alpha_{1} + \dots + \alpha_{n-2}+\alpha_{n-1} <
 \alpha_{1} +\dots + \alpha_{n-2} + \alpha_{n} \nonumber\\
 & < \alpha_{1} + \dots + \alpha_{n}< \alpha_{1} + \dots + \alpha_{n-3} + 2\alpha_{n-2} + \alpha_{n-1} + \alpha_{n} < \cdots\nonumber \\
 & < \alpha_{1} + 2\alpha_{2} + \dots + 2\alpha_{n-2} + \alpha_{n-1} + \alpha_{n} \nonumber\\
 &
 < \alpha_{2} < \dots < \alpha_{n-2}<\alpha_{n-2}+\alpha_{n-1}\nonumber\\
 & <\alpha_{n-2}+\alpha_{n}< \alpha_{n-2} + \alpha_{n-1} + \alpha_{n} <\alpha_{n-1}<\alpha_{n}.\label{eq:D_order}
\end{align}
\end{itemize}

\begin{Remark}\label{rem:tower-like}
An important feature of these convex orders on root systems of type $X_n$ (with $X=A,B,C,D$) is their telescopic structure, that is, erasing all roots containing $\alpha_1$ provides the order alike on the rank $1$ smaller root system of type $X_{n-1}$. This will significantly simplify our calculations in Section~\ref{ssec:R-Matrix-Computation}.
\end{Remark}

\begin{Remark}
It is a classical result, due to~\cite{P}, that the convex orders on $\Phi^+$ are in bijection with the reduced decompositions of the longest element $w_0$ of the Weyl group $W$ of the root system $\Phi$. In particular, the convex orders~\eqref{eq:A_order}--\eqref{eq:D_order} correspond respectively to the following reduced decompositions of $w_0$
(with $s_{i} = s_{\alpha_{i}}$ denoting the simple reflections):
\begin{align*}
 &w_{0} = (s_{1}s_{2}\dots s_{n})(s_{1}s_{2}\dots s_{n-1})\dots (s_{1}s_{2})(s_{1}), \\
 &w_{0} = (s_{1}\dots s_{n - 1}s_{n}s_{n-1}\dots s_{1})(s_{2}\dots s_{n-1}s_{n}s_{n-1}\dots s_{2})\dots (s_{n-1}s_{n}s_{n-1})(s_{n}), \\
 &w_{0} = (s_{1}\dots s_{n-1}s_{n}s_{n-1}\dots s_{1})(s_{2}\dots s_{n-1}s_{n}s_{n-1}\dots s_{2})\dots (s_{n-1}s_{n}s_{n-1})(s_{n}), \\
 &w_{0} =
 \begin{cases}
 (s_{1}\dots s_{n-1}s_{n}s_{n-2}\dots s_{1})(s_{2}\dots s_{n-2}s_{n}s_{n-1}\dots s_{2})\dots& \\
 \qquad\times{}(s_{n-2}s_{n-1}s_{n}s_{n-2})(s_{n}s_{n-1}) &\text{if $n$ is odd}, \\
 (s_{1}\dots s_{n-1}s_{n}s_{n-2}\dots s_{1})(s_{2}\dots s_{n-2}s_{n}s_{n-1}\dots s_{2})\dots &\\
 \qquad\times{}(s_{n-2}s_{n}s_{n-1}s_{n-2})(s_{n-1}s_{n}) &\text{if $n$ is even}.
 \end{cases}
\end{align*}
\end{Remark}

Let us set up convenient notation for the positive roots in each type, and identify the minimal pairs arising through the costandard factorization of the corresponding standard Lyndon words, see Proposition~\ref{prop:minimial-via-standard}.

\emph{Type $A_n$}.
Let $\gamma_{ij} = \alpha_{i} + \dots + \alpha_{j}$ for $1 \le i \le j \le n$, so that
$\Phi^+=\{\gamma_{ij}\}_{i\leq j}$. For $\gamma=\gamma_{ij}$ with $i < j$, the minimal pair arising through the costandard factorization of $\ell(\gamma)$ is $(\alpha,\beta)=(\gamma_{i,j - 1},\alpha_{j})$.

\emph{Type $B_n$}.
Let $\gamma_{ij} = \alpha_{i} + \dots + \alpha_{j}$ for $1 \le i \le j \le n$, and let $\beta_{ij} = \alpha_{i} + \dots + \alpha_{j-1} + 2\alpha_{j} + \dots + 2\alpha_{n}$ for $1 \le i < j \le n$. Let us now indicate the minimal pairs arising through the costandard factorization of $\ell(\gamma)$:
\begin{itemize}\itemsep=0pt
\item
for the roots $\gamma=\gamma_{ij}$ with $i < j$, the minimal pair is $(\alpha,\beta)=(\gamma_{i,j - 1},\alpha_{j})$;

\item
for the roots $\gamma=\beta_{in}$ with $1\leq i<n$, the minimal pair is $(\alpha,\beta)=(\gamma_{in},\alpha_{n})$;

\item
for the roots $\gamma=\beta_{ij}$ with $i < j < n$, the minimal pair is $(\alpha,\beta)=(\beta_{i,j + 1},\alpha_{j})$.
\end{itemize}

\emph{Type $C_n$}.
Let $\gamma_{ij} = \alpha_{i} + \dots + \alpha_{j}$ for $1 \le i \le j \le n$, and let $\beta_{ij} = \alpha_{i} + \dots + \alpha_{j - 1} + 2\alpha_{j} + \dots + 2\alpha_{n-1} + \alpha_{n}$ for $1 \le i \le j<n$. Let us now indicate the minimal pairs arising through the costandard factorization of $\ell(\gamma)$:
\begin{itemize}\itemsep=0pt
\item
for the roots $\gamma=\gamma_{ij}$ with $i<j$, the minimal pair is $(\alpha,\beta)=(\gamma_{i,j - 1},\alpha_{j})$;

\item
for the roots $\gamma=\beta_{ij}$ with $i<j<n-1$, the minimal pair is $(\alpha,\beta)=(\beta_{i,j + 1},\alpha_{j})$;

\item
for the roots $\gamma=\beta_{ii}$ with $1\leq i<n$, the minimal pair is $(\alpha,\beta)=(\gamma_{i,n - 1},\gamma_{in})$;

\item for the roots $\gamma=\beta_{i,n-1}$ with $1\leq i<n-1$, the minimal pair is $(\alpha,\beta)=(\gamma_{in},\alpha_{n-1})$.

\end{itemize}

\emph{Type $D_n$}.
Following~\eqref{eq:D-system}, let $\gamma_{ij} = \varepsilon_{i} - \varepsilon_{j + 1}$ for $1 \le i \le j <n$, and let $\beta_{ij} = \varepsilon_{i} + \varepsilon_{j}$ for $1 \le i < j \le n$. Then $\gamma_{ij} = \alpha_{i} + \dots + \alpha_{j}$ for $1 \le i \leq j < n$, $\beta_{in} = \alpha_{i} + \dots + \alpha_{n-2} + \alpha_{n}$, $\beta_{i,n -1} = \alpha_{i} + \dots + \alpha_{n}$, and $\beta_{ij} = \alpha_{i} + \dots + \alpha_{j - 1} + 2\alpha_{j} + \dots + 2\alpha_{n-2} + \alpha_{n-1} + \alpha_{n}$ for $j < n - 1$. Let us indicate the minimal pairs arising through the costandard factorization of $\ell(\gamma)$:
\begin{itemize}\itemsep=0pt
\item
for the roots $\gamma=\gamma_{ij}$ with $i<j$, the minimal pair is $(\alpha,\beta)=(\gamma_{i,j - 1},\alpha_{j})$;

\item
for the roots $\gamma=\beta_{in}$ with $i\leq n-2$, the minimal pair is $(\alpha,\beta)=(\gamma_{i,n - 2},\alpha_{n})$;

\item
for the roots $\gamma=\beta_{ij}$ with $i < j < n$, the minimal pair is $(\alpha,\beta)=(\beta_{i,j + 1},\alpha_{j})$.
\end{itemize}

\subsection{Root vectors and the PBW theorem}\label{ssec:PBW theorem and bases}

The following construction of (quantum) \textit{root vectors} $e_{\gamma}$, $f_{\gamma}$ goes back to~\cite{K1,K2,L,Ro} in the one-parameter setup, to~\cite{CHW} in the super setup, and to~\cite{BH,BKL,HW1,HW2} in the two-parameter setup. For $\gamma=\alpha_i\in \Pi$, we set
$ e_{\alpha_i}=e_i$, $f_{\alpha_i}=f_i$.
By induction on the height of a root, for any $\gamma\in \Phi^+\setminus \Pi$, define (cf.\ \eqref{eqn:quantum-root-vectors})
\begin{equation}\label{eq:root_vector}
 e_{\gamma} = e_{\alpha}e_{\beta} - \bigl(\omega_{\beta}',\omega_{\alpha}\bigr)e_{\beta}e_{\alpha} , \qquad
 f_{\gamma} = f_{\beta}f_{\alpha} - \bigl(\omega_{\alpha}',\omega_{\beta}\bigr)^{-1}f_{\alpha}f_{\beta},
\end{equation}
where the minimal pair $(\alpha,\beta)$ for $\gamma$ corresponds to the costandard factorization of $\ell(\gamma)$, see Proposition~\ref{prop:minimial-via-standard}.

We now state a two-parameter version of a classical result of~\cite{K1,K2,L,Ro} (which was also adapted to the super case in~\cite{CHW}). A detailed proof of Theorem~\ref{thm:PBW-general}\,(a),~(b), and~\eqref{eq:root-vector-pairing} is presented in~\cite[Theorem~7.1]{MT1}, and the special case of~\eqref{eq:pairing_recursion} specified in~\eqref{eq:A-parity-constants}--\eqref{eq:D-parity-constants} is presented in~\cite[Theorem~7.2]{MT1}. Both arguments make crucial use of an embedding of $U_{r,s}^{+}$ into an appropriate quantum shuffle algebra. A proof of the recursive formula~\eqref{eq:pairing_recursion} in general is given in~\cite{MT2}, along with an alternative proof of the rest of the theorem.

\begin{Theorem}\label{thm:PBW-general}\samepage \quad
\begin{enumerate}\itemsep=0pt
\item[$(a)$] The ordered products
\begin{displaymath}
 \Biggl\{ \overset{\longleftarrow}{\underset{\gamma \in \Phi^{+}}{\prod}} e_{\gamma}^{m_{\gamma}} \mid m_{\gamma} \ge 0 \Biggr\}\qquad \text{and}\qquad \Biggl\{\overset{\longleftarrow}{\underset{\gamma \in \Phi^{+}}{\prod}}f_{\gamma}^{m_{\gamma}} \mid m_{\gamma} \ge 0\Biggr\}
\end{displaymath}
are bases for $U_{r,s}^{+}(\mathfrak{g})$ and $U_{r,s}^{-}(\mathfrak{g})$, respectively. Here and below, the arrow $\leftarrow$ over the product signs refers to the total order~\eqref{eqn:induces} on $\Phi^+$.

\item[$(b)$] The Hopf pairing~\eqref{eq:Hopf-parity} is orthogonal with respect to these bases. More explicitly, we have
\begin{displaymath}
 \Biggl(\overset{\longleftarrow}{\underset{\gamma \in \Phi^{+}}{\prod}} f_{\gamma}^{n_{\gamma}}, \overset{\longleftarrow}{\underset{\gamma \in \Phi^{+}}{\prod}} e_{\gamma}^{m_{\gamma}}\Biggr)=
 \prod_{\gamma\in \Phi^+} \bigl(\delta_{n_\gamma,m_\gamma} \bigl(f_{\gamma}^{m_\gamma},e^{m_\gamma}_\gamma\bigr)\bigr).
\end{displaymath}

\item[$(c)$] For each $\gamma \in \Phi^{+}$ and $m \ge 0$, we have $($cf.\ notation~\eqref{eq:rs_i}$)$
\begin{equation}\label{eq:root-vector-pairing}
 \bigl(f_{\gamma}^{m},e_{\gamma}^{m}\bigr) = s_{\gamma}^{-m(m - 1)/2} (f_{\gamma},e_{\gamma})^{m}[m]_{r_{\gamma},s_{\gamma}}!
\end{equation}
Moreover, we have the following recursive formula for the pairing $(f_\gamma,e_\gamma)$:
\begin{gather}
 (f_{\gamma},e_{\gamma}) = (f_{\alpha},e_{\alpha}) (f_{\beta},e_{\beta})\nonumber \\
 \hphantom{(f_{\gamma},e_{\gamma}) =}{}\times
 \frac{r_{\gamma}\bigl(\omega_{\alpha}',\omega_{\beta}\bigr)^{-1}(rs)^{-p_{\beta,\alpha}}[p_{\beta,\alpha} + 1]_{r,s}^{2}(r_{\alpha} - s_{\alpha})(r_{\beta} - s_{\beta})}
 {r_{\alpha}r_{\beta}(s_{\gamma} - r_{\gamma})},\label{eq:pairing_recursion}
\end{gather}
where $(\alpha,\beta)$ is the minimal pair of~Proposition~$\ref{prop:minimial-via-standard}$ corresponding to the costandard factorization of $\ell(\gamma)$,~and
\begin{displaymath}
 p_{\beta,\alpha} = \max \{k \ge 0\mid \beta - k\alpha \in \Phi\}.
\end{displaymath}
\end{enumerate}
\end{Theorem}

As an immediate corollary, we obtain the following \textit{factorization formula}.

\begin{Theorem}\label{cor:Theta-factorization} The operator $\Theta$ of~\eqref{eq:Theta} can be factorized as follows:
\begin{equation}\label{eq:Theta-factorized}
 \Theta=\overset{\longleftarrow}{\underset{\gamma \in \Phi^{+}}{\prod}}
 \biggl( \sum_{m\geq 0} \frac{1}{\bigl(f^m_{\gamma},e^m_{\gamma}\bigr)} f^m_\gamma\otimes e^m_\gamma \biggr),
\end{equation}
with $e_\gamma$, $f_\gamma$ defined in~\eqref{eq:root_vector} and $\bigl(f_{\gamma}^m,e_{\gamma}^m\bigr)$ evaluated in \eqref{eq:root-vector-pairing}--\eqref{eq:pairing_recursion}.
\end{Theorem}

\begin{Remark}
In fact, $\Theta$ can be expressed in an even more compact form by using \eqref{eq:Theta-factorized} in~conjunction with \textit{$(r,s)$-exponential}
\begin{displaymath}
 \exp_{r,s}(z) = \sum_{m \ge 0}s^{m(m-1)/2} \frac{z^{m}}{[m]_{r,s}!}.
\end{displaymath}
Indeed, for each $\gamma \in \Phi^{+}$, we have
\begin{equation}\label{eq:local-exp}
 \sum_{m \ge 0} \frac{1}{\bigl(f_{\gamma}^{m},e_{\gamma}^{m}\bigr)}f_{\gamma}^{m} \otimes e_{\gamma}^{m} =
 \sum_{m \ge 0}\frac{s^{m(m - 1)/2}_\gamma}{(f_{\gamma},e_{\gamma})^{m}[m]_{r_{\gamma},s_{\gamma}}!}f_{\gamma}^{m} \otimes e_{\gamma}^{m}
 = \exp_{r_{\gamma},s_{\gamma}} \biggl(\frac{f_{\gamma} \otimes e_{\gamma}}{(f_{\gamma},e_{\gamma})}\biggr).
\end{equation}
Thus, the factorization formula~\eqref{eq:Theta-factorized} simplifies as follows:
\begin{displaymath}
 \Theta = \overset{\longleftarrow}{\underset{\gamma \in \Phi^{+}}{\prod}}
 \exp_{r_{\gamma},s_{\gamma}}\biggl(\frac{f_{\gamma} \otimes e_{\gamma}}{(f_{\gamma},e_{\gamma})}\biggr).
\end{displaymath}
\end{Remark}

We conclude this subsection by presenting explicit formulas for $\bigl(f_{\gamma}^{m},e_{\gamma}^{m}\bigr)$ in all classical types:
\begin{itemize}\itemsep=0pt
\item
\emph{Type $A_{n}$}
\begin{equation}\label{eq:A-parity-constants}
 \bigl(f_{\gamma_{ij}}^{m},e_{\gamma_{ij}}^{m}\bigr) = (-1)^{m}s^{-m(m-1)/2}\frac{[m]_{r,s}!}{(r - s)^{m}}
 \qquad \mathrm{for} \ 1 \le i \le j \le n.
\end{equation}

\item
\emph{Type $B_{n}$}
\begin{gather}
 \bigl(f_{\gamma_{ij}}^{m},e_{\gamma_{ij}}^{m}\bigr) = (-1)^{m}s^{-m(m - 1)}\frac{[m]_{r^{2},s^{2}}!}{\bigl(r^{2} - s^{2}\bigr)^{m}}
 \qquad \mathrm{for} \ 1 \le i \le j < n, \nonumber\\
 \bigl(f_{\gamma_{in}}^{m},e_{\gamma_{in}}^{m}\bigr) = (-1)^{m}s^{-m(m-1)/2}\frac{[m]_{r,s}!}{(r - s)^{m}}
 \qquad \mathrm{for} \ 1 \le i \le n,\label{eq:B-parity-constants}\\
 \bigl(f_{\beta_{ij}}^{m},e_{\beta_{ij}}^{m}\bigr) =
 (-1)^{m}[2]_{r,s}^{2m}(rs)^{-2m(n - j)}s^{-m(m - 1)}\frac{[m]_{r^{2},s^{2}}!}{\bigl(r^{2} - s^{2}\bigr)^{m}}
 \qquad \mathrm{for} \ 1 \le i < j \le n.\nonumber
\end{gather}

\item
\emph{Type $C_{n}$}
\begin{gather}
 \bigl(f_{\gamma_{ij}}^{m},e_{\gamma_{ij}}^{m}\bigr) = (-1)^{m}s^{-m(m - 1)/2}\frac{[m]_{r,s}!}{(r - s)^{m}}
 \qquad \mathrm{for} \ 1 \le i \le j \le n \ \mathrm{with} \ (i,j) \neq (n,n),\nonumber \\
 \bigl(f_{\gamma_{nn}}^{m},e_{\gamma_{nn}}^{m}\bigr) = (-1)^{m}s^{-m(m - 1)}\frac{[m]_{r^{2},s^{2}}!}{\bigl(r^{2} - s^{2}\bigr)^{m}},\nonumber \\
 \bigl(f_{\beta_{ii}}^{m},e_{\beta_{ii}}^{m}\bigr) =
 (-1)^{m}[2]_{r,s}^{2m}s^{-m(m - 1)}\frac{[m]_{r^{2},s^{2}}!}{\bigl(r^{2} - s^{2}\bigr)^{m}}
 \qquad \mathrm{for} \ 1 \le i < n,\nonumber\\
 \bigl(f_{\beta_{ij}}^{m},e_{\beta_{ij}}^{m}\bigr) = (-1)^{m}(rs)^{-m(n - j)}s^{-m(m - 1)/2}\frac{[m]_{r,s}!}{(r - s)^{m}}
 \qquad \mathrm{for} \ 1 \le i < j < n.\label{eq:C-parity-constants}
\end{gather}

\item
\emph{Type $D_{n}$}
\begin{gather}
 \bigl(f_{\gamma_{ij}}^{m},e_{\gamma_{ij}}^{m}\bigr) = (-1)^{m}s^{-m(m - 1)/2}\frac{[m]_{r,s}!}{(r - s)^{m}}
 \qquad \mathrm{for} \ 1 \le i \le j < n, \nonumber\\
 \bigl(f_{\beta_{ij}}^{m},e_{\beta_{ij}}^{m}\bigr) = (-1)^{m}(rs)^{-m(n - j)}s^{-m(m - 1)/2}\frac{[m]_{r,s}!}{(r - s)^{m}}
 \qquad \mathrm{for} \ 1 \le i < j \le n.\label{eq:D-parity-constants}
\end{gather}
\end{itemize}

\subsection[R-matrix computation]{$\boldsymbol{R}$-matrix computation}\label{ssec:R-Matrix-Computation}

We shall now use the factorization formula~\eqref{eq:Theta-factorized} to compute $\Theta$ and $\hat{R}_{VV}$ for all classical types.
Throughout this subsection, we will use the more convenient notation (cf.\ \eqref{eq:local-exp})
\begin{equation}\label{eq:local-Theta}
 \Theta_{\gamma} = \sum_{m \ge 0}\frac{1}{\bigl(f_{\gamma}^{m},e_{\gamma}^{m}\bigr)}f_{\gamma}^{m} \otimes e_{\gamma}^{m}
 \qquad \mathrm{for\ any} \ \gamma \in \Phi^{+},
\end{equation}
so that equation \eqref{eq:Theta-factorized} becomes
\begin{displaymath}
 \Theta = \overset{\longleftarrow}{\underset{\gamma \in \Phi^{+}}{\prod}}\Theta_{\gamma}.
\end{displaymath}

\subsubsection[Factorized formula in type A\_n]{Factorized formula in type $\boldsymbol{A_{n}}$}

We start by evaluating the action of $\{e_\gamma,f_\gamma\}_{\gamma\in \Phi^+}$ on the $U_{r,s}(\ssl_{n+1})$-representation $V$ from Proposition~\ref{prp:A_rep}.

\begin{Lemma}
$\rho(e_{\gamma_{ij}}) = E_{i,j + 1}$ and $\rho(f_{\gamma_{ij}}) = E_{j + 1,i}$ for any $1\leq i\leq j\leq n$.
\end{Lemma}

\begin{proof}
The proof is straightforward and proceeds by induction on the height $|\gamma_{ij}|=j-i+1$, where we use
\begin{gather*}
 e_{\gamma_{ij}}=e_{\gamma_{i,j-1}}e_{j} - \bigl(\omega_j',\omega_i \dots \omega_{j-1}\bigr) e_{j}e_{\gamma_{i,j-1}},\\
 f_{\gamma_{ij}} = f_{j}f_{\gamma_{i,j-1}}  - \bigl(\omega'_i\dots \omega'_{j-1},\omega_j\bigr)^{-1} f_{\gamma_{i,j-1}}f_{j}
\end{gather*}
and the explicit $U_{r,s}(\ssl_{n+1})$-action on $V$, cf.\ \eqref{eq:root_vector} and Proposition~\ref{prp:A_rep}.
\end{proof}

Clearly, $\rho(e_{\gamma_{ij}})^{2} = \rho(f_{\gamma_{ij}})^{2} = 0$ for all $i \le j$, so that
\begin{displaymath}
 \Theta_{\gamma_{ij}} = 1 \otimes 1 + (s - r)E_{j + 1,i} \otimes E_{i,j + 1}
\end{displaymath}
on the $U_{r,s}(\ssl_{n+1})$-representation $V \otimes V$. To evaluate $\Theta$, let us first set
\begin{displaymath}
 \Theta_{i} = \Theta_{\gamma_{in}}\Theta_{\gamma_{i,n-1}}\cdots \Theta_{\gamma_{i,i+1}}\Theta_{\alpha_{i}},
\end{displaymath}
so that
\begin{displaymath}
 \Theta = \Theta_{n}\Theta_{n-1}\cdots \Theta_{1}.
\end{displaymath}
Since $E_{j + 1,i} \otimes E_{i,j + 1} \cdot E_{k + 1,i} \otimes E_{i,k + 1} = 0$ for any $i\leq k < j$, we thus obtain
\begin{displaymath}
 \Theta_{i} = 1 \otimes 1 + (s - r)\sum_{j = i + 1}^{n} E_{j + 1,i} \otimes E_{i,j + 1}.
\end{displaymath}
Moreover, since $E_{k + 1,j} \otimes E_{j,k + 1} \cdot E_{\ell + 1,i} \otimes E_{i,\ell + 1} = 0$ for any $k \geq j > i$ and $\ell \geq i$, we get
\begin{equation}\label{eq:A-Theta}
 \Theta = \Theta_{n}\cdots \Theta_{1} = 1 \otimes 1 + (s - r)\sum_{1 \le i < j \le n + 1} E_{ji} \otimes E_{ij}.
\end{equation}

For type $A_{n}$, the function $f$ satisfying~\eqref{eq:f} can be chosen as follows (see~\cite[Lemma~4.4]{BW1}):
\begin{displaymath}
 f(\varepsilon_{i},\varepsilon_{j}) =
 \begin{cases}
 s^{-1} & \text{if}\ i < j, \\
 1 & \text{if}\ i = j, \\
 r & \text{if}\ i > j.
 \end{cases}
\end{displaymath}
Combining this formula with the formula~\eqref{eq:A-Theta} and the flip map $\tau$, we recover the explicit formula~\eqref{eq:A_RMatrix} for $\hat{R}_{VV}=\Theta \circ \widetilde{f} \circ \tau\colon V \otimes V \to V \otimes V$, thus providing an alternative proof of Theorem~\ref{thm:A_RMatrix}, and giving a factorization of the $R$-matrix from~\cite{BW2} into ``local'' operators, parametrized by the set of positive roots.

\subsubsection[Factorized formula in type B\_n]{Factorized formula in type $\boldsymbol{B_{n}}$}

We start by evaluating the action of $\{e_\gamma,f_\gamma\}_{\gamma\in \Phi^+}$ on the $U_{r,s}(\sso_{2n+1})$-module $V$ from Proposition~\ref{prp:B_rep}.

\begin{Lemma}\label{lem:B-root-action}\quad\samepage
\begin{enumerate}\itemsep=0pt
\item[$(a)$] For $i\leq j\leq n$, we have
\begin{displaymath}
\begin{split}
 & \rho(e_{\gamma_{ij}}) = E_{i,j+1} - s^{2(j - i)}E_{(j + 1)',i'}, \\
 & \rho(f_{\gamma_{ij}}) =
 \begin{cases}
 E_{j + 1,i} - s^{2(i - j)}(rs)^{-2}E_{i',(j + 1)'} & \text{if}\ j < n ,\\
 \bigl(r^{-1} + s^{-1}\bigr)(E_{n + 1,i} - s^{2(i - n)}E_{i',n+1})& \text{if}\ j = n.
 \end{cases}
\end{split}
\end{displaymath}

\item[$(b)$] For $i<j\leq n$, we have
\begin{displaymath}
\begin{split}
 & \rho(e_{\beta_{ij}}) = (-1)^{n+ 1-j}\bigl(E_{ij'} - s^{2(n - i)}r^{-2(n - j)}(rs)E_{ji'}\bigr), \\
 & \rho(f_{\beta_{ij}}) = (-1)^{n+1-j}\bigl(r^{-1} + s^{-1}\bigr)^2\cdot s^{2(j-n)}
 \bigl(r^{2(j - n)}E_{j'i} - s^{2(i - n)}(rs)E_{i'j}\bigr).
\end{split}
\end{displaymath}
\end{enumerate}
\end{Lemma}

\begin{proof}
The proof is straightforward and proceeds by an increasing induction on $j$ (from $j=i$ up to $j=n$) for the roots $\gamma=\gamma_{ij}$, and then by a descending induction on $j$ (from $j=n$ till $j=i+1$) for the roots $\gamma=\beta_{ij}$. Here, we use the explicit $U_{r,s}(\sso_{2n+1})$-action on $V$ from Proposition~\ref{prp:B_rep}, the explicit list of minimal pairs $(\alpha,\beta)$ as specified in Section~\ref{ssec:ABCD-orders}, and finally the inductive construction~\eqref{eq:root_vector}.
\end{proof}

According to Lemma~\ref{lem:B-root-action}\,(a), we have $\rho(e_{\gamma_{ij}})^{2} = \rho(f_{\gamma_{ij}})^{2} = 0$ for $1 \le i \le j < n$, so that
\begin{gather}
 \Theta_{\gamma_{ij}} = 1 \otimes 1\nonumber\\
 \hphantom{\Theta_{\gamma_{ij}} =}{}
 + \bigl(s^{2} - r^{2}\bigr)\bigl(E_{j + 1,i} - s^{2(i - j)}(rs)^{-2}E_{i',(j + 1)'}\bigr) \otimes \bigl(E_{i,j+1} - s^{2(j - i)}E_{(j + 1)',i'}\bigr)\label{eq:Theta_gamma_<n}
\end{gather}
for any $1 \le i \le j < n$. In contrast, for $j = n$, we have
\begin{gather*}
 \rho(e_{\gamma_{in}})^{2} = -s^{2(n - i)}E_{ii'}, \qquad
 \rho(f_{\gamma_{in}})^{2} = -s^{2(i - n)}(r + s)^{2}(rs)^{-2}E_{i'i}, \\
 \rho(e_{\gamma_{in}})^{3} = \rho(f_{\gamma_{in}})^{3} = 0.
\end{gather*}
Therefore, we obtain
\begin{gather}
 \Theta_{\gamma_{in}} = 1 \otimes 1
 + \bigl(s^{2} - r^{2}\bigr)(rs)^{-1} \bigl( E_{n + 1,i} - s^{2(i - n)}E_{i',n+1} \bigr) \otimes \bigl(E_{i,n+1} - s^{2(n - i)}E_{n + 1,i'}\bigr) \nonumber\\
 \hphantom{\Theta_{\gamma_{in}} =}{}
 + r^{-2}s^{-1}\bigl(r^{2} - s^{2}\bigr)(r - s) E_{i'i} \otimes E_{ii'}.\label{eq:Theta_gamma_n}
\end{gather}

According to Lemma~\ref{lem:B-root-action}\,(b), we also have $\rho(e_{\beta_{ij}})^{2} = \rho(f_{\beta_{ij}})^{2} = 0$ for all $i < j\leq n$, so that
\begin{gather*}
 \Theta_{\beta_{ij}} = 1 \otimes 1\\
 \hphantom{\Theta_{\beta_{ij}} =}{}
 +\bigl(s^{2} - r^{2}\bigr)(rs)^{-2} \bigl( r^{2(j - n)}E_{j'i} - s^{2(i - n)}(rs)E_{i'j} \bigr) \otimes \bigl(r^{2(n - j)}E_{ij'} - s^{2(n - i)}(rs)E_{ji'}\bigr).
\end{gather*}

To evaluate $\Theta$, let us first set
\begin{displaymath}
 \Theta_{i}'' = \Theta_{\beta_{i,i + 1}}\Theta_{\beta_{i,i + 2}}\cdots \Theta_{\beta_{in}}, \qquad
 \Theta_{i}' = \Theta_{\gamma_{in}}\Theta_{\gamma_{i,n-1}} \cdots \Theta_{\gamma_{i,i+1}}\Theta_{\alpha_{i}}, \qquad
 \Theta_{i} = \Theta_{i}''\Theta_{i}' ,
\end{displaymath}
so that
\begin{displaymath}
 \Theta = \Theta_{n}\Theta_{n-1}\cdots \Theta_{1}.
\end{displaymath}

We will now show by induction on $n - k$ that{\samepage
\begin{gather}
 \Theta^{(k)} := \Theta_{n}\Theta_{n-1}\cdots \Theta_{k}\nonumber\\
 \hphantom{\Theta^{(k)}{:}}{} =
 1 \otimes 1 + c\sum_{i = k}^{n}\bigl(r^{-1}s - r^{2(n - i)}s^{2(i-n)}\bigr)E_{i'i} \otimes E_{ii'} \nonumber\\
 \hphantom{\Theta^{(k)} :=}{}
 + c \sum_{k \le i < j \le n}\bigl(rsE_{ji} \otimes E_{ij} - r^{2(j - i)-1}sE_{ji} \otimes E_{j'i'} - r^{-1}s^{2(i - j) + 1}E_{i'j'} \otimes E_{ij}\nonumber\\
 \hphantom{\Theta^{(k)} :=+ c \sum_{k \le i < j \le n}\bigl(}{}
 + (rs)^{-1}E_{i'j'} \otimes E_{j'i'}\bigr)\nonumber \\
 \hphantom{\Theta^{(k)} :=}{}
 + c\sum_{i = k}^{n}\bigl(E_{n + 1,i} \otimes E_{i,n+1} - r^{2(n - i)}E_{n + 1,i} \otimes E_{n + 1,i'} - s^{2(i-n)}E_{i',n+1} \otimes E_{i,n+1}\nonumber\\
 \hphantom{\Theta^{(k)} :=+ c\sum_{i = k}^{n}\bigl(}{}
 + E_{i',n+1} \otimes E_{n + 1,i'} \bigr) \nonumber\\
 \hphantom{\Theta^{(k)} :=}{}
 + c \sum_{k \le i < j \le n}\bigl( (rs)^{-1}E_{j'i} \otimes E_{ij'} - r^{2(n-i)}s^{2(j - n)}E_{j'i} \otimes E_{ji'}\bigr) \nonumber\\
 \hphantom{\Theta^{(k)} :=}{}
 +c\sum_{k \le i < j \le n}\bigl( - r^{2(n - j)}s^{2(i - n)}E_{i'j} \otimes E_{ij'} + rsE_{i'j} \otimes E_{ji'} \bigr),\label{eq:Theta_k}
\end{gather}
where $c = \bigl(s^{2} - r^{2}\bigr)(rs)^{-1}$. For $k=1$, this provides the desired formula for $\Theta=\Theta^{(1)}$.}

Let us start by computing a single $\Theta_{i}$. First, it is easy to see from \eqref{eq:Theta_gamma_<n} and \eqref{eq:Theta_gamma_n} that the only non-zero products in $\Theta_{\gamma_{in}}\Theta_{\gamma_{i,n-1}}\cdots \Theta_{\alpha_{i}}$ involve $1 \otimes 1$, and thus
\begin{gather*}
 \Theta_{i}' = 1 \otimes 1 + c\bigl(r^{-1}s-1\bigr)E_{i'i} \otimes E_{ii'} \\
 \hphantom{\Theta_{i}' =}{}
 + c\sum_{j = i + 1}^{n}\bigl(rsE_{ji} \otimes E_{ij} - rs^{2(j - i) - 1}E_{ji} \otimes E_{j'i'} - r^{-1}s^{2(i-j) + 1}E_{i'j'} \otimes E_{ij}
 \\
 \hphantom{\Theta_{i}' =+ c\sum_{j = i + 1}^{n}\bigl(}{}
 + (rs)^{-1}E_{i'j'} \otimes E_{j'i'}\bigr) \\
 \hphantom{\Theta_{i}' =}{}
 + c\bigl(E_{n + 1,i} \otimes E_{i,n+1} - s^{2(n - i)}E_{n + 1,i} \otimes E_{n+1, i'} - s^{2(i - n)}E_{i',n + 1} \otimes E_{i,n + 1} \\
 \hphantom{\Theta_{i}' =+ c\bigl(}{}
 + E_{i',n + 1} \otimes E_{n + 1,i'} \bigr).
\end{gather*}
For similar reasons, we have
\begin{align*}
 \Theta_{i}'' &{}= 1 \otimes 1 + c\sum_{j = i + 1}^{n}\bigl((rs)^{-1}E_{j'i} \otimes E_{ij'} - r^{2(j - n)}s^{2(n - i)}E_{j'i} \otimes E_{ji'}\bigr) \\
 &\quad{} +c\sum_{j = i + 1}^{n}\bigl(- r^{2(n - j)}s^{2(i - n)}E_{i'j} \otimes E_{ij'} + rsE_{i'j} \otimes E_{ji'} \bigr).
\end{align*}
On the other hand, when computing the product $\Theta_{i}''\Theta_{i}'$, there are some non-zero products not involving $1 \otimes 1$, namely $cr^{2(n - j)}s^{2(i - n)}E_{i'j} \otimes E_{ij'} \cdot crs^{2(j - i) - 1}E_{ji} \otimes E_{j'i'}$ for $i < j \le n$. The overall contribution of those equals
\begin{displaymath}
\begin{split}
 c^{2}r^{2n + 1}s^{-2n - 1}\Biggl(\sum_{j = i + 1}^{n}r^{-2j}s^{2j} \Biggr)E_{i'i} \otimes E_{ii'} &= c^{2}rs^{2(i - n) + 1}[n - i]_{r^{2},s^{2}} E_{i'i} \otimes E_{ii'} \\
 & = c\bigl(1-r^{2(n - i)}s^{2(i - n)} \bigr)E_{i'i} \otimes E_{ii'},
\end{split}
\end{displaymath}
and therefore we have
\begin{gather}
 \Theta_{i} = \Theta_{i}''\Theta_{i}'
 = 1 \otimes 1 + c\bigl(r^{-1}s-r^{2(n - i)}s^{2(i - n)}\bigr) E_{i'i} \otimes E_{ii'} \nonumber\\
 \hphantom{\Theta_{i} = \Theta_{i}''\Theta_{i}'=}{}
 + c\sum_{j = i + 1}^{n}\bigl(rsE_{ji} \otimes E_{ij} - rs^{2(j - i) - 1}E_{ji} \otimes E_{j'i'} - r^{-1}s^{2(i-j) + 1}E_{i'j'} \otimes E_{ij}\nonumber\\
 \hphantom{\Theta_{i} =\Theta_{i}''\Theta_{i}'= + c\sum_{j = i + 1}^{n}\bigl(}{}
 + (rs)^{-1}E_{i'j'} \otimes E_{j'i'}\bigr) \nonumber\\
 \hphantom{\Theta_{i} =\Theta_{i}''\Theta_{i}'= }{}
 + c\bigl(E_{n + 1,i} \otimes E_{i,n+1} - s^{2(n - i)}E_{n + 1,i} \otimes E_{n+1, i'} - s^{2(i - n)}E_{i',n + 1} \otimes E_{i,n + 1}\nonumber\\
 \hphantom{\Theta_{i} =\Theta_{i}''\Theta_{i}'= + c\bigl(}{}
 + E_{i',n + 1} \otimes E_{n + 1,i'} \bigr) \nonumber\\
 \hphantom{\Theta_{i} =\Theta_{i}''\Theta_{i}' =}{}
 + c\sum_{j = i + 1}^{n}\bigl((rs)^{-1}E_{j'i} \otimes E_{ij'} - r^{2(j - n)}s^{2(n - i)}E_{j'i} \otimes E_{ji'} \bigr)\nonumber\\
 \hphantom{\Theta_{i} =\Theta_{i}''\Theta_{i}'= }{}
 +c\sum_{j = i + 1}^{n}\bigl(- r^{2(n - j)}s^{2(i - n)}E_{i'j} \otimes E_{ij'} + rsE_{i'j} \otimes E_{ji'} \bigr).\label{eq:B-individual-Theta}
\end{gather}

In particular, for $i=n$ we get
\begin{align*}
 \Theta_{n}={}&1 \otimes 1 + c\bigl(r^{-1}s-1\bigr) E_{i'i} \otimes E_{ii'}
 \\
 &{}+ c\bigl(E_{n + 1,n} \otimes E_{n,n+1} - E_{n + 1,n} \otimes E_{n+1, n'} - E_{n',n + 1} \otimes E_{n,n + 1} + E_{n',n + 1} \otimes E_{n + 1,n'} \bigr),
\end{align*}
which agrees with the claimed formula~\eqref{eq:Theta_k} when $k = n$, thus establishing the base of our induction.

Let
us now prove the step of induction in~\eqref{eq:Theta_k}. It suffices to treat the $k=1$ case, due to the telescopic structure of the action on $V$ and order on $\Phi^+$, see Remark~\ref{rem:tower-like}. It thus remains to evaluate $\Theta^{(1)}=\Theta^{(2)}\cdot \Theta_1$, where $\Theta^{(2)}$ is given by~\eqref{eq:Theta_k} and $\Theta_1$ was just evaluated above. In~addition to the terms of this product that involve $1 \otimes 1$, we get the following extra summands:\looseness=1
\begin{gather}
 \Biggl(-c\sum_{i = 2}^{n} s^{2(i -n)}E_{i',n + 1} \otimes E_{i,n + 1}\Biggr) \bigl(-cs^{2(n - 1)}E_{n + 1,1} \otimes E_{n + 1,1'} \bigr) \nonumber\\
 \qquad{}= \bigl(r^{2} - s^{2}\bigr)^{2}\sum_{i = 2}^{n}r^{-2}s^{2i - 4}E_{i'1} \otimes E_{i1'},\label{eq:term_a}
\\
 \Biggl(-c\sum_{2 \le i < j \le n}r^{2(n - i)}s^{2(j - n)}E_{j'i} \otimes E_{ji'}\Biggr) \Biggl(-c\sum_{i = 2}^{n}rs^{2i - 3}E_{i1} \otimes E_{i'1'}\Biggr) \nonumber\\
 \qquad{}=
 c^{2}\sum_{i = 2}^{n-1}\sum_{j = i + 1}^{n} r^{2(n - i) + 1}s^{2(j + i - n -1) - 1}E_{j'1} \otimes E_{j1'} \nonumber\\
 \qquad{}=
 c^{2}\sum_{j = 3}^{n}\Biggl( \sum_{i = 2}^{j - 1}r^{2(n - i) + 1}s^{2(j + i - n -1) - 1}\Biggr)E_{j'1} \otimes E_{j1'} \nonumber\\
 \qquad{}=
 \bigl(r^{2} - s^{2}\bigr)^{2}\sum_{j = 2}^{n}r^{2(n - j) + 1}s^{2(j - n) - 1}[j - 2]_{r^{2},s^{2}}E_{j'1} \otimes E_{j1'},\label{eq:term_b}
\\
 \Biggl(-c\sum_{2 \le i < j \le n}r^{2(n - j)}s^{2(i - n)}E_{i'j} \otimes E_{ij'}\Biggr)\Biggl(-c\sum_{j = 2}^{n}rs^{2j - 3}E_{j1} \otimes E_{j'1'}\Biggr) \nonumber\\
 \qquad{}=
 \bigl(r^{2} - s^{2}\bigr)^{2}\sum_{i = 2}^{n-1}\Biggl( \sum_{j = i + 1}^{n}r^{2(n - j)-1}s^{2(i + j - n - 2) - 1}\Biggr) E_{i'1} \otimes E_{i1'} \nonumber\\
 \qquad{}=
 \bigl(r^{2} -s^{2}\bigr)^{2}\sum_{i = 2}^{n}r^{-1}s^{4i - 2n - 3}[n - i]_{r^{2},s^{2}}E_{i'1} \otimes E_{i1'},\label{eq:term_c}
\\
 \Biggl(c\sum_{i = 2}^{n}\bigl(r^{-1}s - r^{2(n - i)}s^{2(i - n)}\bigr)E_{i'i} \otimes E_{ii'} \Biggr)\Biggl(-c\sum_{i = 2}^{n}rs^{2(i - 1) -1}E_{i1} \otimes E_{i'1'}\Biggr) \nonumber\\
 \qquad{}=
 \bigl(r^{2} - s^{2}\bigr)^{2}\sum_{i = 2}^{n}\bigl(r^{2(n - i) -1}s^{4i - 2n -5} -r^{-2}s^{2(i - 2)}\bigr)E_{i'1} \otimes E_{i1'},\label{eq:term_d}
\\
 \Biggl(-c\sum_{2 \le i < j \le n}r^{-1}s^{2(i - j) + 1}E_{i'j'} \otimes E_{ij} \Biggr)\Biggl(-c\sum_{j = 2}^{n}r^{2(j - n)}s^{2(n - 1)}E_{j'1} \otimes E_{j1'} \Biggr) \nonumber\\
 \qquad{}=
 \bigl(r^{2} - s^{2}\bigr)^{2}\sum_{i = 2}^{n-1}\Biggl(\sum_{j = i + 1}^{n}r^{2(j - n - 1) - 1}s^{2(n + i - j - 1) - 1}\Biggr)E_{i'1} \otimes E_{i1'} \nonumber\\
 \qquad{}=
 \bigl(r^{2} - s^{2}\bigr)^{2}\sum_{i = 2}^{n}r^{2(i - n) - 1}s^{2i - 3}[n - i]_{r^{2},s^{2}}E_{i'1} \otimes E_{i1'},\label{eq:term_e}
\\
 \Biggl(-c\sum_{2 \le i < j \le n}r^{2(j - i) - 1}sE_{ji} \otimes E_{j'i'} \Biggr)\Biggl(-c\sum_{i = 2}^{n}rs^{2i - 3}E_{i1} \otimes E_{i'1'}\Biggr) \nonumber\\
 \qquad{}=
 \bigl(r^{2} - s^{2}\bigr)^{2}\sum_{i = 2}^{n}\sum_{j = i + 1}^{n}r^{2(j - i - 1)}s^{2(i -2)}E_{j1} \otimes E_{j'1'} \nonumber\\
 \qquad{}=
 \bigl(r^{2} - s^{2}\bigr)^{2}\sum_{j = 3}^{n}\Biggl( \sum_{i = 2}^{j - 1}r^{2(j - i - 1)}s^{2(i - 2)}\Biggr) E_{j1} \otimes E_{j'1'} \nonumber\\
 \qquad{}=
 \bigl(r^{2} - s^{2}\bigr)^{2}\sum_{j = 2}^{n}[j - 2]_{r^{2},s^{2}}E_{j1} \otimes E_{j'1'},\label{eq:term_f}
\\
 \Biggl(-c\sum_{i = 2}^{n}r^{2(n - i)}E_{n + 1,i} \otimes E_{n + 1,i'} \Biggr) \Biggl(-c\sum_{i = 2}^{n}rs^{2i - 3}E_{i1} \otimes E_{i'1'}\Biggr) \nonumber\\
 \qquad{}=
 \bigl(r^{2} - s^{2}\bigr)^{2}\Biggl(\sum_{i = 2}^{n}r^{2(n - i) - 1}s^{2(i - 2) - 1}\Biggr)E_{n+1,1} \otimes E_{n + 1,1'} \nonumber\\
 \qquad{}=
 \bigl(r^{2} - s^{2}\bigr)^{2}(rs)^{-1}[n - 1]_{r^{2},s^{2}}E_{n + 1,1} \otimes E_{n + 1,1'}.\label{eq:term_g}
\end{gather}

Thus, the overall contribution of the terms $\{E_{j1} \otimes E_{j'1'}\}_{1<j<n+1}$ into $\Theta$ equals
\begin{gather*}
 \sum_{j = 2}^{n} \bigl(r^{2} - s^{2}\bigr)s^{2(j-2)} E_{j1} \otimes E_{j'1'} +
 \sum_{j = 3}^{n}\bigl(r^{2} - s^{2}\bigr)^{2}[j - 2]_{r^{2},s^{2}} E_{j1} \otimes E_{j'1'}\\
 \qquad{} = \bigl(r^{2} - s^{2}\bigr)\sum_{j = 2}^{n}r^{2(j - 2)}E_{j1} \otimes E_{j'1'},
\end{gather*}
where the first summand arises from $1 \otimes 1 \cdot \Theta_{1}$ and the second from~\eqref{eq:term_f}.
Likewise, the overall coefficient of $E_{n + 1,1} \otimes E_{n + 1,1'}$ in $\Theta$ is
\begin{displaymath}
 \bigl(r^{2} - s^{2}\bigr)(rs)^{-1}s^{2(n-1)} + \bigl(r^{2} - s^{2}\bigr)^{2}(rs)^{-1}[n-1]_{r^{2},s^{2}}=
 \bigl(r^{2} - s^{2}\bigr)(rs)^{-1}r^{2(n - 1)},
\end{displaymath}
where the first summand arises from $1 \otimes 1 \cdot \Theta_{1}$ and the second from~\eqref{eq:term_g}.
Finally, the remaining terms~\eqref{eq:term_a}--\eqref{eq:term_e} contribute to the coefficients of $\{E_{j'1} \otimes E_{j1'}\}_{2\leq j\leq n}$. Combining these with the corresponding terms from $1 \otimes 1 \cdot \Theta_{1}$, one eventually arrives at
\begin{displaymath}
 \bigl(r^{2} - s^{2}\bigr)(rs)^{-1} \sum_{i = 2}^{n} r^{2(n - 1)}s^{2(i-n)} E_{i'1} \otimes E_{i1'}.
\end{displaymath}

Thus,~\eqref{eq:term_a}--\eqref{eq:term_g} and their counterparts from $1 \otimes 1 \cdot \Theta_{1}$ match the corresponding part in the right-hand side of~\eqref{eq:Theta_k} for $k=1$. On the other hand, it is easy to see that the remaining terms arising from $\Theta^{(2)} \cdot 1 \otimes 1$ and $1 \otimes 1 \cdot \Theta_{1}$ exactly match with the remaining terms in the right-hand side of~\eqref{eq:Theta_k} for $k=1$. This completes our proof of the induction step, thus establishing~\eqref{eq:Theta_k}.

To derive the formula for $\hat{R}_{VV}$ it only remains to compute the values of $f$ from~\eqref{eq:f} on the weights of $V$. In accordance with~\eqref{eq:f}, we have
\begin{displaymath}
 f(\lambda,\mu) = \bigl(\omega_{\mu}',\omega_{\lambda}\bigr)^{-1},
\end{displaymath}
where we extend the Hopf pairing to the weight lattice as in~\eqref{eq:weight-pairing}.
From the formulas~\eqref{eq:B-pairing}, the equality $\varepsilon_k=\alpha_{k}+\dots+\alpha_n$  for $1\leq k\leq n$, and the basic properties of the Hopf pairing, we derive
\begin{equation}\label{eq:B-f-1}
 f(\varepsilon_{i},\varepsilon_{j}) = f(-\varepsilon_{i},-\varepsilon_{j}) =
 \begin{cases}
 (rs)^{-1} & \text{if}\ i < j ,\\
 r^{-1}s & \text{if}\ i = j, \\
 rs & \text{if}\ i > j,
 \end{cases}
\end{equation}
while the remaining values are then determined by
\begin{equation}\label{eq:B-f-2}
 f(\varepsilon_{i},-\varepsilon_{j}) = f(-\varepsilon_{i},\varepsilon_{j}) = f(\varepsilon_{i},\varepsilon_{j})^{-1},\qquad
 f(0,0) = f(0, \pm \varepsilon_{i}) = f(\pm \varepsilon_{i},0) = 1.
\end{equation}
We also note that if we set $\varepsilon_{n+1} = 0$ and $\varepsilon_{i'} = -\varepsilon_{i}$ for all $1 \le i \le n$, then we obtain
\begin{equation}\label{eq:f_aij_B}
 f(\varepsilon_{i},\varepsilon_{j}) = a_{ij} \qquad \mathrm{for\ all} \ i \neq j,j',
\end{equation}
where $a_{ij}$ are given by \eqref{eq:B_coeffs}. Combining the formulas~\eqref{eq:B-f-1} and~\eqref{eq:B-f-2} with the formula~\eqref{eq:Theta_k} for $k=1$ and the flip map $\tau$, we recover the explicit formula~\eqref{eq:B_RMatrix} for $\smash{\hat{R}_{VV}=\Theta \circ \widetilde{f}} \circ \tau\colon V \otimes V \to V \otimes V$, thus providing an alternative proof of Theorem~\ref{thm:B_RMatrix}.

\subsubsection[Factorized formula in type C\_n]{Factorized formula in type $\boldsymbol{C_{n}}$}

The calculations in this case are very similar to those for type $B_{n}$, so we highlight only the main points.
As before, we start by deriving explicit formulas for the action of $\{e_\gamma,f_\gamma\}_{\gamma\in \Phi^+}$ on the $U_{r,s}(\ssp_{2n})$-representation $V$ from Proposition~\ref{prp:C_rep}.

\begin{Lemma}\label{lem:C-allroots}\quad
\begin{enumerate}\itemsep=0pt
\item[$(a)$] For $i \le j \le n$, we have
\begin{displaymath}
\begin{split}
 \rho(e_{\gamma_{ij}}) &=
 \begin{cases}
 E_{i,j + 1} - s^{j - i}E_{(j + 1)',i'} & \text{if}\ j < n, \\
 E_{in'} + s^{n + 1 - i}E_{ni'} & \text{if}\ j = n,\\
 \end{cases} \\
 \rho(f_{\gamma_{ij}}) &=
 \begin{cases}
 E_{j + 1,i} - s^{i - j}(rs)^{-1}E_{i',(j+1)'} & \text{if}\ j < n ,\\
 (rs)^{-1}E_{n'i} + s^{i - n - 1}E_{i'n} & \text{if}\ j = n.
 \end{cases}
\end{split}
\end{displaymath}

\item[$(b)$] For $1 \le i < n$, we have
\begin{displaymath}
\begin{split}
 & \rho(e_{\beta_{ii}}) = s^{n-i}(r +s) E_{ii'}, \\
 & \rho(f_{\beta_{ii}}) = s^{i - n}\bigl(r^{-1} + s^{-1}\bigr)E_{i'i}.
\end{split}
\end{displaymath}

\item[$(c)$] For $1 \le i < j < n$, we have
\begin{displaymath}
\begin{split}
 & \rho(e_{\beta_{ij}}) = (-1)^{n-j}\bigl(E_{ij'} + r^{j-n}s^{n + 1 - i}E_{ji'}\bigr) , \\
 & \rho(f_{\beta_{ij}}) = (-s)^{j-n}\bigl(r^{j - n}(rs)^{-1}E_{j'i} + s^{i - n- 1}E_{i'j}\bigr) .
\end{split}
\end{displaymath}
\end{enumerate}
\end{Lemma}

Similarly to the above treatment of $B_{n}$-type, we define
\begin{displaymath}
 \Theta_{i} = \Theta_{\beta_{i,i + 1}}\cdots \Theta_{\beta_{i,n-1}}
 \Theta_{\gamma_{in}} \Theta_{\beta_{ii}}\Theta_{\gamma_{i,n-1}}\cdots \Theta_{\gamma_{i,i+1}}\Theta_{\alpha_{i}}
 \qquad \mathrm{for} \ 1\leq i\leq n,
\end{displaymath}
where each factor is evaluated through~\eqref{eq:local-Theta} by using Lemma~\ref{lem:C-allroots} and formulas~\eqref{eq:C-parity-constants}.
The following formula is derived completely analogously to~\eqref{eq:B-individual-Theta}:
\begin{gather}
 \Theta_{i} = 1 \otimes 1 + c\bigl(r^{n - i + 1}s^{i - n} + s\bigr) E_{i'i} \otimes E_{ii'} \nonumber\\
 \hphantom{\Theta_{i} =}{}
 + c\sum_{j = i + 1}^{n}\bigl(rsE_{ji} \otimes E_{ij} - rs^{j - i}E_{ji} \otimes E_{j'i'} - s^{i - j + 1}E_{i'j'} \otimes E_{ij} + E_{i'j'} \otimes E_{j'i'}\bigr)\nonumber \\
 \hphantom{\Theta_{i} =}{}
 + c\sum_{j = i + 1}^{n}\bigl(E_{j'i} \otimes E_{ij'} + r^{j - n}s^{n + 1 - i}E_{j'i} \otimes E_{ji'} + r^{n - j + 1}s^{i - n}E_{i'j} \otimes E_{ij'}\nonumber\\
 \hphantom{\Theta_{i} =+ c\sum_{j = i + 1}^{n}\bigl(}{}
 + rsE_{i'j} \otimes E_{ji'}\bigr),\label{eq:C-individual-Theta}
\end{gather}
where $c = (s - r)(rs)^{-1}$. Consider the following family of operators $\Theta^{(k)}=\Theta_n \Theta_{n-1} \cdots \Theta_{k}$, so that $\Theta=\Theta^{(1)}$. We claim that $\Theta^{(k)}$ is explicitly given by the following formula:
\begin{gather}
 \Theta^{(k)} = 1 \otimes 1 + c\sum_{i = k}^{n}\bigl(r^{n - i + 1}s^{i - n} + s\bigr)E_{i'i} \otimes E_{ii'} \nonumber \\
 \hphantom{\Theta^{(k)} =}{}
 + c\sum_{k \le i < j \le n}\bigl(rsE_{ji} \otimes E_{ij} - r^{j - i}sE_{ji} \otimes E_{j'i'} - s^{i - j + 1}E_{i'j'} \otimes E_{ij} + E_{i'j'} \otimes E_{j'i'}\bigr) \nonumber\\
 \hphantom{\Theta^{(k)} =}{}
 + c\sum_{k \le i < j \le n}\bigl(E_{j'i} \otimes E_{ij'} + r^{n + 1 - i}s^{j - n}E_{j'i} \otimes E_{ji'} + r^{n - j + 1}s^{i - n}E_{i'j} \otimes E_{ij'} \nonumber\\
 \hphantom{\Theta^{(k)} = + c\sum_{k \le i < j \le n}\bigl(}{}
 + rsE_{i'j} \otimes E_{ji'}\bigr).\label{eq:C-Theta_k}
\end{gather}
The proof proceeds by an induction on $n-k$, with the base case $k=n$ following from~\eqref{eq:C-individual-Theta}. As per the step of induction, we note that when opening brackets in $\Theta^{(1)}=\Theta^{(2)}\Theta_1$, besides for the summands where one of the terms is $1\otimes 1$, we get the following additional terms:
\begin{gather*}
 (s - r)\sum_{j = 2}^{n}r^{n - 1}s^{j - n - 1}\bigl(1 - (r^{-1}s)^{j - 2}\bigr)E_{j'1} \otimes E_{j1'},
\\
 (s - r)\sum_{j = 2}^{n}r^{n - j}s^{2j - 2 - n}\bigl(1 - r^{j - n}s^{n- j}\bigr)E_{j'1} \otimes E_{j1'},
\\
 -(r - s)^{2}\sum_{j =2}^{n}s^{j - 3}\bigl(r^{n - j}s^{j - n} + r^{-1}s\bigr)E_{j'1} \otimes E_{j1'},
\\
 (s - r)\sum_{j = 2}^{n}r^{j - n - 1}s^{n - 1}\bigl(r^{n - j}s^{j - n} - 1\bigr)E_{j'1} \otimes E_{j1'},
\\
 (s - r)\sum_{j = 2}^{n}r^{j - 2}\bigl(r^{2 - j}s^{j - 2} - 1\bigr)E_{j1} \otimes E_{j'1'}.
\end{gather*}
Combining these summands with the appropriate terms from $1 \otimes 1 \cdot \Theta_{1}$ and $\Theta^{(2)}\cdot 1 \otimes 1$ matches precisely the right-hand side of~\eqref{eq:C-Theta_k} for $k=1$, thus providing the formula for $\Theta$.

From~\eqref{eq:C-pairing}, the equality \smash{$\varepsilon_k=\alpha_{k}+\dots+\alpha_{n-1}+\frac{1}{2}\alpha_n$}  for $1\leq k\leq n$, and basic properties of the Hopf pairing, we derive
\begin{equation}\label{eq:C-wght-pair-1}
 f(\varepsilon_{i},\varepsilon_{j}) = f(-\varepsilon_{i},-\varepsilon_{j}) =
 \begin{cases}
 (rs)^{-1/2} & \text{if} \ i < j, \\
 \bigl(r^{-1}s\bigr)^{1/2} & \text{if}\ i = j, \\
 (rs)^{1/2} & \text{if}\ i > j,
 \end{cases}
\end{equation}
while the remaining values are then determined by
\begin{equation}\label{eq:C-wght-pair-2}
 f(\varepsilon_{i},-\varepsilon_{j}) = f(-\varepsilon_{i},\varepsilon_{j}) = f(\varepsilon_{i},\varepsilon_{j})^{-1}.
\end{equation}
Furthermore, setting $\varepsilon_{i'} = -\varepsilon_{i}$ for all $1 \le i \le n$, we obtain
\begin{equation}\label{eq:f_aij_CD}
 f(\varepsilon_{i},\varepsilon_{j}) = a_{ij} \qquad \mathrm{for\ all} \ i \neq j,j',
\end{equation}
where $a_{ij}$ are given by \eqref{eq:C_coeffs}. Combining the formulas~\eqref{eq:C-wght-pair-1} and~\eqref{eq:C-wght-pair-2} with the formula~\eqref{eq:C-Theta_k} for $k=1$ and the flip map $\tau$, we recover the explicit formula~\eqref{eq:C_RMatrix} for $\smash{\hat{R}_{VV}=\Theta \circ \widetilde{f}} \circ \tau\colon V \otimes V \to V \otimes V$, thus providing an alternative proof of Theorem~\ref{thm:C_RMatrix}.

\subsubsection[Factorized formula in type D\_n]{Factorized formula in type $\boldsymbol{D_{n}}$}

As in the previous types, we start by deriving explicit formulas for the action of $\{e_\gamma,f_\gamma\}_{\gamma\in \Phi^+}$ on the $U_{r,s}(\sso_{2n})$-representation $V$ from Proposition~\ref{prp:D_rep}.

\begin{Lemma}\quad\samepage
\begin{enumerate}
\item[$(a)$] For $1 \le i \le j < n$, we have
\begin{displaymath}
\begin{split}
 &\rho(e_{\gamma_{ij}})= E_{i,j + 1} - s^{j - i}E_{(j + 1)',i'},\\
 &\rho(f_{\gamma_{ij}})= E_{j + 1,i} - s^{i - j}(rs)^{-1}E_{i',(j + 1)'}.
\end{split}
\end{displaymath}

\item[$(b)$] For $1 \le i < j \le n$, we have
\begin{displaymath}
\begin{split}
 &\rho(e_{\beta_{ij}})= (-1)^{n - j}\bigl((rs)^{-1}E_{ij'} - r^{j - n}s^{n - i - 1} E_{ji'}\bigr), \\
 &\rho(f_{\beta_{ij}})= (-1)^{n-j}\bigl((rs)^{j - n}E_{j'i} - s^{i + j + 1 - 2n}E_{i'j}\bigr).
\end{split}
\end{displaymath}
\end{enumerate}
\end{Lemma}

Similarly to the previous types, we define
\begin{displaymath}
 \Theta_{i} = \Theta_{\beta_{i,i+1}}\cdots \Theta_{\beta_{i,n-1}}\Theta_{\beta_{in}}
 \Theta_{\gamma_{i,n-1}}\Theta_{\gamma_{i,n-2}}\cdots \Theta_{\gamma_{i,i+1}}\Theta_{\alpha_{i}}
 \qquad \mathrm{for} \ 1 \le i < n,
\end{displaymath}
 (in particular, $\Theta_{n-1}=\Theta_{\alpha_n}\Theta_{\alpha_{n-1}}$) and derive the following counterpart of the formulas~\eqref{eq:B-individual-Theta} and~\eqref{eq:C-individual-Theta}:
\begin{gather}
 \Theta_{i} = 1 \otimes 1 + c(s - r^{n - i}s^{i + 1 - n})E_{i'i} \otimes E_{ii'}\label{eq:D-individual-Theta} \\
 \hphantom{\Theta_{i} =}{}
 \!+ c\!\sum_{j = i + 1}^{n}\!(rsE_{ji} \otimes E_{ij} - rs^{j - i}E_{ji} \otimes E_{j'i'} - s^{i - j + 1}E_{i'j'} \otimes E_{ij} + E_{i'j'} \otimes E_{j'i'})\nonumber \\
 \hphantom{\Theta_{i} =}{}
 \!+ c\!\sum_{j = i + 1}^{n}\!(E_{j'i} \otimes E_{ij'} - r^{j + 1 - n}s^{n - i}E_{j'i} \otimes E_{ji'} - r^{n - j}s^{i + 1 - n}E_{i'j} \otimes E_{ij'} + rsE_{i'j} \otimes E_{ji'}),\nonumber
\end{gather}
where $c = (s - r)(rs)^{-1}$. Consider the following family of operators $\Theta^{(k)}=\Theta_{n-1} \cdots \Theta_{k}$, so that $\Theta=\Theta^{(1)}$. We claim that $\Theta^{(k)}$ is explicitly given by the following formula:
{\samepage\begin{gather}
 \Theta^{(k)} = 1 \otimes 1 + c \sum_{i = k}^{n}\bigl(s - r^{n - i}s^{i + 1 - n}\bigr)E_{i'i} \otimes E_{ii'}\nonumber \\
 \hphantom{\Theta^{(k)} =}{}
 + c\sum_{k \le i < j \le n}\bigl(rsE_{ji} \otimes E_{ij} - r^{j - i}s E_{ji} \otimes E_{j'i'} - s^{i - j + 1}E_{i'j'} \otimes E_{ij} + E_{i'j'} \otimes E_{j'i'}\bigr)\nonumber \\
 \hphantom{\Theta^{(k)} =}{}
 + c\sum_{k \le i < j \le n}\bigl(E_{j'i} \otimes E_{ij'} - r^{n - i}s^{j + 1- n}E_{j'i} \otimes E_{ji'} - r^{n - j}s^{i+ 1 - n}E_{i'j} \otimes E_{ij'}\nonumber\\
 \hphantom{\Theta^{(k)} =+ c\sum_{k \le i < j \le n}\bigl(}{}
 + rsE_{i'j} \otimes E_{ji'}\bigr).\label{eq:D-Theta_k}
\end{gather}
The proof proceeds by an induction on $n-k - 1$, with the base case $k=n-1$ following from~\eqref{eq:D-individual-Theta}.}
As per the step of induction, we note that when opening brackets in $\Theta^{(1)}=\Theta^{(2)}\Theta_1$, besides for the summands where one of the terms is $1\otimes 1$, we get the following additional terms:\looseness=1
\begin{gather*}
 \sum_{j = 2}^{n}(r - s)r^{n - 2}s^{j - n}\bigl(1 - r^{2-j}s^{j - 2}\bigr)E_{j'1} \otimes E_{j1'},
\\
 \sum_{j = 2}^{n}(r-s)r^{n- j- 1}s^{2j - 1 - n}\bigl(1 - r^{j-n}s^{n - j}\bigr)E_{j'1} \otimes E_{j1'},
\\
 \sum_{j = 2}^{n}(r - s)^{2}r^{-1}s^{j - 2}\bigl(r^{n - j}s^{j - n} - 1\bigr)E_{j'1} \otimes E_{j1'},
\\
 \sum_{j = 2}^{n}(r - s)r^{j - n}s^{n - 2}\bigl(r^{n-j}s^{j - n} - 1\bigr)E_{j'1} \otimes E_{j1'},
\\
 \sum_{j = 2}^{n}(r - s)r^{j - 2}\bigl(1 - r^{2 - j}s^{j - 2}\bigr)E_{j1} \otimes E_{j'1'}.
\end{gather*}
Combining these summands with the appropriate terms from $1 \otimes 1 \cdot \Theta_{1}$ and $\Theta^{(2)}\cdot 1\otimes 1$ matches precisely the right-hand side of~\eqref{eq:D-Theta_k} for $k=1$, thus providing the formula for $\Theta$.

Since $f$ is again given by~\eqref{eq:C-wght-pair-1}, we recover the explicit formula~\eqref{eq:D_RMatrix} for
\[
\hat{R}_{VV}=\Theta \circ \widetilde{f} \circ \tau\colon \ V \otimes V \to V \otimes V,
\]
 thus providing an alternative proof of Theorem~\ref{thm:D_RMatrix}.

\section[R-matrices with a spectral parameter]{$\boldsymbol{R}$-matrices with a spectral parameter}\label{sec:affine-R}

In this section, we generalize the previous constructions to the affine setup of two-parameter quantum affine algebras, which were introduced and studied in the literature case-by-case, see~\cite{GHZ, HRZ,HZ2,HZ1}.

\subsection{Two-parameter quantum affine algebras and evaluation modules}\label{ssec:2-param affine}

We start with a uniform definition of $\UU$. To this end, let $\theta\in \Phi^+$ be the highest root of~$\fg$, so that the new simple root $\alpha_0$ of $\wfg$ is given by $\alpha_0=\delta-\theta$. Define the matrix of structural constants \smash{$\Omega=(\Omega_{ij})_{i,j=0}^n$} by
\begin{alignat*}{3}
 &\Omega_{ij}=\bigl(\omega'_i,\omega_j\bigr), \qquad&& \Omega_{0j}=\bigl(\omega'_{-\theta},\omega_j\bigr), & \\
 &\Omega_{i0}=\bigl(\omega'_i, \omega_{-\theta}\bigr), \qquad&& \Omega_{00}=\bigl(\omega'_{-\theta},\omega_{-\theta}\bigr)
 \qquad \mathrm{for\ all} \ 1\leq i,j\leq n.&
\end{alignat*}
We also recall that the Cartan matrix \smash{$\hat{C}=(c_{ij})_{i,j = 0}^{n}$} for $\wfg$ is given by
\begin{alignat*}{3}
 &c_{ij}=\frac{2(\alpha_{i},\alpha_{j})}{(\alpha_{i},\alpha_{i})}, \qquad&&
 c_{0j} = \frac{2(-\theta,\alpha_{j})}{(\theta,\theta)}, &\\
 &c_{i0} = \frac{2(\alpha_{i},-\theta)}{(\alpha_{i},\alpha_{i})}, \qquad&&
 c_{00}=2 \qquad \mathrm{for\ all} \ 1\leq i,j\leq n.&
\end{alignat*}
Thus, $\hat{C}$ is the extended Cartan matrix of the Cartan matrix $C$ for $\fg$ from Section~\ref{ssec:general-def}.

\begin{Remark}
For the reader's convenience, we specify the new values of $\hat{C}$ and $\Omega$ in classical types:
\begin{itemize}\itemsep=0pt
\item Type $A_{n}$ $(n\geq 2)$:
\begin{gather*}
\begin{aligned}
 &\Omega_{00}=rs^{-1},\quad&& \Omega_{01}=r^{-1},\qquad&& \Omega_{0n}=s, \\
 &\Omega_{10}=s,\quad&& \Omega_{n0}=r^{-1}, \qquad&& \Omega_{0i}=\Omega_{i0}=1\qquad \text{for}\ 1<i<n;
 \end{aligned}
\\
c_{01}=c_{0n}=c_{10}=c_{n0}=-1,\qquad c_{0i}=c_{i0}=0 \qquad \text{for}\ 1<i<n.
\end{gather*}

\item Type $B_n$ $(n\geq 3)$:
\begin{gather*}
\begin{aligned}
 &\Omega_{00}=r^2s^{-2} ,\qquad&& \Omega_{01}=r^{-2}s^{-2} , \quad&& \Omega_{02}=r^{-2} , \\
 &\Omega_{0n}=r^2s^2 ,\qquad &&\Omega_{10}=r^2s^2 , \quad&& \Omega_{20}=s^2 ,\\
 & \Omega_{n0}=r^{-2}s^{-2} ,\qquad&& \Omega_{0i}=\Omega_{i0}=1 \qquad&& \text{for}\ 2<i<n;
 \end{aligned}
\\
 c_{02}=c_{20}=-1 , \qquad c_{01}=c_{10}=0 , \qquad c_{0i}=c_{i0}=0 \qquad\text{for}\ 2<i\leq n .
\end{gather*}

\item Type $C_n$ $(n\geq 2)$:
\begin{gather*}
\begin{aligned}
 &\Omega_{00}=r^2s^{-2}, \quad&& \Omega_{01}=r^{-2},\quad&& \Omega_{0n}=r^2s^2,\\
 & \Omega_{10}=s^2, \quad&& \Omega_{n0}=r^{-2}s^{-2} , \quad&& \Omega_{0i}=\Omega_{i0}=1 \qquad \text{for}\ 1<i<n ;
 \end{aligned}
\\
 c_{01}=-1,\qquad c_{10}=-2 , \qquad c_{0i}=c_{i0}=0 \qquad \text{for}\ 1<i\leq n .
\end{gather*}

\item Type $D_n$ $(n\geq 4)$:
\begin{gather*}
\begin{aligned}
 &\Omega_{00}=rs^{-1},\qquad&& \Omega_{01}=r^{-1}s^{-1},\qquad&& \Omega_{02}=r^{-1},\\
 &\Omega_{0n}=r^2s^2, \qquad&& \Omega_{10}=rs,\qquad&& \Omega_{20}=s,\\
 &\Omega_{n0}=r^{-2}s^{-2} , \qquad&& \Omega_{0i}=\Omega_{i0}=1 \qquad&& \text{for}\ 2<i<n ;
 \end{aligned}
\\
 c_{02}=c_{20}=-1 ,\qquad c_{01}=c_{10}=0 , \qquad c_{0i}=c_{i0}=0 \qquad \text{for}\ 2<i\leq n .
\end{gather*}
\end{itemize}
\end{Remark}

\begin{Definition}\label{def:2_param_aff}
The \textit{two-parameter quantum affine algebra} $\UU$ is the associative $\mathbb{K}$-algebra generated by
 \smash{$\bigl\{e_{i},f_{i},\omega_{i}^{\pm 1},\bigl(\omega_{i}'\bigr)^{\pm 1}\bigr\}_{i=0}^{n}\cup
 \bigl\{\gamma^{\pm 1},(\gamma')^{\pm 1}\bigr\} \cup \bigl\{D^{\pm 1}, (D')^{\pm 1}\bigr\}$}
with the following defining relations:
\begin{gather}
 D^{\pm 1}\cdot D^{\mp 1}=1, \qquad (D')^{\pm 1}\cdot (D')^{\mp 1}=1, \qquad DD'=D'D, \nonumber\\
 [D,\omega_i]=0, \qquad \bigl[D,\omega'_i\bigr]=0, \qquad \bigl[D',\omega_i\bigr]=0, \qquad \bigl[D',\omega'_i\bigr]=0, \nonumber\\
 D e_i D^{-1}=r_i^{\delta_{0i}}e_i, \quad D f_i D^{-1}=r_i^{-\delta_{0i}}f_i, \quad
 D' e_i \bigl(D'\bigr)^{-1}=s_i^{\delta_{0i}}e_i, \quad D' f_i \bigl(D'\bigr)^{-1}=s_i^{-\delta_{0i}}f_i,\nonumber
\\ \label{eq:aR0}
 \gamma=\omega_{\delta}=\omega_0\omega_{\theta},\qquad \gamma'=\omega'_{\delta}=\omega'_0\omega'_{\theta}
\quad \text{are central\ elements},
\\ \label{eq:aR1}
 [\omega_i,\omega_j]=\big[\omega_i,\omega'_j\big]=\big[\omega'_i,\omega'_j\big]=0, \qquad
 \omega_{i}^{\pm 1}\omega_{i}^{\mp 1} = 1 = \bigl(\omega_{i}'\bigr)^{\pm 1}\bigl(\omega_{i}'\bigr)^{\mp 1},
\\ \label{eq:aR2}
 \omega_{i}e_{j} = \Omega_{ji} e_{j}\omega_{i}, \qquad
 \omega_{i}f_{j} = \Omega^{-1}_{ji} f_{j}\omega_{i},
\\ \label{eq:aR3}
 \omega_{i}'e_{j} = \Omega^{-1}_{ij} e_{j}\omega_{i}', \qquad
 \omega_{i}'f_{j} = \Omega_{ij} f_{j}\omega_{i}',
\\ \label{eq:aR4}
 e_{i}f_{j} - f_{j}e_{i} = \delta_{ij}\frac{\omega_i-\omega'_i}{r_{i} - s_{i}},
\\
\begin{aligned}
& \sum_{k = 0}^{1 - c_{ij}} (-1)^k \qbinom{1 - c_{ij}}{k}_{r_{i},s_{i}}
 (r_{i}s_{i})^{\frac{1}{2}k(k - 1)} \Omega_{ji}^{k} s_{i}^{kc_{ij}} e_{i}^{1 - c_{ij} - k}e_{j}e_{i}^{k} = 0
 \qquad \forall i\ne j, \\
& \sum_{k = 0}^{1 - c_{ij}} (-1)^k \qbinom{1 - c_{ij}}{k}_{r_{i},s_{i}}
 (r_{i}s_{i})^{\frac{1}{2}k(k - 1)} \Omega_{ji}^{k} s_{i}^{kc_{ij}} f_{i}^{k}f_{j}f_{i}^{1 - c_{ij} - k} = 0
 \qquad \forall i\ne j.
\end{aligned} \label{eq:aR5}
\end{gather}
\end{Definition}

 \begin{Remark}
We note that the relations~\eqref{eq:R5} are compatible with~\eqref{eq:aR5} as $\Omega_{ji}s_i^{c_{ij}}=(rs)^{\langle \alpha_j,\alpha_i\rangle}$.
\end{Remark}

It is often more convenient to work with the version of $\UU$ without the degree generators~$D$,~$D'$.

\begin{Definition}
Let $\Uu$ be the associative $\mathbb{K}$-algebra generated by
$\smash{\bigl\{e_{i},f_{i},\omega_{i}^{\pm 1},\bigl(\omega_{i}'\bigr)^{\pm 1}\bigr\}_{i=0}^{n}}\cup \smash{\bigl\{\gamma^{\pm 1},(\gamma')^{\pm 1}\bigr\}}$
with the defining relations~\eqref{eq:aR0}--\eqref{eq:aR5}.
\end{Definition}

We now extend the $\uu$-modules $V$ from Section~\ref{sec:column_repr} to $\Uu$-modules $V(u)$ with $u\in\BC^\times$, depending on two additional parameters $a,b\in \BC^\times$. We shall further extend them to $\UU$-modules $V(u)$ with $u$ viewed as an indeterminate. We call the resulting $\Uu$- and $\UU$-modules as \emph{evaluation representations} \smash{\smash{$\rho^{a,b}_u$}}.

\begin{Proposition}[$A_n$, $n\geq 1$]\label{prop:A_aff_rep}
For any nonzero $a,b\in \BC$ set $c=rs\cdot ab$. Then the $U_{r,s}(\ssl_{n+1})$-action $\rho$ on $V$ from Proposition~$\ref{prp:A_rep}$ can be extended to a $U'_{r,s}(\widehat{\ssl}_{n+1})$-action \smash{$\rho^{a,b}_u$} on the vector space $V(u)=V$ by setting
\begin{displaymath}
 \rho^{a,b}_u(x)=\rho(x) \qquad \mathrm{for\ all} \ x\in \bigl\{e_i,f_i,\omega_i,\omega'_i\bigr\}_{i=1}^{n}
\end{displaymath}
and defining the action of $e_0$, $f_0$, $\omega_0$, $\omega'_0$, $\gamma$, $\gamma'$ via
\begin{displaymath}
\begin{split}
 & \rho^{a,b}_u(e_0)=au\cdot E_{n+1,1}, \qquad \rho^{a,b}_u(f_0)=bu^{-1}\cdot E_{1,n+1},
 \qquad \rho^{a,b}_u(\gamma)=c \, \Id=\rho^{a,b}_u\bigl(\gamma'\bigr), \\
 & \rho^{a,b}_u(\omega_0)=c\Biggl( r^{-1}E_{11}+r^{-1}s^{-1}\sum_{i=2}^{n}E_{ii}+s^{-1}E_{n+1,n+1} \Biggr), \\
 & \rho^{a,b}_u\bigl(\omega'_0\bigr)=c\Biggl( s^{-1}E_{11}+r^{-1}s^{-1}\sum_{i=2}^{n}E_{ii}+r^{-1}E_{n+1,n+1} \Biggr).
\end{split}
\end{displaymath}
\end{Proposition}

\begin{Proposition}[$B_n$, $n\geq 2$]\label{prop:B_aff_rep}
For any nonzero $a,b\in \BC$ set $c=(rs)^2 ab$. Then the $U_{r,s}(\sso_{2n+1})$-action $\rho$ on $V$ from Proposition~$\ref{prp:B_rep}$ can be extended to a $U'_{r,s}(\widehat{\sso}_{2n+1})$-action \smash{$\rho^{a,b}_u$} on the vector space $V(u)=V$ by~setting
\begin{displaymath}
 \rho^{a,b}_u(x)=\rho(x) \qquad \mathrm{for\ all} \ x\in \bigl\{e_i,f_i,\omega_i,\omega'_i\bigr\}_{i=1}^{n}
\end{displaymath}
and defining the action of $e_0$, $f_0$, $\omega_0$, $\omega'_0$, $\gamma$, $\gamma'$ via
\begin{align*}
 & \rho^{a,b}_u(e_0)=au \cdot \bigl(E_{1' 2} - r^2s^2 E_{2' 1} \bigr), \qquad
 \rho^{a,b}_u(f_0)=bu^{-1} \cdot \bigl(E_{2 1'}-E_{1 2'}), \\
 &\rho^{a,b}_u(\gamma)=c \,\Id=\rho^{a,b}_u(\gamma') , \\
 & \rho^{a,b}_u(\omega_0)=
 c\Biggl( s^2E_{11}+r^{-2}E_{22}+\sum_{i=3}^{n} \bigl( r^{-2}s^{-2}E_{ii} + r^2s^2 E_{i'i'}\bigr) +
 r^2E_{2' 2'} \\
 & \hphantom{\rho^{a,b}_u(\omega_0)=c\Biggl(}{}
 + s^{-2}E_{1' 1'} + E_{n+1,n+1} \Biggr), \\
 & \rho^{a,b}_u\big(\omega'_0\big)=
 c\Biggl( r^2E_{11}+s^{-2}E_{22}+\sum_{i=3}^{n} \bigl( r^{-2}s^{-2}E_{ii} + r^2s^2 E_{i'i'}\bigr) +s^2E_{2' 2'}\\
 & \hphantom{\rho^{a,b}_u(\omega'_0)=c\Biggl(}{}
 + r^{-2}E_{1' 1'}+E_{n+1,n+1} \Biggr).
\end{align*}
\end{Proposition}

\begin{Proposition}[$C_n$, $n\geq 2$]\label{prop:C_aff_rep}
For any nonzero $a,b\in \BC$ set $c=rs\cdot ab$. Then the $U_{r,s}(\ssp_{2n})$-action $\rho$ on $V$ from Proposition~$\ref{prp:C_rep}$ can be extended to a $U'_{r,s}(\widehat{\ssp}_{2n})$-action \smash{$\rho^{a,b}_u$} on the vector space $V(u)=V$ by setting
\begin{displaymath}
 \rho^{a,b}_u(x)=\rho(x) \qquad \mathrm{for\ all} \ x\in \bigl\{e_i,f_i,\omega_i,\omega'_i\bigr\}_{i=1}^{n}
\end{displaymath}
and defining the action of $e_0$, $f_0$, $\omega_0$, $\omega'_0$, $\gamma$, $\gamma'$ via
\begin{displaymath}
\begin{split}
 & \rho^{a,b}_u(e_0)=au \cdot E_{1' 1}, \qquad \rho^{a,b}_u(f_0)=bu^{-1} \cdot E_{1 1'}, \qquad
 \rho^{a,b}_u(\gamma)=c\, \Id=\rho^{a,b}_u(\gamma'),\\
 & \rho^{a,b}_u(\omega_0)=
 c\Biggl( r^{-1}sE_{11}+\sum_{i=2}^{n} \bigl( r^{-1}s^{-1}E_{ii} + rs E_{i'i'}\bigr) + rs^{-1}E_{1' 1'} \Biggr), \\
 & \rho^{a,b}_u\big(\omega'_0\big)=
 c\Biggl( rs^{-1}E_{11}+\sum_{i=2}^{n} \bigl( r^{-1}s^{-1}E_{ii} + rs E_{i'i'}\bigr) + r^{-1}sE_{1' 1'} \Biggr).
\end{split}
\end{displaymath}
\end{Proposition}

\begin{Proposition}[$D_n$, $n\geq 3$]\label{prop:D_aff_rep}
For any nonzero $a,b\in \BC$ set $c=rs\cdot ab$. Then the $U_{r,s}(\sso_{2n})$-action $\rho$ on $V$ from Proposition~$\ref{prp:D_rep}$ can be extended to a $U'_{r,s}\bigl(\widehat{\sso}_{2n}\bigr)$-action \smash{$\rho^{a,b}_u$} on the vector space $V(u)=V$ by setting
\begin{displaymath}
 \rho^{a,b}_u(x)=\rho(x) \qquad \mathrm{for\ all} \ x\in \bigl\{e_i,f_i,\omega_i,\omega'_i\bigr\}_{i=1}^{n}
\end{displaymath}
and defining the action of $e_0$, $f_0$, $\omega_0$, $\omega'_0$, $\gamma$, $\gamma'$ via
\begin{displaymath}
\begin{split}
 & \rho^{a,b}_u(e_0)=au \cdot \bigl(E_{1' 2} - rs E_{2' 1} \bigr), \qquad
 \rho^{a,b}_u(f_0)=bu^{-1} \cdot \bigl(E_{2 1'}-E_{1 2'}), \\
 &\rho^{a,b}_u(\gamma)=c \, \Id=\rho^{a,b}_u(\gamma') , \\
 & \rho^{a,b}_u(\omega_0)=
 c\Biggl( sE_{11}+r^{-1}E_{22}+\sum_{i=3}^{n} \bigl( r^{-1}s^{-1}E_{ii} + rs E_{i'i'}\bigr) +
 rE_{2' 2'} + s^{-1}E_{1' 1'} \Biggr), \\
 & \rho^{a,b}_u\big(\omega'_0\big)=
 c\Biggl( rE_{11}+s^{-1}E_{22}+\sum_{i=3}^{n} \bigl( r^{-1}s^{-1}E_{ii} + rs E_{i'i'}\bigr) +
 sE_{2' 2'} + r^{-1}E_{1' 1'} \Biggr).
\end{split}
\end{displaymath}
\end{Proposition}

These evaluation $\Uu$-modules \smash{$\rho^{a,b}_u$} can be naturally upgraded to $\UU$-modules in a~standard way.

\begin{Proposition}\label{prop:affine-actions-withD}
Let $u$ be an indeterminate and redefine the vector space via $V(u)=V\otimes_\BC \BC\big[u,u^{-1}\big]$ accordingly. Then, the formulas defining \smash{$\rho^{a,b}_u$} on the generators from Propositions~{\rm \ref{prop:A_aff_rep}--\ref{prop:D_aff_rep}} together with
\begin{gather*}
 \rho^{a,b}_u(D) \bigl(v\otimes u^k\bigr) = r_0^k \cdot v\otimes u^k , \qquad
 \rho^{a,b}_u\bigl(D'\bigr) \bigl(v\otimes u^k\bigr) = s_0^k \cdot v\otimes u^k \qquad
 \mathrm{for\ all} \ v\in V,\ k\in \BZ
\end{gather*}
give rise to the same-named action \smash{$\rho^{a,b}_u$} of $\UU$ on $V(u)$.
\end{Proposition}

The proofs of all these results are straightforward, cf.\ our proof of Proposition~\ref{prp:B_rep}.

\subsection[Affine R-matrices]{Affine $\boldsymbol{R}$-matrices}

Let $U^{\geq}_{r,s}(\widehat{\fg})$ be the subalgebra of $\UU$ generated by \smash{$\bigl\{e_i,\omega_i^{\pm 1},\gamma^{\pm 1},D^{\pm 1}\bigr\}_{i=0}^n$} and similarly let $U^{\leq}_{r,s}(\widehat{\fg})$ be the subalgebra generated by \smash{$\bigl\{f_i,(\omega'_i)^{\pm 1},(\gamma')^{\pm 1},(D')^{\pm 1}\bigr\}_{i=0}^n$}. Likewise, we define the~subalgebras \smash{$U^{',\geq}_{r,s}(\widehat{\fg})$} and \smash{$U^{',\leq}_{r,s}(\widehat{\fg})$} of $\Uu$. We note that the same formulas as in Section~\ref{ssec:general-def} for $\uu$ can be used to define the Hopf algebra structures on both $\UU$ and $\Uu$, so that \smash{$U^{\geq}_{r,s}(\widehat{\fg})$}, \smash{$U^{\leq}_{r,s}(\widehat{\fg})$} and \smash{$U^{',\geq}_{r,s}(\widehat{\fg})$}, \smash{$U^{',\leq}_{r,s}(\widehat{\fg})$} are also Hopf subalgebras of $\UU$ and $\Uu$. Finally, similarly to Proposition~\ref{prop:pairing_2param}, one has bilinear Hopf pairings
\begin{align}
 & (\cdot,\cdot)\colon\ U^{\leq}_{r,s}(\widehat{\fg})\times U^{\geq}_{r,s}(\widehat{\fg}) \to \mathbb{K} ,\nonumber\\
 & (\cdot,\cdot)\colon\ U^{',\leq}_{r,s}(\widehat{\fg}) \times U^{',\geq}_{r,s}(\widehat{\fg}) \to \mathbb{K}.\label{eq:parity-affine}
\end{align}
We note that the second of these pairings is actually degenerate as $\bigl(\smash{\gamma'-1,U^{',\geq}_{r,s}(\widehat{\fg})}\bigr)=0=\smash{\bigl(U^{',\leq}_{r,s}(\widehat{\fg}),\gamma-1\bigr)}$. On the other hand \big(which is one of the key reasons to add the generators~$D$,~$D'$\big), the first pairing in~\eqref{eq:parity-affine} is non-degenerate, and hence allows to realize the two-parameter quantum affine algebra $\UU$ as a Drinfel'd double of its Hopf subalgebras $U^{\leq}_{r,s}(\widehat{\fg})$, $U^{\geq}_{r,s}(\widehat{\fg})$ with respect to the pairing above.

The above discussion yields the universal $R$-matrix for $\UU$, which induces intertwiners \smash{$V\otimes W\iso W\otimes V$} for suitable $\UU$-modules $V,W$, akin to Section~\ref{ssec:universal-R}. In order to not overburden the paper, we choose to skip the detailed presentation on this standard but rather technical discussion. Instead, we shall now proceed directly to the main goal of this paper--the evaluation of such intertwiners when $V=\rho^{a,b}_u$, $W=\rho^{a,b}_v$ are the modules from Section~\ref{ssec:2-param affine}. In this context, we are looking for $\UU$-module intertwiners $\hat{R}(u/v)$ satisfying
\begin{equation}\label{eq:affine-intertwiner}
 \hat{R}(u/v)\circ \bigl(\rho^{a,b}_u\otimes \rho^{a,b}_v\bigr)(x) =
 \bigl(\rho^{a,b}_v\otimes \rho^{a,b}_u\bigr)(x)\circ \hat{R}(u/v)
\end{equation}
for all $x\in \UU$ (equivalently, one can rather request $x\in \Uu$ in the context of $\Uu$-modules). According to~\cite[Proposition~2]{Jim}, it suffices to check the validity of~\eqref{eq:affine-intertwiner} only for $x=f_i$, $0\leq i\leq n$. In fact, the space of such solutions is one-dimensional, see~\cite[Proposition~1]{Jim}, which is essentially due to the irreducibility of the tensor product $\rho^{a,b}_u\otimes \rho^{a,b}_v$ (which still holds when viewing them as $\Uu$-modules as long as $u$, $v$ are \emph{generic}), in contrast to Proposition~\ref{prop:struct}. As an immediate corollary, see~\cite[Proposition~3]{Jim}, the operator \smash{$R(u/v)=\hat{R}(u/v)\circ \tau$} satisfies the aforementioned Yang--Baxter relation with a spectral parameter:
\begin{align}
 & R_{12}(v/w) R_{13}(u/w) R_{23}(u/v) = R_{23}(u/v) R_{13}(u/w) R_{12}(v/w),\nonumber \\
 & \hat{R}_{12}(v/w) \hat{R}_{23}(u/w) \hat{R}_{12}(u/v) =
 \hat{R}_{23}(u/v) \hat{R}_{12}(u/w) \hat{R}_{23}(v/w),\label{eq:qYB-affine}
\end{align}
with notation as in Section~\ref{ssec:universal-R}. We shall now present explicit formulas for such $\hat{R}(z)$ in all classical types, generalizing~\cite{Jim} for the one-parameter setup.\footnote{We note that a twist of~\cite{R} is needed to recover the formulas of~\cite[Section~3]{Jim} due to a different coproduct~\cite[equation~(2.10)]{Jim} on $U_q(\widehat{\fg})$, cf.~Remark~\ref{remark6.14}.} We note that the origin of these formulas will be explained in the next section, where they will be derived through the so-called \emph{Yang--Baxterization} technique of~\cite{GWX}. However, once the formulas are guessed, the above discussions imply that it suffices to check that they satisfy~\eqref{eq:affine-intertwiner} for $x=f_i$, $0\leq i\leq n$.

We start with the simplest case of $A$-type (part~(b) of which goes back to~\cite[Section~2]{JL2}).

\begin{Theorem}[type $A_{n}$]\label{thm:A_AffRMatrix}\quad
\begin{enumerate}
\item[$(a)$] Let $z=u/v$. For \smash{$U_{r,s}\bigl(\widehat{\ssl}_{n+1}\bigr)$}-modules \smash{$\rho_u^{a,b}$} and \smash{$\rho_v^{a,b}$} from Proposition~$\ref{prop:A_aff_rep}$, the following operator $\hat{R}(z)$ satisfies~\eqref{eq:affine-intertwiner} whenever $ab = (rs)^{-1}$:
\begin{align}
 \hat{R}(z) &=
 \bigl(1 -zrs^{-1}\bigr)\sum_{i = 1}^{n + 1} E_{ii} \otimes E_{ii} + (1 - z)r\sum_{i > j} E_{ij} \otimes E_{ji} +
 (1 - z)s^{-1}\sum_{i < j}E_{ij} \otimes E_{ji} \nonumber\\
 & \quad + \bigl(1 - rs^{-1}\bigr)\sum_{i > j}E_{ii} \otimes E_{jj} +
 \bigl(1 - rs^{-1}\bigr)z\sum_{i < j}E_{ii} \otimes E_{jj}.\label{eq:A_AffRMatrix}
\end{align}

\item[$(b)$] The operator $R(z)=\hat{R}(z)\circ \tau$ given explicitly by
\begin{displaymath}
\begin{split}
 R(z) &=
 \bigl(1 -zrs^{-1}\bigr)\sum_{i = 1}^{n + 1} E_{ii} \otimes E_{ii} + (1 - z)r\sum_{i > j} E_{ii} \otimes E_{jj} +
 (1 - z)s^{-1}\sum_{i < j} E_{ii} \otimes E_{jj} \\
 & \quad + \bigl(1 - rs^{-1}\bigr)\sum_{i > j} E_{ij} \otimes E_{ji} +
 \bigl(1 - rs^{-1}\bigr)z\sum_{i < j} E_{ij} \otimes E_{ji}
\end{split}
\end{displaymath}
satisfies the Yang--Baxter equation with a spectral parameter~\eqref{eq:qYB-affine}.
\end{enumerate}
\end{Theorem}

The main results of this section generalize the above theorem to the other classical series.

\begin{Theorem}[type $B_{n}$]\label{thm:B_AffRMatrix}\quad
\begin{enumerate}
\item[$(a)$] Let $z=u/v$. For \smash{$U_{r,s}\bigl(\widehat{\sso}_{2n+1}\bigr)$}-modules \smash{$\rho_u^{a,b}$} and \smash{$\rho_v^{a,b}$} from Proposition~$\ref{prop:B_aff_rep}$, the following operator $\hat{R}(z)$ satisfies~\eqref{eq:affine-intertwiner} whenever $ab = (rs)^{-2}$:
\begin{gather}
 \hat{R}(z) = \bigl(z - r^{-2}s^{2}\bigr)(z - \xi) \sum_{1 \le i \le 2n + 1}^{i \neq n + 1} E_{ii} \otimes E_{ii} + \sum_{1 \le i,j \le 2n + 1}^{j \neq i,i'} a_{ij}(z) E_{ij} \otimes E_{ji} \nonumber\\
 \hphantom{\hat{R}(z) =}{}
 + \bigl(1 - r^{-2}s^{2}\bigr)(z - \xi) \sum_{i > j}^{j \neq i'} E_{ii} \otimes E_{jj} + \bigl(1 - r^{-2}s^{2}\bigr)z(z - \xi) \sum_{i < j}^{j \neq i'} E_{ii} \otimes E_{jj} \nonumber\\
 \hphantom{\hat{R}(z) =}{}
 + \sum_{i,j = 1}^{2n + 1} b_{ij}(z) E_{i'j} \otimes E_{ij'},\label{eq:B_AffRMatrix}
\end{gather}
where $\xi = r^{-2n +1}s^{2n - 1}$,
\begin{gather*}
 a_{ij}(z) =
 \begin{cases}
 r^{-1}s(z - 1)(z - \xi)(rs)^{-\sigma_{i}\sigma_{j}} &
 \text{if}\ i < j,j'\ \text{or}\ i > j,j', \\
 r^{-1}s(z - 1)(z - \xi)(rs)^{\sigma_{i}\sigma_{j}} &
 \text{if}\ j < i < j'\ \text{or}\ j' < i < j,
 \end{cases}
\\
 b_{ij}(z) =
 \begin{cases}
 \bigl(r^{-2}s^{2}z - \xi\bigr)(z - 1) & \text{if}\ j = i,\ i \neq n + 1 ,\\
 r^{-1}s(z - 1)(z - \xi) + \bigl(r^{-2}s^{2} - 1\bigr)(\xi - 1)z & \text{if}\ i = j = n+ 1 ,\\
 \bigl(r^{-2}s^{2} - 1\bigr)\bigl(\xi t_{i}t_{j}^{-1}(z - 1) - \delta_{ij'}(z - \xi)\bigr) & \text{if}\ i < j ,\\
 \bigl(r^{-2}s^{2} - 1\bigr)z\bigl(t_{i}t_{j}^{-1}(z - 1) - \delta_{ij'}(z - \xi)\bigr) & \text{if}\ i > j,
 \end{cases}
\end{gather*}
with $t_{i}$, $\sigma_{i}$ precisely as in~\eqref{eq:B_coeffs}.

\item[$(b)$] The operator $R(z)=\hat{R}(z)\circ \tau$ is a solution of the Yang--Baxter equation with a spectral parameter~\eqref{eq:qYB-affine}.
\end{enumerate}
\end{Theorem}

\begin{Theorem}[type $C_{n}$]\label{thm:C_AffRMatrix}\quad
\begin{enumerate}
\item[$(a)$] Let $z=u/v$. For \smash{$U_{r,s}\big(\widehat{\ssp}_{2n}\big)$}-modules \smash{$\rho_u^{a,b}$} and \smash{$\rho_v^{a,b}$} from Proposition~$\ref{prop:C_aff_rep}$, the following operator $\hat{R}(z)$ satisfies~\eqref{eq:affine-intertwiner} whenever $ab = (rs)^{-1}$:
\begin{gather}
 \hat{R}(z) = \bigl(z - r^{-1}s\bigr)(z - \xi)\sum_{i = 1}^{2n}E_{ii} \otimes E_{ii} + \sum_{1 \le i,j \le 2n}^{j \neq i,i'}a_{ij}(z)E_{ij} \otimes E_{ji}\nonumber \\
 \hphantom{\hat{R}(z) =}{} + \bigl(1 - r^{-1}s\bigr)(z - \xi)\sum_{i > j}^{j \neq i'}E_{ii} \otimes E_{jj} +\bigl(1 - r^{-1}s\bigr)z(z - \xi)\sum_{i < j}^{j \neq i'}E_{ii} \otimes E_{jj} \nonumber\\
 \hphantom{\hat{R}(z) =}{} + \sum_{i,j = 1}^{2n} b_{ij}(z)E_{i'j} \otimes E_{ij'},\label{eq:C_AffRMatrix}
\end{gather}
where $\xi = r^{-n - 1}s^{n + 1}$,
\begin{gather*}
 a_{ij}(z) =
 \begin{cases}
 r^{-1/2}s^{1/2}(z - 1)(z - \xi)(rs)^{-\frac{1}{2}\sigma_{i}\sigma_{j}} &
 \text{if}\ i < j,j'\ \text{or}\ i > j,j', \\
 r^{-1/2}s^{1/2}(z - 1)(z - \xi)(rs)^{\frac{1}{2}\sigma_{i}\sigma_{j}} &
 \text{if}\ j < i < j'\ \text{or}\ j' < i < j,
 \end{cases}
\\
 b_{ij}(z) =
 \begin{cases}
 \bigl(r^{-1}sz - \xi\bigr)(z - 1) & \text{if}\ j = i, \\
 \bigl(r^{-1}s - 1\bigr)\bigl(\xi t_{i}t_{j}^{-1}(z - 1) - \delta_{ij'}(z - \xi)\bigr) & \text{if}\ i < j ,\\
 \bigl(r^{-1}s - 1\bigr)z\bigl(t_{i}t_{j}^{-1}(z - 1) - \delta_{ij'}(z - \xi)\bigr) & \text{if}\ i > j,
 \end{cases}
\end{gather*}
with $t_{i}$, $\sigma_{i}$ precisely as in~\eqref{eq:C_coeffs}.

\item[$(b)$] The operator $R(z)=\hat{R}(z)\circ \tau$ is a solution of the Yang--Baxter equation with a spectral parameter~\eqref{eq:qYB-affine}.
\end{enumerate}
\end{Theorem}

\begin{Theorem}[type $D_{n}$]\label{thm:D_AffRMatrix} \quad
\begin{enumerate}
\item[$(a)$] Let $z = u/v$. For \smash{$U_{r,s}\bigl(\widehat{\sso}_{2n}\bigr)$}-modules \smash{$\rho_u^{a,b}$} and \smash{$\rho_v^{a,b}$} from Proposition~$\ref{prop:D_aff_rep}$, the following operator $\hat{R}(z)$ satisfies~\eqref{eq:affine-intertwiner} whenever $ab = (rs)^{-1}$:
\begin{gather}
 \hat{R}(z) =
 \bigl(z - r^{-1}s\bigr)(z - \xi)\sum_{i = 1}^{2n}E_{ii} \otimes E_{ii} + \sum_{1 \le i,j \le 2n}^{j \neq i,i'}a_{ij}(z)E_{ij} \otimes E_{ji} \nonumber\\
 \hphantom{ \hat{R}(z) =}{}
 + \bigl(1 - r^{-1}s\bigr)(z - \xi)\sum_{i > j}^{j \neq i'}E_{ii} \otimes E_{jj} + \bigl(1 - r^{-1}s\bigr)z(z - \xi)\sum_{i < j}^{j \neq i'}E_{ii} \otimes E_{jj} \nonumber\\
 \hphantom{ \hat{R}(z) =}{}
 + \sum_{i,j = 1}^{2n} b_{ij}(z)E_{i'j} \otimes E_{ij'},\label{eq:D_AffRMatrix}
\end{gather}
where $\xi = r^{-n + 1}s^{n - 1}$,
\begin{gather*}
 a_{ij}(z) =
 \begin{cases}
 r^{-1/2}s^{1/2}(z - 1)(z - \xi)(rs)^{-\frac{1}{2}\sigma_{i}\sigma_{j}} &
 \text{if}\ i < j,j'\ \text{or}\ i > j,j', \\
 r^{-1/2}s^{1/2}(z - 1)(z - \xi)(rs)^{\frac{1}{2}\sigma_{i}\sigma_{j}} &
 \text{if}\ j < i < j'\ \text{or}\ j' < i < j,
 \end{cases}
\\
 b_{ij}(z) =
 \begin{cases}
 \bigl(r^{-1}sz - \xi\bigr)(z - 1) & \text{if}\ j = i ,\\
 \bigl(r^{-1}s - 1\bigr)\bigl(\xi t_{i}t_{j}^{-1}(z - 1) - \delta_{ij'}(z - \xi)\bigr) & \text{if}\ i < j ,\\
 \bigl(r^{-1}s - 1\bigr)z\bigl(t_{i}t_{j}^{-1}(z - 1) - \delta_{ij'}(z - \xi)\bigr) & \text{if}\ i > j,
 \end{cases}
\end{gather*}
with $t_{i}$, $\sigma_{i}$ precisely as in~\eqref{eq:D_coeffs}.

\item[$(b)$] The operator $R(z)=\hat{R}(z)\circ \tau$ is a solution of the Yang--Baxter equation with a spectral parameter~\eqref{eq:qYB-affine}.
\end{enumerate}
\end{Theorem}

\begin{Remark}\label{remark6.14}
The careful reader may have noticed that for $r = q$ and $s = q^{-1}$, the formula~\eqref{eq:B_AffRMatrix} does not precisely match with the corresponding $R$-matrix of~\cite{Jim}. This discrepancy occurs because~\cite{Jim} uses a coproduct $\Delta'$ on $U_q(\fg)$ that differs from the more standard $\Delta$ of~\cite{J}. One may check that
\begin{displaymath}
 \Delta'(x) = \tilde{f}^{1/2}\Delta(x)\tilde{f}^{-1/2}\qquad \text{for all}\ x \in U_{q}(\fg),
\end{displaymath}
where both sides are regarded as operators on $V\otimes W$ for any finite-dimensional (type $1$) $U_{q}(\fg)$-modules $V$ and $W$ (see~\cite[Chapter~7]{J} for the notation) and $\tilde{f}$ denotes the $r = q$, $s = q^{-1}$ specialization of the corresponding map from Section~\ref{ssec:universal-R}. Using this relationship, one may then show that
\begin{displaymath}
 \hat{R}' = \tilde{f}^{1/2} \hat{R} \tilde{f}^{-1/2}
\end{displaymath}
is the corresponding intertwiner relative to the coproduct $\Delta'$ of~\cite{Jim}. The Yang--Baxterization results of Section~\ref{sec:Baxterization}, then imply that the one-parameter specialization of our~\eqref{eq:B_AffRMatrix}--\eqref{eq:D_AffRMatrix} should be related to the corresponding formulas in~\cite{Jim} via conjugation by $\tilde{f}^{-1/2}$. The latter can be verified directly case-by-case.
\end{Remark}

\subsection{Proofs of explicit formulas}

We shall now present the proofs of the theorems from the previous subsection. We start with $A$-type.

\begin{proof}[Proof of Theorem~\ref{thm:A_AffRMatrix}]
We shall only prove (a), since (b) follows from it (see also our discussion in Section~\ref{sec:Baxterization}).

As follows from Section~\ref{sec:Baxterization}, the operator $\hat{R}\bigl(uv^{-1}\bigr)$ of~\eqref{eq:A_AffRMatrix} is a linear combination of $\hat{R}$ from~\eqref{eq:A_RMatrix} and its inverse \smash{$\hat{R}^{-1}$}, hence it satisfies \eqref{eq:affine-intertwiner} for $x = e_{i},f_{i}$, $1 \le i \le n$. For the remaining generators, it will be helpful to record the explicit action of $\hat{R}\bigl(uv^{-1}\bigr)$ on the basis vectors $v_{i} \otimes v_{j}$ of $V(u) \otimes V(v)$:
\begin{equation}\label{eq:R-action-A}
 \hat{R}\bigl(uv^{-1}\bigr) v_{i} \otimes v_{j} =
 \begin{cases}
 \bigl(1 - uv^{-1}rs^{-1}\bigr)v_{i} \otimes v_{i} & \text{if}\ i = j, \\
 \bigl(1 - rs^{-1}\bigr)v_{i} \otimes v_{j} + s^{-1}(1 - uv^{-1})v_{j} \otimes v_{i} & \text{if}\ i > j, \\
 uv^{-1}\bigl(1 - rs^{-1}\bigr)v_{i} \otimes v_{j} + r(1 - uv^{-1})v_{j} \otimes v_{i} & \text{if}\ i < j.
 \end{cases}
\end{equation}
In particular, $\hat{R}\bigl(uv^{-1}\bigr)$ preserves the weight subspaces, so that \eqref{eq:affine-intertwiner} holds also for $x = \omega_{i},\omega_{i}'$, $0 \le i \le n$.

We shall now verify~\eqref{eq:affine-intertwiner} for $x = e_{0}$. First, we note that
\begin{gather}\label{eq:e0-action-A}
 a^{-1}(\rho_{u} \otimes \rho_{v})(e_{0})(v_{i} \otimes v_{j}) =
 \begin{cases}
 0 & \text{if}\ i \neq 1,\ j \neq 1, \\
 u v_{n+1} \otimes v_{j} & \text{if}\ i = 1,\ j \neq 1, \\
 u v_{n+1} \otimes v_{1} + vr^{-1} v_{1} \otimes v_{n+1} & \text{if}\ i = j = 1, \\
 vs^{-1} v_{n+1} \otimes v_{n+1} & \text{if}\ i = n + 1,\ j = 1 ,\\
 vr^{-1}s^{-1} v_{i} \otimes v_{n+1} & \text{if}\ i \neq 1,n+1,\ j = 1,
 \end{cases}
\end{gather}
and we have a similar formula for $a^{-1}(\rho_{v} \otimes \rho_{u})(e_{0})(v_{i} \otimes v_{j})$ with $u$ and $v$ interchanged. As~short\-hand, we shall use the notation $E_{uv}$ and $E_{vu}$ to denote the operators $a^{-1}(\rho_{u} \otimes \rho_{v})(e_{0})$ and $a^{-1}(\rho_{v} \otimes \rho_{u})(e_{0})$, respectively. According to~\eqref{eq:R-action-A} and~\eqref{eq:e0-action-A}, we clearly have
\[
\hat{R}\bigl(uv^{-1}\bigr)(E_{uv}(v_{i} \otimes v_{j})) = 0 =  E_{vu}\bigl(\hat{R}\bigl(uv^{-1}\bigr)(v_{i} \otimes v_{j})\bigr)
\qquad \text{for $i,j \neq 1$.}
\]
 We now consider the remaining cases:
\begin{itemize}\itemsep=0pt
\item
If $i = 1$ and $j \neq 1,n+1$, then
\begin{displaymath}
 \hat{R}\bigl(uv^{-1}\bigr)(E_{uv}(v_{1} \otimes v_{j})) =
 s^{-1}u\bigl(1 - uv^{-1}\bigr)v_{j} \otimes v_{n+ 1} + u\bigl(1 - rs^{-1}\bigr)v_{n+1} \otimes v_{j},
\end{displaymath}
which is equal to
\begin{displaymath}
\begin{split}
 E_{vu}\hat{R}\bigl(uv^{-1}\bigr)(v_{1} \otimes v_{j})
 & = E_{vu}\bigl(r\bigl(1 - uv^{-1}\bigr)v_{j} \otimes v_{1} + uv^{-1}\bigl(1 - rs^{-1}\bigr)v_{1} \otimes v_{j}\bigr) \\
 & = s^{-1}u\bigl(1 - uv^{-1}\bigr)v_{j} \otimes v_{n + 1} + u\bigl(1 - rs^{-1}\bigr)v_{n+1} \otimes v_{j}.
\end{split}
\end{displaymath}

\item
If $i = 1$ and $j = n+ 1$, then
\begin{displaymath}
 \hat{R}\bigl(uv^{-1}\bigr)(E_{uv}(v_{1} \otimes v_{n+1})) = u\bigl(1 - uv^{-1}rs^{-1}\bigr)v_{n+ 1} \otimes v_{n+ 1},
\end{displaymath}
which is equal to
\begin{displaymath}
\begin{split}
 E_{vu}\hat{R}\bigl(uv^{-1}\bigr)(v_{1} \otimes v_{n+1}) &=
 rs^{-1}u\bigl(1 - uv^{-1}\bigr)v_{n+1} \otimes v_{n + 1} + u\bigl(1 - rs^{-1}\bigr)v_{n +1} \otimes v_{n+1} \\
 &= u\bigl(1 - uv^{-1}rs^{-1}\bigr)v_{n+1} \otimes v_{n+1}.
\end{split}
\end{displaymath}

\item
If $i = j = 1$, then we have
\begin{displaymath}
\begin{split}
 \hat{R}\bigl(uv^{-1}\bigr)(E_{uv}(v_{1} \otimes v_{1})) &{}=
 \hat{R}\bigl(uv^{-1}\bigr)\bigl(uv_{n+1} \otimes v_{1} + vr^{-1}v_{1} \otimes v_{n+1}\bigr) \\
 & = us^{-1}\bigl(1 - uv^{-1}\bigr)v_{1} \otimes v_{n+1} + u\bigl(1 - rs^{-1}\bigr)v_{n+1} \otimes v_{1} \\
 &\quad{}+ ur^{-1}\bigl(1 - rs^{-1}\bigr)v_{1} \otimes v_{n+1} + v\bigl(1 - uv^{-1}\bigr) v_{n+1} \otimes v_{1} \\
 & = ur^{-1}\bigl(1 - uv^{-1}rs^{-1}\bigr)v_{1} \otimes v_{n+1} + \bigl(v - urs^{-1}\bigr)v_{n+1} \otimes v_{1},
\end{split}
\end{displaymath}
which is equal to
\begin{displaymath}
\begin{split}
 E_{vu}\hat{R}\bigl(uv^{-1}\bigr)(v_{1} \otimes v_{1}) &=
 E_{vu}\bigl(1 - uv^{-1}rs^{-1}\bigr)v_{1} \otimes v_{1} \\
 &= \bigl(v - urs^{-1}\bigr)v_{n+1} \otimes v_{1} + ur^{-1}\bigl(1 - uv^{-1}rs^{-1}\bigr)v_{1} \otimes v_{n+1}.
\end{split}
\end{displaymath}

\item
If $i = n + 1$ and $j = 1$, then
\begin{displaymath}
 \hat{R}\bigl(uv^{-1}\bigr)(E_{uv}(v_{n+ 1} \otimes v_{1})) = vs^{-1}\bigl(1 - uv^{-1}rs^{-1}\bigr)v_{n+1} \otimes v_{n+1},
\end{displaymath}
which is equal to
\begin{align*}
 &E_{vu}\hat{R}\bigl(uv^{-1}\bigr)(v_{n+1} \otimes v_{1})
= E_{vu}\bigl(s^{-1}\bigl(1 - uv^{-1}\bigr)v_{1} \otimes v_{n+1} + \bigl(1 - rs^{-1}\bigr)v_{n+1} \otimes v_{1}\bigr) \\
 &\qquad{}= vs^{-1}\bigl(1 - uv^{-1}\bigr)v_{n+1} \otimes v_{n+1} + us^{-1}\bigl(1 - rs^{-1}\bigr)v_{n + 1} \otimes v_{n+1} \\
 &\qquad{}= vs^{-1}\bigl(1 - uv^{-1}rs^{-1}\bigr)v_{n + 1} \otimes v_{n+1}.
\end{align*}

\item
Finally, if $i \neq 1,n + 1$ and $j = 1$, then
\begin{gather*}
 \hat{R}\bigl(uv^{-1}\bigr)(E_{uv}(v_{i} \otimes v_{1})) =
 vs^{-1}\bigl(1 - uv^{-1}\bigr)v_{n+1} \otimes v_{i} + ur^{-1}s^{-1}\bigl(1 - rs^{-1}\bigr)v_{i} \otimes v_{n+1},
\end{gather*}
which is equal to
\begin{displaymath}
\begin{split}
 E_{vu}\hat{R}\bigl(uv^{-1}\bigr)(v_{i} \otimes v_{1})
 & = E_{vu}\bigl(s^{-1}\bigl(1 - uv^{-1}\bigr)v_{1} \otimes v_{i} + \bigl(1 - rs^{-1}\bigr)v_{i} \otimes v_{1}\bigr) \\
 & = vs^{-1}\bigl(1 - uv^{-1}\bigr)v_{n+1} \otimes v_{i} + ur^{-1}s^{-1}\bigl(1 - rs^{-1}\bigr)v_{i} \otimes v_{n+1}.
\end{split}
\end{displaymath}

\end{itemize}
This completes our verification of~\eqref{eq:affine-intertwiner} for $x = e_{0}$.

The verification for $x = f_{0}$ is completely analogous. This completes the proof of Theorem~\ref{thm:A_AffRMatrix}\,(a).
\end{proof}

Let us now present the proof of Theorem \ref{thm:B_AffRMatrix} (the proofs of Theorems~\ref{thm:C_AffRMatrix} and~\ref{thm:D_AffRMatrix} are completely analogous).

\begin{proof}[Proof of Theorem~\ref{thm:B_AffRMatrix}]
As follows from Section~\ref{sec:Baxterization}, the operator \smash{$\hat{R}\bigl(uv^{-1}\bigr)$} of~\eqref{eq:B_AffRMatrix} is a linear combination of $\hat{R}$ from~\eqref{eq:B_RMatrix}, its inverse \smash{$\hat{R}^{-1}$}, and the identity operator $\Id$. Hence, it satisfies~\eqref{eq:affine-intertwiner} for $x=e_i,f_i$, $1\leq i\leq n$. Moreover, \smash{$\hat{R}\bigl(uv^{-1}\bigr)$} clearly preserves the weight subspaces, so it also satisfies \eqref{eq:affine-intertwiner} for $x = \omega_{i},\omega_{i}'$, $0 \le i \le n$.

We shall now verify~\eqref{eq:affine-intertwiner} for $x=f_0$, proceeding similarly to our proof of Lemma~\ref{lem:B_Intertwining}. To~this end, let \smash{$F_{uv} = b^{-1}\bigl(\rho_{u}^{a,b} \otimes \rho_{v}^{a,b}\bigr)(f_{0})$} and \smash{$F_{vu} = b^{-1}\bigl(\rho_{v}^{a,b} \otimes \rho_{u}^{a,b}\bigr)(f_{0})$}, where we assume that $ab = (rs)^{-2}$, so that $c = (rs)^{2}ab = 1$. Then we need to verify that
\[
\hat{R}\bigl(uv^{-1}\bigr)F_{uv}= F_{vu}\hat{R}\bigl(uv^{-1}\bigr)\in \End(V \otimes V).
\]
First, let us record the explicit formula for $F_{uv}\in \End(V \otimes V)$:
\begin{gather*}
 F_{uv} = v^{-1} \cdot 1 \otimes E_{21'} - v^{-1}\cdot 1 \otimes E_{12'} + u^{-1}r^{2}E_{21'} \otimes E_{11} + u^{-1}s^{-2}E_{21'} \otimes E_{22} \\
 \hphantom{ F_{uv} =}{}
 + u^{-1}s^{2} E_{21'} \otimes E_{2'2'}+ u^{-1}r^{-2}E_{21'} \otimes E_{1'1'} + u^{-1}E_{21'} \otimes E_{n + 1,n + 1} \\
 \hphantom{ F_{uv} =}{}
+ \sum_{i = 3}^{n}\bigl(u^{-1}(rs)^{-2}E_{21'} \otimes E_{ii} + u^{-1}(rs)^{2}E_{21'} \otimes E_{i'i'}\bigr) - u^{-1}r^{2}E_{12'} \otimes E_{11} \\
 \hphantom{ F_{uv} =}{}
- u^{-1}s^{-2}E_{12'} \otimes E_{22} - u^{-1}s^{2}E_{12'} \otimes E_{2'2'} - u^{-1}r^{-2}E_{12'} \otimes E_{1'1'} \\
 \hphantom{ F_{uv} =}{}
 -u^{-1}E_{12'} \otimes E_{n+1,n+1}- \sum_{i = 3}^{n}\bigl(u^{-1}(rs)^{-2}E_{12'} \otimes E_{ii} + u^{-1}(rs)^{2}E_{12'} \otimes E_{i'i'}\bigr),
\end{gather*}
while $F_{vu}$ is given by the same formula with $u$ and $v$ interchanged. As in the proof of Lemma~\ref{lem:B_Intertwining}, it will be helpful to break the operator $\hat{R}\bigl(uv^{-1}\bigr)$ from~\eqref{eq:B_AffRMatrix} into the following six terms:
\begin{gather*}
 R_{1}\bigl(uv^{-1}\bigr) = \bigl(uv^{-1} - r^{-2}s^{2}\bigr)\bigl(uv^{-1} - \xi\bigr) \sum_{1 \le i \le 2n + 1}^{i \neq n + 1} E_{ii} \otimes E_{ii} \\
 \qquad{}
 + \bigl(r^{-1}s\bigl(uv^{-1} - 1\bigr)\bigl(uv^{-1} - \xi\bigr) + \bigl(r^{-2}s^{2} - 1\bigr)(\xi - 1)uv^{-1}\bigr)E_{n+1,n + 1} \otimes E_{n+ 1,n + 1}, \\
 R_{2}\bigl(uv^{-1}\bigr) = \bigl(r^{-2}s^{2}uv^{-1} - \xi\bigr)\bigl(uv^{-1} - 1\bigr)\sum_{1 \le i \le 2n + 1}^{i \neq n + 1}E_{i'i} \otimes E_{ii'}, \\
 R_{3}\bigl(uv^{-1}\bigr) = \sum_{1 \le i,j \le 2n + 1}^{j \neq i,i'}a_{ij}\bigl(uv^{-1}\bigr)E_{ij} \otimes E_{ji}, \\
 R_{4}\bigl(uv^{-1}\bigr) = \bigl(1 - r^{-2}s^{2}\bigr)\bigl(uv^{-1} - \xi\bigr)\sum_{i > j}^{j \neq i'}E_{ii} \otimes E_{jj}, \\
 R_{5}\bigl(uv^{-1}\bigr) = \bigl(1 - r^{-2}s^{2}\bigr)uv^{-1}\bigl(uv^{-1} - \xi\bigr)\sum_{i < j}^{j \neq i'}E_{ii} \otimes E_{jj}, \\
 R_{6}\bigl(uv^{-1}\bigr) = \sum_{1 \le i,j \le 2n + 1}^{i \neq j}b_{ij}\bigl(uv^{-1}\bigr)E_{i'j} \otimes E_{ij'}.
\end{gather*}
Then a direct computation yields
\begin{gather*}
 R_{1}\bigl(uv^{-1}\bigr)F_{uv}= v^{-1}\bigl(uv^{-1} - r^{-2}s^{2}\bigr)\bigl(uv^{-1} - \xi\bigr)E_{22} \otimes E_{21'} \\
\hphantom{ R_{1}\bigl(uv^{-1}\bigr)F_{uv}=}{}
- v^{-1}\bigl(uv^{-1} - r^{-2}s^{2}\bigr)\bigl(uv^{-1} - \xi\bigr)E_{11} \otimes E_{12'} \\
\hphantom{ R_{1}\bigl(uv^{-1}\bigr)F_{uv}=}{}
 +u^{-1}s^{-2}\bigl(uv^{-1} - r^{-2}s^{2}\bigr)\bigl(uv^{-1} - \xi\bigr)E_{21'} \otimes E_{22} \\
\hphantom{ R_{1}\bigl(uv^{-1}\bigr)F_{uv}=}{}
 - u^{-1}r^{2}\bigl(uv^{-1} - r^{-2}s^{2}\bigr)\bigl(uv^{-1} - \xi\bigr)E_{12'} \otimes E_{11},
\\
 F_{vu}R_{1}\bigl(uv^{-1}\bigr) = u^{-1}\bigl(uv^{-1} - r^{-2}s^{2}\bigr)\bigl(uv^{-1} - \xi\bigr)E_{1'1'} \otimes E_{21'}\\
 \hphantom{F_{vu}R_{1}\bigl(uv^{-1}\bigr) =}{}
 - u^{-1}\bigl(uv^{-1} - r^{-2}s^{2}\bigr)\bigl(uv^{-1} - \xi\bigr)E_{2'2'} \otimes E_{12'} \\
 \hphantom{F_{vu}R_{1}\bigl(uv^{-1}\bigr) =}{}
 +v^{-1}r^{-2}\bigl(uv^{-1} - r^{-2}s^{2}\bigr)\bigl(uv^{-1} - \xi\bigr)E_{21'} \otimes E_{1'1'} \\
 \hphantom{F_{vu}R_{1}\bigl(uv^{-1}\bigr) =}{}
 - v^{-1}s^{2}\bigl(uv^{-1} - r^{-2}s^{2}\bigr)\bigl(uv^{-1} - \xi\bigr)E_{12'} \otimes E_{2'2'},
\\
 R_{2}\bigl(uv^{-1}\bigr)F_{uv} = v^{-1}\bigl(r^{-2}s^{2}uv^{-1} - \xi\bigr)\bigl(uv^{-1} - 1\bigr)E_{22'} \otimes E_{2'1'} \\
 \hphantom{R_{2}\bigl(uv^{-1}\bigr)F_{uv} = }{}
 - v^{-1}\bigl(r^{-2}s^{2}uv^{-1} - \xi\bigr)\bigl(uv^{-1} - 1\bigr)E_{11'} \otimes E_{1'2'} \\
 \hphantom{R_{2}\bigl(uv^{-1}\bigr)F_{uv} = }{}
 +u^{-1}s^{2}\bigl(r^{-2}s^{2}uv^{-1} - \xi\bigr)\bigl(uv^{-1} - 1\bigr)E_{2'1'} \otimes E_{22'} \\
 \hphantom{R_{2}\bigl(uv^{-1}\bigr)F_{uv} = }{}
 - u^{-1}r^{-2}\bigl(r^{-2}s^{2}uv^{-1} - \xi\bigr)\bigl(uv^{-1} - 1\bigr)E_{1'2'} \otimes E_{11'},
\\
 F_{vu}R_{2}\bigl(uv^{-1}\bigr) = u^{-1}\bigl(r^{-2}s^{2}uv^{-1} - \xi\bigr)\bigl(uv^{-1} - 1\bigr)E_{11'} \otimes E_{21}\\
 \hphantom{ F_{vu}R_{2}\bigl(uv^{-1}\bigr) =}{}
 - u^{-1}\bigl(r^{-2}s^{2}uv^{-1} - \xi\bigr)\bigl(uv^{-1} - 1\bigr)E_{22'} \otimes E_{12} \\
 \hphantom{ F_{vu}R_{2}\bigl(uv^{-1}\bigr) =}{}
 + v^{-1}r^{2}\bigl(r^{-2}s^{2}uv^{-1} - \xi\bigr)\bigl(uv^{-1} - 1\bigr)E_{21} \otimes E_{11'} \\
 \hphantom{ F_{vu}R_{2}\bigl(uv^{-1}\bigr) =}{}
 - v^{-1}s^{-2}\bigl(r^{-2}s^{2}uv^{-1} - \xi\bigr)\bigl(uv^{-1} - 1\bigr)E_{12} \otimes E_{22'},
\\
 R_{3}\bigl(uv^{-1}\bigr)F_{uv} - F_{vu}R_{3}\bigl(uv^{-1}\bigr) \\
\quad{} = s^{2}v^{-1}\bigl(uv^{-1} - 1\bigr)\bigl(uv^{-1} - \xi\bigr)E_{21} \otimes E_{11'}+ r^{-2}v^{-1}\bigl(uv^{-1} - 1\bigr)\bigl(uv^{-1} - \xi\bigr)E_{21'} \otimes E_{1'1'} \\
\qquad{} -r^{-2}v^{-1}\bigl(uv^{-1} - 1\bigr)\bigl(uv^{-1} - \xi\bigr)E_{12} \otimes E_{22'} - s^{2}v^{-1}\bigl(uv^{-1} - 1\bigr)\bigl(uv^{-1} - \xi\bigr)E_{12'} \otimes E_{2'2'} \\
\qquad{} + u^{-1}\bigl(uv^{-1} - 1\bigr)\bigl(uv^{-1} - \xi\bigr)E_{11'} \otimes E_{21}+r^{-2}s^{2}u^{-1}\bigl(uv^{-1} - 1\bigr)\bigl(uv^{-1} - \xi\bigr)E_{1'1'} \otimes E_{21'} \\
\qquad{}- r^{-2}u^{-1}\bigl(uv^{-1} - 1\bigr)\bigl(uv^{-1} - \xi\bigr)E_{21'} \otimes E_{22} - s^{2}u^{-1}\bigl(uv^{-1} - 1\bigr)\bigl(uv^{-1} - \xi\bigr)E_{2'1'} \otimes E_{22'} \\
\qquad{}- u^{-1}\bigl(uv^{-1} - 1\bigr)\bigl(uv^{-1} - \xi\bigr)E_{22'} \otimes E_{12} -r^{-2}s^{2}u^{-1}\bigl(uv^{-1} -1\bigr)\bigl(uv^{-1} -\xi\bigr)E_{2'2'} \otimes E_{12'} \\
\qquad{}+s^{2}u^{-1}\bigl(uv^{-1} - 1\bigr)\bigl(uv^{-1} - \xi\bigr)E_{12'} \otimes E_{11} + r^{-2}u^{-1}\bigl(uv^{-1} -1\bigr)\bigl(uv^{-1} - \xi\bigr)E_{1'2'} \otimes E_{11'} \\
\qquad{}-v^{-1}\bigl(uv^{-1} - 1\bigr)\bigl(uv^{-1} - \xi\bigr)E_{22} \otimes E_{21'} - r^{-2}s^{2}v^{-1}\bigl(uv^{-1} - 1\bigr)\bigl(uv^{-1} - \xi\bigr)E_{22'} \otimes E_{2'1'} \\
\qquad{}+v^{-1}\bigl(uv^{-1} - 1\bigr)\bigl(uv^{-1} - \xi\bigr)E_{11} \otimes E_{12'} + v^{-1}r^{-2}s^{2}\bigl(uv^{-1} - 1\bigr)\bigl(uv^{-1} - \xi\bigr)E_{11'} \otimes E_{1'2'},
\\
 F_{vu}R_{5}\bigl(uv^{-1}\bigr) - R_{4}\bigl(uv^{-1}\bigr)F_{uv} \\
 \quad{}= \bigl(r^{-2}s^{2} - 1\bigr)\bigl(uv^{-1} - \xi\bigr)v^{-1}E_{1'1'} \otimes E_{21'}+ \bigl(1 - r^{-2}s^{2}\bigr)\bigl(uv^{-1} - \xi\bigr)v^{-1}E_{22} \otimes E_{21'} \\
 \qquad{}+ \bigl(1 - r^{-2}s^{2}\bigr)\bigl(uv^{-1} - \xi\bigr)v^{-1}E_{2'2'} \otimes E_{21'}+ \bigl(1 - r^{-2}s^{2}\bigr)\bigl(uv^{-1} - \xi\bigr)v^{-1}E_{22} \otimes E_{12'}
 \\
 \qquad{}+\bigl(1 - r^{-2}s^{2}\bigr)\bigl(uv^{-1} - \xi\bigr)v^{-1}E_{2'2'} \otimes E_{12'}+ \bigl(r^{-2}s^{2} - 1\bigr)\bigl(uv^{-1} - \xi\bigr)v^{-1}E_{11} \otimes E_{12'}
 \\
 \qquad{}+\bigl(s^{2} - r^{2}\bigr)\bigl(uv^{-1} - \xi\bigr)u^{-1}E_{21'} \otimes E_{11}+ \bigl(r^{-4}s^{2} - r^{-2}\bigr)\bigl(uv^{-1} - \xi\bigr)uv^{-2}E_{12'} \otimes E_{1'1'},
\\
 F_{vu}R_{4}\bigl(uv^{-1}\bigr) - R_{5}\bigl(uv^{-1}\bigr)F_{uv} \\
 \quad{}= \bigl(r^{-2}s^{2} - 1\bigr)\bigl(uv^{-1} - \xi\bigr)uv^{-2}E_{11} \otimes E_{21'} + \bigl(r^{-4}s^{2} - r^{-2}\bigr)\bigl(uv^{-1} - \xi\bigr)v^{-1}E_{21'} \otimes E_{1'1'} \\
 \qquad{}+ \bigl(s^{-2} - r^{-2}\bigr)\bigl(uv^{-1} - \xi\bigr)v^{-1}E_{12'} \otimes E_{22}+\bigl(s^{2} - r^{-2}s^{4}\bigr)\bigl(uv^{-1} - \xi\bigr)v^{-1}E_{12'} \otimes E_{2'2'} \\
 \qquad{}+ \bigl(r^{-2}s^{2} - 1\bigr)\bigl(uv^{-1} - \xi\bigr)u^{-1}E_{1'1'} \otimes E_{12'}+\bigl(s^{-2} - r^{-2}\bigr)\bigl(uv^{-1} - \xi\bigr)v^{-1}E_{21'} \otimes E_{22} \\
 \qquad{}+ \bigl(s^{2} - r^{-2}s^{4}\bigr)\bigl(uv^{-1} - \xi\bigr)v^{-1}E_{21'} \otimes E_{2'2'}+\bigl(s^{2} - r^{2}\bigr)\bigl(uv^{-1} - \xi\bigr)v^{-1}E_{12'} \otimes E_{11},
\\
 R_{6}\bigl(uv^{-1}\bigr)F_{uv} \\
 \quad{}= \bigl(r^{-4}s^{2} - r^{-2}\bigr)\bigl(uv^{-1} - 1\bigr)v^{-1}E_{1'2'} \otimes E_{11'} - \bigl(r^{-2}s^{2} - 1\bigr)\bigl(uv^{-1} - 1\bigr)v^{-1}\xi E_{22'} \otimes E_{2'1'} \\
 \qquad{}-\bigl(r^{-2}s^{2} - 1\bigr)\bigl(uv^{-1} - \xi\bigr)v^{-1}E_{2'2'} \otimes E_{21'} - \bigl(r^{-2}s^{4} - s^{2}\bigr)\bigl(uv^{-1} - \xi\bigr)v^{-1}E_{21'} \otimes E_{2'2'} \\
 \qquad{}-\bigl(r^{-2}s^{4} - s^{2}\bigr)\bigl(uv^{-1} - 1\bigr)v^{-1}E_{2'1'} \otimes E_{22'} + \bigl(r^{-2}s^{2} - 1\bigr)\bigl(uv^{-1} - 1\bigr)\xi v^{-1}E_{11'} \otimes E_{1'2'} \\
 \qquad{}-\bigl(r^{-2}s^{2} - 1\bigr)\bigl(uv^{-1} - 1\bigr)\bigl(uv^{-1} - \xi\bigr)u^{-1}E_{1'1'} \otimes E_{12'}\\
 \qquad{} + \bigl(r^{-2}s^{2} - 1\bigr)\bigl(uv^{-1} - \xi\bigr)v^{-1}E_{1'1'} \otimes E_{12'} \\
 \qquad{} +\bigl(r^{-4}s^{2} - r^{-2}\bigr)\bigl(uv^{-1} - 1\bigr)\bigl(uv^{-1} - \xi\bigr)v^{-1}E_{12'} \otimes E_{1'1'} \\
 \qquad{} + \bigl(r^{-4}s^{2} - r^{-2}\bigr)\bigl(uv^{-1} - \xi\bigr)v^{-1}E_{12'} \otimes E_{1'1'},
\\
 F_{vu}R_{6}\bigl(uv^{-1}\bigr) \\
 \quad{}= \bigl(r^{2} - s^{2}\bigr)\bigl(uv^{-1} - \xi\bigr)u^{-1}E_{21'} \otimes E_{11} + \bigl(r^{-2}s^{2} - 1\bigr)\bigl(uv^{-1} - 1\bigr)v^{-1}E_{22'} \otimes E_{12}\\
 \qquad{} +\bigl(r^{-2}s^{2} - 1\bigr)\bigl(uv^{-1} - \xi\bigr)v^{-1}E_{22} \otimes E_{12'} + \bigl(r^{2} - s^{2}\bigr)\bigl(uv^{-1} - 1\bigr)\xi v^{-1}E_{21} \otimes E_{11'} \\
 \qquad{}+\bigl(r^{-2} - s^{-2}\bigr)\bigl(uv^{-1} - \xi\bigr)v^{-1}E_{12'} \otimes E_{22} + \bigl(r^{-2} - s^{-2}\bigr)\bigl(uv^{-1} - 1\bigr)\xi v^{-1}E_{12} \otimes E_{22'} \\
 \qquad{}-\bigl(r^{-2}s^{2} - 1\bigr)\bigl(uv^{-1} - 1\bigr)v^{-1}E_{11'} \otimes E_{21} - \bigl(r^{-2}s^{2} - 1\bigr)\bigl(uv^{-1} - \xi\bigr)uv^{-2}E_{11} \otimes E_{21'}.
\end{gather*}
Combining the expressions above, we get
\begin{gather}
 R_{3}\bigl(uv^{-1}\bigr)F_{uv} - F_{vu}R_{3}\bigl(uv^{-1}\bigr) + R_{1}\bigl(uv^{-1}\bigr)F_{uv}\label{eq:R123_affine}\\
\quad{} + R_{2}\bigl(uv^{-1}\bigr)F_{uv} - F_{vu}R_{1}\bigl(uv^{-1}\bigr) - F_{vu}R_{2}\bigl(uv^{-1}\bigr) \nonumber\\
 \quad{} = \bigl(1 - r^{-2}s^{2}\bigr)\bigl(uv^{-1} - \xi\bigr)v^{-1}E_{22} \otimes E_{21'} + \bigl(r^{-2}s^{2} - 1\bigr)\bigl(uv^{-1} - \xi\bigr)v^{-1}E_{11} \otimes E_{12'} \nonumber\\
 \qquad{} + \bigl(s^{-2} - r^{-2}\bigr)\bigl(uv^{-1} - \xi\bigr)v^{-1}E_{21'} \otimes E_{22} + \bigl(s^{2} - r^{2}\bigr)\bigl(uv^{-1} - \xi\bigr)v^{-1}E_{12'} \otimes E_{11} \nonumber\\
 \qquad{} + \bigl(r^{-2}s^{2} - 1\bigr)\bigl(uv^{-1} - 1\bigr)\xi v^{-1}E_{22'} \otimes E_{2'1'} + \bigl(1 - r^{-2}s^{2}\bigr)\bigl(uv^{-1} - 1\bigr)\xi v^{-1}E_{11'} \otimes E_{1'2'} \nonumber\\
 \qquad{} + \bigl(r^{-2}s^{4} - s^{2}\bigr)\bigl(uv^{-1} - 1\bigr)v^{-1}E_{2'1'} \otimes E_{22'} + \bigl(r^{-2} - r^{-4}s^{2}\bigr)\bigl(uv^{-1} - 1\bigr)v^{-1}E_{1'2'} \otimes E_{11'} \nonumber\\
 \qquad{} + \bigl(r^{-2}s^{2} - 1\bigr)\bigl(uv^{-1} - \xi\bigr)v^{-1} E_{1'1'} \otimes E_{21'} + \bigl(1 - r^{-2}s^{2}\bigr)\bigl(uv^{-1} - \xi\bigr)v^{-1}E_{2'2'} \otimes E_{12'} \nonumber\\
 \qquad{} + \bigl(r^{-4}s^{2} - r^{-2}\bigr)\bigl(uv^{-1} - \xi\bigr)v^{-1}E_{21'} \otimes E_{1'1'}+ \bigl(s^{2} - r^{-2}s^{4}\bigr)\bigl(uv^{-1} - \xi\bigr)v^{-1}E_{12'} \otimes E_{2'2'} \nonumber\\
 \qquad{} + \bigl(1 - r^{-2}s^{2}\bigr)\bigl(uv^{-1} - 1\bigr)v^{-1}E_{11'} \otimes E_{21} +\bigl(r^{-2}s^{2} - 1\bigr)\bigl(uv^{-1} - 1\bigr)v^{-1}E_{22'} \otimes E_{12} \nonumber\\
 \qquad{} + \bigl(r^{2} - s^{2}\bigr)\bigl(uv^{-1} - 1\bigr)\xi v^{-1}E_{21} \otimes E_{11'} +\bigl(r^{-2} - s^{-2}\bigr)\bigl(uv^{-1} - 1\bigr)\xi v^{-1}E_{12} \otimes E_{22'}\nonumber
\end{gather}
and
\begin{gather}
 \bigl(F_{vu}R_{5}\bigl(uv^{-1}\bigr) - R_{4}\bigl(uv^{-1}\bigr)F_{uv}\bigr) + \bigl(F_{vu}R_{4}\bigl(uv^{-1}\bigr) - R_{5}\bigl(uv^{-1}\bigr)F_{uv}\bigr) \label{eq:R456_affine}\\
 \quad {}+ F_{vu}R_{6}\bigl(uv^{-1}\bigr) - R_{6}\bigl(uv^{-1}\bigr)F_{uv} \nonumber\\
 \quad{}= \bigl(r^{-2}s^{2} - 1\bigr)\bigl(uv^{-1} - \xi\bigr)v^{-1}E_{1'1'} \otimes E_{21'} + \bigl(1 - r^{-2}s^{2}\bigr)\bigl(uv^{-1} - \xi\bigr)v^{-1}E_{22} \otimes E_{21'} \nonumber\\
 \qquad{} + \bigl(1 - r^{-2}s^{2}\bigr)\bigl(uv^{-1} - \xi\bigr)v^{-1}E_{2'2'} \otimes E_{12'} + \bigl(r^{-2}s^{2} - 1\bigr)\bigl(uv^{-1} - \xi\bigr)v^{-1}E_{11} \otimes E_{12'} \nonumber\\
 \qquad{} + \bigl(r^{-4}s^{2} - r^{-2}\bigr)\bigl(uv^{-1} -\xi\bigr)v^{-1}E_{21'} \otimes E_{1'1'} + \bigl(s^{2} - r^{-2}s^{4}\bigr)\bigl(uv^{-1} - \xi\bigr)v^{-1}E_{12'} \otimes E_{2'2'} \nonumber\\
 \qquad{} + \bigl(s^{-2} - r^{-2}\bigr)\bigl(uv^{-1} - \xi\bigr)v^{-1}E_{21'} \otimes E_{22} + \bigl(s^{2} - r^{2}\bigr)\bigl(uv^{-1} - \xi\bigr)v^{-1}E_{12'} \otimes E_{11} \nonumber\\
 \qquad{} +\bigl(r^{-2}s^{2} - 1\bigr)\bigl(uv^{-1} - 1\bigr)v^{-1}E_{22'} \otimes E_{12} +\bigl(r^{2} - s^{2}\bigr)\bigl(uv^{-1} - 1\bigr)\xi v^{-1}E_{21} \otimes E_{11'} \nonumber\\
 \qquad{} + \bigl(r^{-2} - s^{-2}\bigr)\bigl(uv^{-1} - 1\bigr)\xi v^{-1}E_{12} \otimes E_{22'} + \bigl(1 - r^{-2}s^{2}\bigr)\bigl(uv^{-1} - 1\bigr)v^{-1}E_{11'} \otimes E_{21} \nonumber\\
 \qquad{} +\bigl(r^{-2} - r^{-4}s^{2}\bigr)\bigl(uv^{-1} - 1\bigr)v^{-1}E_{1'2'} \otimes E_{11'} + \bigl(r^{-2}s^{2} - 1\bigr)\bigl(uv^{-1} - 1\bigr)\xi v^{-1}E_{22'} \otimes E_{2'1'} \nonumber\\
 \qquad{} + \bigl(r^{-2}s^{4} - s^{2}\bigr)\bigl(uv^{-1} - 1\bigr)v^{-1}E_{2'1'} \otimes E_{22'} + \bigl(1 - r^{-2}s^{2}\bigr)\bigl(uv^{-1} - 1\bigr)\xi v^{-1}E_{11'} \otimes E_{1'2'}.\nonumber
\end{gather}
The right-hand sides of \eqref{eq:R123_affine} and \eqref{eq:R456_affine} are clearly equal, which completes the proof of~\eqref{eq:affine-intertwiner} for $x = f_{0}$.

In the same way, one can verify \eqref{eq:affine-intertwiner} for $x = e_{0}$ (although according to~\cite{Jim}, it is actually sufficient to check~\eqref{eq:affine-intertwiner} only for $x = f_{i}$). This completes the proof of part (a) of Theorem~\ref{thm:B_AffRMatrix}. The proof of part (b) actually follows from~(a), as in~\cite[Proposition~3]{Jim}.
\end{proof}

\begin{Remark}
Similarly to Remark~\ref{rem:easy-finite}, we note that the above proofs of Theorems~\ref{thm:B_AffRMatrix}--\ref{thm:D_AffRMatrix} are quite elementary, but they require knowing the correct formulas for $\hat{R}(z)$ at the first place. In the next section, we present the origin of these formulas, by using the Yang--Baxterization technique of~\cite{GWX}.
\end{Remark}

\section{Yang--Baxterization}\label{sec:Baxterization}

In this last section, we present a natural derivation of the rather complicated formulas~\eqref{eq:A_AffRMatrix}--\eqref{eq:D_AffRMatrix} from their finite counterparts~\eqref{eq:A_RMatrix},~\eqref{eq:B_RMatrix},~\eqref{eq:C_RMatrix} and~\eqref{eq:D_RMatrix}. This is based on a so-called \emph{Yang--Baxterization} technique of~\cite{GWX}, which produces $\hat{R}(z)$ satisfying~\eqref{eq:qYB-affine} from $\hat{R}$ satisfying~\eqref{eq:qYB-two} when the latter has $2$ or $3$ eigenvalues.

\textbf{Yang--Baxterization in $\boldsymbol{A}$-type.}
For a uniform exposition, we start by recalling the derivation of~\eqref{eq:A_AffRMatrix} via that technique. As noted in~\cite{BW2}, the $R$-matrix $\hat{R}$ of~\eqref{eq:A_RMatrix} is diagonalizable with two eigenvalues $\lambda_1=-rs^{-1}$, $\lambda_2=1$, so that
\begin{gather}
 \hat{R}^{-1}=-\lambda_1^{-1}\lambda_2^{-1}\hat{R}+\bigl(\lambda^{-1}_1+\lambda^{-1}_2\bigr)\Id.\label{eq:R_inverse_A}
\end{gather}
In that setup, the Yang--Baxterization of~\cite[equation~(3.15)]{GWX} produces the following solution of~\eqref{eq:qYB-affine}:
\begin{equation}\label{eq:Baxterization-simple}
 \hat{R}(z)=\lambda_2^{-1}\hat{R} + z\lambda_1 \hat{R}^{-1}.
\end{equation}
Plugging explicit formulas for $\hat{R},\hat{R}^{-1}$ from~\eqref{eq:A_RMatrix} and~\eqref{eq:R_inverse_A} into the right-hand side of~\eqref{eq:Baxterization-simple},
we~derive precisely the operator~\eqref{eq:A_AffRMatrix}.

\textbf{Yang--Baxterization in $\boldsymbol{BCD}$-types.}
Let us treat the other three classical series. To this end, we recall that the $R$-matrices $\hat{R}$ of~\eqref{eq:B_RMatrix},~\eqref{eq:C_RMatrix} and~\eqref{eq:D_RMatrix} have three distinct eigenvalues~$\lambda_1$, $\lambda_2$, $\lambda_3$, in accordance with Lemmas~\ref{lem:B_eigen},~\ref{lem:C_eigen} and~\ref{lem:D_eigen}. In that setup, the Yang--Baxterization of~\cite[equations~(3.29) and~(3.31)]{GWX} produces the following two solutions to~\eqref{eq:qYB-affine}:
\begin{gather}\label{eq:affinization_a}
 \hat{R}(z) =
 \lambda_{1}z(z - 1)\hat{R}^{-1} +
 \biggl(1 + \frac{\lambda_{1}}{\lambda_{2}} + \frac{\lambda_{1}}{\lambda_{3}} + \frac{\lambda_{2}}{\lambda_{3}}\biggr) z \,\Id
 - \frac{1}{\lambda_{3}}(z - 1)\hat{R},
\\ \label{eq:affinization_b}
 \hat{R}(z) = \lambda_{1}z(z - 1)\hat{R}^{-1} + \biggl(1 + \frac{\lambda_{1}}{\lambda_{2}} + \frac{\lambda_{1}}{\lambda_{3}} + \frac{\lambda_{1}^{2}}{\lambda_{2}\lambda_{3}}\biggr)z\,\Id - \frac{\lambda_{1}}{\lambda_{2}\lambda_{3}}(z - 1)\hat{R},
\end{gather}
provided that $\hat{R}$ satisfies the additional relations of~\cite[equation~(3.27)]{GWX} (cf.\ correction~\cite[equation~(A.9)]{GWX}), which, in particular, hold whenever $\hat{R}$ is a representation of a \emph{Birman--Wenzl algebra}.

\begin{Remark}
For the purpose of the present section, we shall not really need to verify these additional relations, since according to Theorems~\ref{thm:B_AffRMatrix}--\ref{thm:D_AffRMatrix} the constructed $\hat{R}(z)$ do manifestly satisfy the relation~\eqref{eq:qYB-affine}.
\end{Remark}

To apply formulas above, it remains to evaluate $\hat{R}^{-1}$. To this end, we consider the $\BC$-algebra involution
\begin{gather}\label{eq:sigma-antiautom}
 \sigma\colon\ \uu \to \uu
\end{gather}
given by
\[
 e_i\mapsto e_i, \qquad f_i\mapsto f_i, \qquad \omega_{i}\mapsto \omega_{i}', \qquad
 \omega_{i}' \mapsto \omega_{i}, \qquad r\mapsto s, \qquad s\mapsto r.
\]
Evoking the notation~\eqref{eq:Theta}, we define
\begin{displaymath}
 \overline{\Theta} = \sum_{\mu\ge 0} \overline{\Theta}_{\mu} \qquad \mathrm{with} \
 \overline{\Theta}_{\mu} = (\sigma \otimes \sigma) (\Theta_{\mu}) \ \mathrm{for\ all} \ \mu.
\end{displaymath}
We also introduce another coproduct homomorphism $\overline{\Delta}\colon \uu \to \uu \otimes \uu$ via
\begin{displaymath}
 \overline{\Delta} = (\sigma \otimes \sigma) \circ \Delta \circ \sigma^{-1}.
\end{displaymath}
Then, for any finite-dimensional $\uu$-modules $V$ and $W$, we have
\begin{displaymath}
 \Delta(u) \circ \Theta = \Theta \circ \overline{\Delta}(u)\colon\ V \otimes W \to V \otimes W
 \qquad \mathrm{for\ all} \ u \in \uu,
\end{displaymath}
cf.\ \cite[Lemma~4.10]{BW1} and~\cite[Lemma~3.3]{BGH2}. Applying $\sigma \otimes \sigma$ to the equality above, we then get $\overline{\Delta}(\sigma(u)) \overline{\Theta} = \overline{\Theta} \Delta(\sigma(u))$. Since $\sigma$ is an automorphism, the last equality can be written as $\overline{\Delta}(u) \overline{\Theta} = \overline{\Theta}\Delta(u)$ for any $u\in \uu$. Let us now also show that $\overline{\Delta}(u) \widetilde{f} = \widetilde{f}\Delta^{\op}(u)$ on $V \otimes W$. It suffices to verify this formula when $u$ is one of the generators. For $u = \omega_{i}$ or $\omega_{i}'$, this is obvious. For $u = e_{i}$ and any $v \in V[\lambda]$, $w \in V[\mu]$:
\begin{displaymath}
 \overline{\Delta}(e_{i}) \widetilde{f} (v \otimes w) =
 f(\lambda,\mu) \bigl(e_{i}v \otimes w + \bigl(\omega_{i}',\omega_{\lambda}\bigr)^{-1} v \otimes e_{i}w\bigr)
\end{displaymath}
and
\begin{displaymath}
 \widetilde{f}(\Delta^{\op}(e_{i})(v \otimes w)) =
 f(\lambda,\mu + \alpha_{i})v \otimes e_{i}w + f(\lambda + \alpha_{i},\mu)\bigl(\omega_{\mu}',\omega_{i}\bigr) e_{i}v \otimes w,
\end{displaymath}
which are equal due to the properties~\eqref{eq:f} satisfied by $f$. The proof for $u = f_{i}$ is completely analogous. Putting all of this together, we find
\begin{displaymath}
 \tau \circ \widetilde{f}^{-1} \circ \overline{\Theta} \circ \Delta(u) =
 \tau \circ \widetilde{f}^{-1} \circ \overline{\Delta}(u)\circ \overline{\Theta} =
 \tau \circ \Delta^{\op}(u) \circ \widetilde{f}^{-1}\circ \overline{\Theta} =
 \Delta(u) \circ \tau \circ \widetilde{f}^{-1} \circ \overline{\Theta}
\end{displaymath}
as linear maps $V \otimes W \to W \otimes V$. Thus $\overline{R} = \tau \circ \widetilde{f}^{-1}\circ \overline{\Theta}\colon V \otimes W \to W \otimes V$ is a $\uu$-module isomorphism.

Specializing now to the case where $V = W$ is one of the representations from Propositions~\ref{prp:B_rep}--\ref{prp:D_rep}, one can easily see from the defining formulas that $\rho(\sigma(u)) = \widetilde{\sigma}(\rho(u))$ for any $u \in \uu$. Here, we regard $\rho(u)$ as an element of $\mathrm{Mat}_{N}(\mathbb{K})$, and $\widetilde{\sigma}\colon \mathrm{Mat}_{N}(\mathbb{K}) \to \mathrm{Mat}_{N}(\mathbb{K})$ is the $\mathbb{C}$-algebra automorphism defined by
\begin{displaymath}
 \widetilde{\sigma}(r^{k}s^{\ell}E_{ij}) = r^{\ell}s^{k}E_{ij} \qquad \forall k,\ell\in \BZ,\ 1\leq i,j\leq N.
\end{displaymath}
By abuse of notation, we shall use $\widetilde{\sigma}$ to denote similar $\BC$-algebra automorphisms $\mathrm{Mat}_{N}(\mathbb{K})^{\otimes 2} \to \mathrm{Mat}_{N}(\mathbb{K})^{\otimes 2}$ and $\mathbb{K}\to \mathbb{K}$. To obtain explicit formulas for $\overline{R}$, which we present in equations \eqref{eq:B_Inv}--\eqref{eq:D_Inv} below, we just need to evaluate \smash{$\tau \circ \widetilde{f}^{-1} \circ \widetilde{\sigma}\bigl(\hat{R} \circ \tau \circ \widetilde{f}^{-1}\bigr)$} from the respective formulas for \smash{$\hat{R}$} given in~\eqref{eq:B_RMatrix},~\eqref{eq:C_RMatrix} and~\eqref{eq:D_RMatrix}, see~\eqref{eq:R-int}. Most of the terms transform easily into the corresponding terms in \eqref{eq:B_Inv}--\eqref{eq:D_Inv}, due to the equalities
\begin{displaymath}
 \widetilde{\sigma}(f(\varepsilon_j,\varepsilon_i))=f(\varepsilon_i,\varepsilon_j)^{-1}
 \qquad \mathrm{for\ all} \ i,j,
\end{displaymath}
but some additional explanation is necessary for two of them. First, the term $\sum_{i \neq j,j'} a_{ij}E_{ij} \otimes E_{ji}$ transforms into the corresponding one in~\eqref{eq:B_Inv}--\eqref{eq:D_Inv} due to the following identities satisfied by~$a_{ij}$ of~\eqref{eq:B_coeffs},~\eqref{eq:C_coeffs} and~\eqref{eq:D_coeffs}:
\begin{displaymath}
 a_{ij} = f(\varepsilon_{i},\varepsilon_{j}) \qquad \mathrm{and} \qquad a_{ji}^{-1} = a_{ij}
 \qquad \mathrm{for\ all} \ i \neq j,j'.
\end{displaymath}
For the former equality, see~\eqref{eq:f_aij_B} and~\eqref{eq:f_aij_CD}.

Second, the term \smash{$\sum_{i < j \neq i'} t_{i}t_{j}^{-1}E_{i'j} \otimes E_{ij'}$} transforms into the corresponding one in \eqref{eq:B_Inv}--\eqref{eq:D_Inv} due to the additional observations that, unless $i = n + 1$ or $j = n + 1$ in type $B_{n}$, we have (with $\varepsilon_{i'}=-\varepsilon_{i}$ as defined prior to~\eqref{eq:f_aij_B} and~\eqref{eq:f_aij_CD})
\begin{displaymath}
 f(-\varepsilon_{i},\varepsilon_{i}) = f(\varepsilon_{j},-\varepsilon_{j})
 \qquad \text{and} \qquad \widetilde{\sigma}\bigl(t_{i'}t_{j'}^{-1}\bigr) = t_{i}t_{j}^{-1}
 \qquad \mathrm{for\ all} \ 1\leq i,j\leq N.
\end{displaymath}
The former equality follows from our formulas~\eqref{eq:B-f-1} and~\eqref{eq:B-f-2} for type $B_n$ and \eqref{eq:C-wght-pair-1} and~\eqref{eq:C-wght-pair-2} for types $C_n$, $D_n$. When $i=n+1$ or $j=n+1$ in type $B_n$, we rather use the following equalities:
\begin{displaymath}
 f(0,0) = 1,\qquad \widetilde{\sigma}(t_{n+1})rs^{-1} = t_{n+1},\qquad
 \widetilde{\sigma}(t_{i'}) = t_{i} \qquad \mathrm{for} \ i \neq n + 1.
\end{displaymath}

These results allow us to prove the following lemmas.

\begin{Lemma}[type $B_{n}$]\label{lem:B_Inverse}
The inverse of the operator $\hat{R}\colon V \otimes V \to V \otimes V$ from Theorem $\ref{thm:B_RMatrix}$ is equal to
\begin{align}
 \overline{R} &= \tau \circ \widetilde{f}^{-1} \circ \overline{\Theta} \nonumber\\
 &= rs^{-1}\sum_{1\leq i\leq 2n+1}^{i\ne n+1} E_{ii} \otimes E_{ii} + E_{n+1,n+1} \otimes E_{n+1,n+1} +
 r^{-1}s \sum_{1\leq i\leq 2n+1}^{i\ne n+1} E_{ii'} \otimes E_{i'i} \nonumber\\
 &\quad + \sum_{1\leq i,j\leq 2n+1}^{j\ne i,i'} a_{ij} E_{ij} \otimes E_{ji}
 + \bigl(s^{2} - r^{2}\bigr)(rs)^{-1}\sum_{i = 1}^{n}\bigl(r^{2(i - n) - 1}s^{2(n-i) + 1} - 1\bigr) E_{ii} \otimes E_{i'i'} \nonumber\\
 &\quad + \bigl(r^{2} - s^{2}\bigr)(rs)^{-1} \sum_{i<j}^{j\ne i'} E_{ii} \otimes E_{jj}
 + \bigl(s^{2} - r^{2}\bigr)(rs)^{-1}\sum_{i > j}^{j \ne i'} t_{i}t_{j}^{-1}E_{i'j} \otimes E_{ij'},\label{eq:B_Inv}
\end{align}
with the constants $t_{i}$ and $a_{ij}$ given explicitly by~\eqref{eq:B_coeffs}.
\end{Lemma}

\begin{Lemma}[type $C_{n}$]
The inverse of the operator $\hat{R}\colon V \otimes V \to V \otimes V$ from Theorem $\ref{thm:C_RMatrix}$ is equal to
\begin{align}
 \overline{R} &= \tau \circ \widetilde{f}^{-1} \circ \overline{\Theta} \nonumber\\
 &=
 r^{1/2}s^{-1/2}\sum_{i = 1}^{2n} E_{ii} \otimes E_{ii}
 + r^{-1/2}s^{1/2}\sum_{i = 1}^{2n} E_{ii'} \otimes E_{i'i} \nonumber\\
 &\quad + \sum_{1\leq i,j\leq 2n}^{j\ne i,i'} a_{ij}E_{ij} \otimes E_{ji}
 + (r - s)(rs)^{-1/2}\sum_{i = 1}^{n}\bigl(r^{i - n- 1}s^{n + 1 - i} + 1\bigr)E_{ii} \otimes E_{i'i'}\nonumber \\
 &\quad + (r -s)(rs)^{-1/2}\sum_{i < j}^{j\neq i'} E_{ii} \otimes E_{jj}
 + (s - r)(rs)^{-1/2}\sum_{i > j}^{j\neq i'} t_{i}t_{j}^{-1}E_{i'j} \otimes E_{ij'},\label{eq:C_Inv}
\end{align}
with the constants $t_{i}$ and $a_{ij}$ given explicitly by~\eqref{eq:C_coeffs}.
\end{Lemma}

\begin{Lemma}[type $D_{n}$]\label{lem:D_Inverse}
The inverse of the operator $\hat{R}\colon V \otimes V \to V \otimes V$ from Theorem $\ref{thm:D_RMatrix}$ is equal to
\begin{align}
 \overline{R} &= \tau \circ \widetilde{f}^{-1} \circ \overline{\Theta} \nonumber\\
 &=r^{1/2}s^{-1/2}\sum_{i = 1}^{2n} E_{ii} \otimes E_{ii}
 + r^{-1/2}s^{1/2}\sum_{i = 1}^{2n} E_{ii'} \otimes E_{i'i} \nonumber\\
 &\quad + \sum_{1\leq i,j\leq 2n}^{j \neq i,i'} a_{ij}E_{ij} \otimes E_{ji}
 + (r - s)(rs)^{-1/2}\sum_{i = 1}^{n}\bigl(1 - r^{i- n}s^{n - i}\bigr) E_{ii} \otimes E_{i'i'} \nonumber\\
 &\quad + (r -s)(rs)^{-1/2}\sum_{i < j}^{j \neq i'} E_{ii} \otimes E_{jj}
 + (s - r)(rs)^{-1/2}\sum_{i > j}^{j \neq i'} t_{i}t_{j}^{-1}E_{i'j} \otimes E_{ij'},\label{eq:D_Inv}
\end{align}
with the constants $t_{i}$ and $a_{ij}$ given explicitly by~\eqref{eq:D_coeffs}.
\end{Lemma}

Since $\overline{R}$ is a $\uu$-module intertwiner by the above discussions, it suffices to verify that the eigenvalues of $\overline{R}$ on the highest weight vectors $w_1$, $w_2$, $w_3$ from our proof of Proposition~\ref{prop:struct} are inverse to those of $\hat{R}$ as specified in Lemmas~\ref{lem:B_eigen}--\ref{lem:D_eigen}. As the arguments are very similar, we~shall only present the proof in type $B_{n}$.

\begin{proof}[Proof of Lemma \ref{lem:B_Inverse}]
For $w_{1} = v_{1} \otimes v_{1}$, the eigenvalue of $\hat{R}$ is $\lambda_1=r^{-1}s$, while we clearly have
\begin{displaymath}
 \overline{R}(w_1)=\overline{R}(v_{1} \otimes v_{1}) = rs^{-1} v_{1} \otimes v_{1}=\lambda^{-1}_1 w_1 .
\end{displaymath}
For $w_{2}$, we have
\begin{displaymath}
\begin{split}
 \overline{R}(w_{2}) = \overline{R}\bigl(v_{1} \otimes v_{2} - rs^{-1}v_{2} \otimes v_{1}\bigr) &=
 \bigl(rs^{-1} - r^{-1}s\bigr)v_{1} \otimes v_{2} + v_{2} \otimes v_{1} - rs^{-1}v_{1} \otimes v_{2} \\
 &=
 -r^{-1}s\bigl(v_{1} \otimes v_{2} - rs^{-1}v_{2} \otimes v_{1}\bigr)
\end{split}
\end{displaymath}
for $n = 1$, and
\begin{displaymath}
\begin{split}
 \overline{R}(w_{2}) = \overline{R}\bigl(v_{1} \otimes v_{2} - r^{2}v_{2} \otimes v_{1}\bigr) &=
 \bigl(rs^{-1} - r^{-1}s\bigr)v_{1} \otimes v_{2} + rsv_{2} \otimes v_{1} - rs^{-1}v_{1} \otimes v_{2} \\
 &=
 -r^{-1}s\bigl(v_{1} \otimes v_{2} - r^{2}v_{2} \otimes v_{1}\bigr)
\end{split}
\end{displaymath}
for $n > 1$.  Thus, we obtain $\overline{R}(w_{2})=\lambda_2^{-1} w_2$ for all~$n$.

Finally, the eigenvalue of the $\overline{R}$-action on $w_3$ equals the ratio of coefficients of $v_{1'}\otimes v_{1}$ in~$\overline{R}(w_3)$ and $w_3$. To compute the former, we note that only the third summand of~\eqref{eq:B_Inv} makes a~nontrivial contribution of $r^{-1}s\cdot v_{1'}\otimes v_{1}$. As the coefficient of $v_{1'}\otimes v_{1}$ in $w_3$ equals $r^{2n - 1}s^{-2n + 1}$, we conclude that \smash{$\overline{R}(w_{3}) = r^{-2n}s^{2n}w_{3}$}. Thus, the eigenvalue of $w_3$ for $\overline{R}$ is indeed equal to \smash{$\lambda_3^{-1}=r^{-2n}s^{2n}$}.
\end{proof}

With these preliminaries out of the way, we can now present our formulas for $\hat{R}(z)$, which we obtain from \eqref{eq:affinization_a} in types $B_{n}$ and $D_{n}$, and from \eqref{eq:affinization_b} in type $C_{n}$ (the key reason to use a~different formula in type $C_n$ is because \eqref{eq:affinization_a} does not produce a solution that also satisfies the intertwining property \eqref{eq:affine-intertwiner}):
\begin{itemize}\itemsep=0pt
\item
\emph{Type $B_n$}:
\begin{displaymath}
 \hat{R}(z) =
 r^{-1}sz(z - 1)\overline{R} + \bigl(1 - r^{-2n + 1}s^{2n - 1}\bigr)\bigl(1 - r^{-2}s^{2}\bigr)z\,\Id - r^{-2n}s^{2n}(z - 1)\hat{R},
\end{displaymath}
which after a direct computation simplifies to~\eqref{eq:B_AffRMatrix}.

\item
\emph{Type $C_n$}:
\begin{gather*}
 \hat{R}(z) =
 r^{-1/2}s^{1/2}z(z-1)\overline{R} + \bigl(1 - r^{-n-1}s^{n + 1}\bigr)\bigl(1 - r^{-1}s\bigr)z\,\Id \\
 \hphantom{\hat{R}(z) =}{}
 - r^{-n - 3/2}s^{n+ 3/2}(z - 1)\hat{R},
\end{gather*}
which after a direct computation simplifies to~\eqref{eq:C_AffRMatrix}.

\item
\emph{Type $D_n$}:
\begin{gather*}
 \hat{R}(z) =
 r^{-1/2}s^{1/2}z(z-1)\overline{R} + \bigl(1 - r^{-n+1}s^{n-1}\bigr)\bigl(1 - r^{-1}s\bigr)z\,\Id \\
 \hphantom{\hat{R}(z) =}{}
 - r^{-n+1/2}s^{n-1/2}(z - 1)\hat{R},
\end{gather*}
which after a direct computation simplifies to~\eqref{eq:D_AffRMatrix}.
\end{itemize}

\subsection*{Acknowledgements}

This note represents a part (followed up by~\cite{MT1,MT2}) of the year-long REU project at Purdue University. We are grateful to Purdue University for support and for the opportunity to present these results at REU math conference in April~2024.

A.T.\ is deeply indebted to Andrei Negu\c{t} for numerous stimulating discussions over the years and sharing the beautiful combinatorial features of quantum groups, to Sarah Witherspoon for a~correspondence on two-parameter quantum groups, to Curtis Wendlandt for bringing attention to~\cite{HM}, to Rouven Fraseek and Daniele Valeri for invitations for research visits in Italy during the summer of 2024.
A.T.\ is grateful to INdAM-GNSAGA and FAR UNIMORE project CUP-E93C2300204000 for the support and wonderful working conditions during his visit to Italy, where the final version of the paper was completed; to IHES for the hospitality and wonderful working conditions in the summer 2025, where the journal version of this paper was prepared. We are very grateful to the referees for their useful suggestions that improved the exposition.

The work of both authors was partially supported by an NSF Grant DMS-$2302661$.

\pdfbookmark[1]{References}{ref}
\LastPageEnding

\end{document}